\documentclass[11pt,leqno]{amsart}
\usepackage{amsthm,amsfonts,amssymb,amsmath,oldgerm,mathabx}
\numberwithin{equation}{section}
\usepackage{fullpage}
\usepackage{color}

\usepackage[pdftex]{graphicx} 
\usepackage{hyperref}
\usepackage{stmaryrd}
\usepackage{comment}

\usepackage{caption}
\usepackage{subcaption}

\renewcommand\d{\partial}
\DeclareMathOperator\eD{e}
\DeclareMathOperator\iD{i}
\newcommand{\mypar}{\shortparallel}
\def\eps{\varepsilon }
\newcommand{\eDr}{\eD_{\rm r}}

\DeclareMathOperator{\Id}{Id}
\DeclareMathOperator{\I}{I}
\DeclareMathOperator{\b0}{\mathbf{0}}
\newcommand{\transp}[1]{{#1}^{{\sf T}}}
\newcommand{\cond}[1]{{\bf (D#1)}}
\newcommand{\case}[1]{{\bf (C#1)}}
\DeclareMathOperator{\Span}{Span}

\DeclareMathOperator{\diag}{diag}

\DeclareMathOperator{\dD}{d}
\DeclareMathOperator{\Div}{\transp{\nabla}}
\DeclareMathOperator{\Curl}{\nabla\wedge}
\DeclareMathOperator{\beD}{\mathbf{e}}
\DeclareMathOperator{\supp}{supp}

\newcommand\br{\begin{remark}}
\newcommand\er{\end{remark}}
\newcommand\bp{\begin{pmatrix}}
\newcommand\ep{\end{pmatrix}}
\newcommand{\be}{\begin{equation}}
\newcommand{\ee}{\end{equation}}
\newcommand{\bes}{\begin{equation*}}
\newcommand{\ees}{\end{equation*}}
\newcommand\ba{\begin{equation}\begin{aligned}}
\newcommand\ea{\end{aligned}\end{equation}}
\newcommand\bas{\begin{equation*}\begin{aligned}}
\newcommand\eas{\end{aligned}\end{equation*}}
\newcommand\ds{\displaystyle}
\newcommand\nn{\nonumber}
\newcommand{\beg}{\begin{example}}
\newcommand{\eeg}{\end{example}}
\newcommand{\bpr}{\begin{proposition}}
\newcommand{\epr}{\end{proposition}}
\newcommand{\bt}{\begin{theorem}}
\newcommand{\et}{\end{theorem}}
\newcommand{\bc}{\begin{corollary}}
\newcommand{\ec}{\end{corollary}}
\newcommand{\bl}{\begin{lemma}}
\newcommand{\el}{\end{lemma}}
\newcommand{\bd}{\begin{definition}}
\newcommand{\ed}{\end{definition}}
\newcommand{\brs}{\begin{remarks}}
\newcommand{\ers}{\end{remarks}}

\newtheorem{theorem}{Theorem}[section]
\newtheorem{proposition}[theorem]{Proposition}
\newtheorem{corollary}[theorem]{Corollary}
\newtheorem{lemma}[theorem]{Lemma}

\theoremstyle{remark}
\newtheorem{remark}[theorem]{Remark}
\theoremstyle{definition}
\newtheorem{definition}[theorem]{Definition}

\newtheorem{example}[theorem]{Example}

\newcommand\R{\mathbf R}
\newcommand\C{\mathbf C}
\newcommand{\N}{\mathbf N}
\newcommand{\Z}{\mathbf Z}

\newcommand{\lnor}{|\!|\!|}
\newcommand{\rnor}{|\!|\!|}

\newcommand{\bcKW}{\bcK^W}
\newcommand{\bPsiW}{\bPsi^W}
\newcommand{\bphiW}{\bfphi^W}
\newcommand{\Wbphi}{\bfphi_W}
\newcommand{\brW}{\bfr^W}
\newcommand{\bPsiaux}{\bPsi^{\rm aux}}
\newcommand{\bphiaux}{\bfphi^{\rm aux}}
\newcommand{\auxbphi}{\bfphi_{\rm aux}}
\newcommand{\braux}{\bfr^{\rm aux}}
\newcommand{\bDW}{\bD^W}
\newcommand{\bDO}{\bD^{(0)}}
\newcommand{\lambO}{\lambda^{(0)}}
\newcommand{\piO}{\pi^{(0)}}
\newcommand{\gin}{\bfg^{\textrm{in}}}
\newcommand{\Vin}{\bV^{\textrm{in}}}
\newcommand{\Vout}{\bV^{\textrm{out}}}
\newcommand{\phiin}{\bfphi^{\textrm{in}}}
\newcommand{\phiout}{\bfphi^{\textrm{out}}}
\newcommand{\Smod}{S_{\textrm{mod}}}
\newcommand{\SigW}{\Sigma_W}
\newcommand{\SigLF}{\Sigma_{\rm LF}}
\newcommand{\SigOLF}{\Sigma^{\rm LF}_{(0)}}
\newcommand{\SigOHF}{\Sigma^{\rm HF}_{(0)}}
\newcommand{\SigO}{\Sigma_{(0)}}

\newcommand{\CC}{{\mathbb C}}

\newcommand\bA{{\mathbf A}}
\newcommand\bB{{\mathbf B}}

\newcommand\bD{{\mathbf D}}
\newcommand\bE{{\mathbf E}}
\newcommand\bF{{\mathbf F}}
\newcommand\bG{{\mathbf G}}
\newcommand\bH{{\mathbf H}}

\newcommand\bK{{\mathbf K}}
\newcommand\bL{{\mathbf L}}
\newcommand\bM{{\mathbf M}}

\newcommand\bP{{\mathbf P}}
\newcommand\bQ{{\mathbf Q}}
\newcommand\bR{{\mathbf R}}
\newcommand\bS{{\mathbf S}}

\newcommand\bU{{\mathbf U}}
\newcommand\bV{{\mathbf V}}
\newcommand\bW{{\mathbf W}}
\newcommand\bX{{\mathbf X}}
\newcommand\bY{{\mathbf Y}}
\newcommand\bZ{{\mathbf Z}}
\newcommand\bGamma{\boldsymbol{\Gamma}}

\newcommand\bLambda{\boldsymbol{\Lambda}}

\newcommand\bSigma{\boldsymbol{\Sigma}}

\newcommand\bPhi{\boldsymbol{\Phi}}
\newcommand\bPsi{\boldsymbol{\Psi}}
\newcommand\bOm{\boldsymbol{\Omega}}

\newcommand\bfa{{\mathbf a}}
\newcommand\bfb{{\mathbf b}}
\newcommand\bfc{{\mathbf c}}
\newcommand\bfd{{\mathbf d}}

\newcommand\bff{{\mathbf f}}
\newcommand\bfg{{\mathbf g}}

\newcommand\bfm{{\mathbf m}}

\newcommand\bfp{{\mathbf p}}
\newcommand\bfq{{\mathbf q}}
\newcommand\bfr{{\mathbf r}}

\newcommand\bfv{{\mathbf v}}

\newcommand\bfx{{\mathbf x}}
\newcommand\bfy{{\mathbf y}}
\newcommand\bfz{{\mathbf z}}
\newcommand\bfalpha{\boldsymbol{\alpha}}

\newcommand\bfzeta{\boldsymbol{\zeta}}
\newcommand\bfeta{\boldsymbol{\eta}}

\newcommand\bfxi{{\boldsymbol \xi}}

\newcommand\bfphi{\boldsymbol{\phi}}
\newcommand\bfvarphi{\boldsymbol{\varphi}}

\newcommand\bfpsi{\boldsymbol{\psi}}

\newcommand\bfell{\boldsymbol{\ell}}

\newcommand\ubK{{\underline \bK}}

\newcommand\ubU{{\underline \bU}}

\newcommand\ubLambda{{\underline \bLambda}}

\newcommand\ubOm{{\underline \bOm}}

\newcommand\ubc{{\underline \bfc}}

\newcommand\tS{\widetilde{S}}

\newcommand\tLambda{\widetilde{\Lambda}}

\newcommand\tbK{\widetilde{\bK}}

\newcommand\tbU{\widetilde{\bU}}

\newcommand\tbW{\widetilde{\bW}}

\newcommand\tbY{\widetilde{\bY}}

\newcommand\tchi{\widetilde{\chi}}
\newcommand\tbp{\widetilde{\bfp}}
\newcommand\tbq{\widetilde{\bfq}}

\newcommand\cB{{\mathcal B}}
\newcommand\cC{{\mathcal C}}
\newcommand\cD{{\mathcal D}}
\newcommand\cE{{\mathcal E}}
\newcommand\cF{{\mathcal F}}
\newcommand\cG{{\mathcal G}}
\newcommand\cH{{\mathcal H}}

\newcommand\cK{{\mathcal K}}

\newcommand\cM{{\mathcal M}}
\newcommand\cN{{\mathcal N}}
\newcommand\cO{{\mathcal O}}
\newcommand\cP{{\mathcal P}}
\newcommand\cQ{{\mathcal Q}}
\newcommand\cR{{\mathcal R}}

\newcommand\cU{{\mathcal U}}

\newcommand\cW{{\mathcal W}}
\newcommand\cX{{\mathcal X}}

\newcommand\bcK{\boldsymbol{\cK}}

\newcommand\bcU{\boldsymbol{\cU}}

\newcommand\bcW{\boldsymbol{\cW}}

\usepackage{enumitem}
\title{
Phase sinks and sources\\around two-dimensional periodic-wave solutions\\of reaction-diffusion-advection systems
}

\author{Benjamin Melinand}
\address{CEREMADE, CNRS, Universit\'e Paris-Dauphine, Universit\'e PSL, 75016 Paris, France}
\email{{\tt melinand@ceremade.dauphine.fr}}
\thanks{}

\author{L.~Miguel Rodrigues}
\address{
Univ Rennes, CNRS, IRMAR - UMR 6625, F-35000 Rennes, France}
\email{{\tt luis-miguel.rodrigues@univ-rennes1.fr}}
\thanks{Research of L.M.R. was supported by the Institut Universitaire de France and partially supported by France 2030 through the programme Centre Henri Lebesgue ANR-11-LABX-0020-01.}

\begin{document}

\begin{abstract}
We develop a complete stability theory for two-dimensional periodic traveling waves of reaction-diffusion systems. More precisely, we identify a diffusive spectral stability assumption, prove that it implies nonlinear stability and provide a sharp asymptotic description of the dynamics resulting from both localized and critically nonlocalized perturbations. In particular, we show that the long-time behavior is governed at leading order by a second-order Whitham modulation system and elucidate how the intertwining of diffusive and dispersive effects may enhance decay rates. The latter requires a non trivial extension of the large-time estimates for constant-coefficient hyperbolic-parabolic operators to some classes of systems with no particular structure, including on one hand systems with a scalar-like --- but not scalar --- hyperbolic part and a cross-diffusion, and on the other hand anisotropic systems with dispersion.

\vspace{0.5em}


\vspace{0.5em}

\end{abstract}

\date{\today}
\maketitle

%\clearpage
\tableofcontents
%\clearpage
%%%%%%%%%%%%%%%%%

\section{Introduction}\label{s:introduction}

We initiate here a general programme aiming at a complete stability theory for genuinely multi-dimensional periodic traveling waves of parabolic systems. By a stability theory, we mean general results --- or at least a systematic approach --- that on one hand convert suitably defined spectral stability into nonlinear asymptotic stability in a suitable sense and on the other hand provide large-time asymptotics for the dynamics about such stable waves. Our goal is to extend the comprehensive theory now available for plane periodic waves to the multidimensional context. Concerning the latter we refer the reader to \cite{JNRZ-conservation} and references therein for a general picture and to \cite{JNRZ-RD1,JNRZ-RD2} for the pieces of work that are the most closely related to the present analysis.

In the present contribution, we focus on the case when spatial variables vary in $\R^2$ and the equations form a reaction-diffusion-advection system. Namely we consider 
\be\label{e:rd-intro}
\bcW_t=\Delta \bcW+\Div\bG(\bcW)+\bff(\bcW)\,,
\ee
for the $\R^n$-valued unknown $\bcW$, $\bcW(t,\bfx)\in\R^n$ (with $n\in\N^\star$), where $t$ denotes time variable and $\bfx\in\R^2$ denotes spatial variable. In \eqref{e:rd-intro}, we identify vectors of $\R^n$ with column vectors --- that is, with elements of $\cM_{n,1}(\R)$---; flux and source nonlinearities $\bG$ and $\bff$ are smooth functions on $\R^n$ valued respectively in $\cM_{2,n}$ and $\R^n\cong \cM_{n,1}(\R)$; the spatial divergence operator $\Div$ acts row-wise and the spatial Laplacian $\Delta$ is scalar\footnote{In the sense that it acts component-wise, with the same action on each component.}. For more details and further conventions concerning vectorial and differential notation we refer the reader to the devoted section, Section~\ref{s:notation}.

We study the general form in System~\eqref{e:rd-intro} as a compromise between generality and readibility. We claim however that, beyond some form of parabolicity, only two features of~\eqref{e:rd-intro} matter: the fact that coefficients depend neither on time nor on space variables; the fact, implicit here but encoded in \cond2, that the source term is non degenerate so that no hidden conservation law stems from~\eqref{e:rd-intro}. Both assumptions are deeply reflected in the structure of periodic traveling waves expounded below. To support the claim that the detailed structure of the original system is almost immaterial, let us anticipate on our analysis and list possible generalizations, by increasing order of difficulty. By a change of spatial variables, one may reduce any scalar symmetric elliptic operator to the scalar Laplacian. Changes required to analyze the general second-order semilinear parabolic case, including cross diffusions, are mostly notational. This simple observation turns out to be crucial to cover many of the examples in the literature that we give below. Studying the general second-order quasilinear parabolic case may be done along the same lines by increasing by one the level of regularity of involved solutions. Similarly the analysis of quasilinear parabolic cases of other orders differs mostly by the level of regularity of solutions. For some detailed examples of adaptations of the plane-wave analysis, the reader is referred to \cite{BJNRZ-KS} on a fourth-order semilinear parabolic equation and to \cite{RZ} on a quasilinear system that is parabolic only in the sense of some averaged version of the Kawashima condition. 

A (uniformly) traveling-wave solution to \eqref{e:rd-intro} is a solution $\bcW$ in the form $\bcW(t,\bfx)=\bcU(\bfx-t\,\bfc)$, where $\bcU$ is the wave profile and $\bfc\in\R^2$ is the wave speed. We say that the wave is periodic if its profile $\bcU$ is periodic. We are specifically interested in the case when $\bcU$ is genuinely multi-dimensional so that its group of periods is discrete, thus may be written as $\bX_1\,\Z+\bX_2\,\Z$ for some basis of $\R^2$, $(\bX_1,\bX_2)$. In the latter situation, we say that $\bcU$ is $(\bX_1,\bX_2)$-periodic. Alternatively one may scale the group of periods to be $\Z^2$ by introducing wave vectors $(\bK_1,\bK_2)$, given as the dual basis of $(\bX_1,\bX_2)$. As a result, a two-dimensional periodic wave is equivalently defined as a solution $\bcW$ of the form 
\be\label{e:wave-intro}
\bcW(t,\bfx)\,=\,\bU\left(\transp{\bK}\,(\bfx-t\,\bfc)\right)\,
=\,\bU\left(\transp{\bK}\,\bfx+t\,\bOm\right)\,, 
\ee
with $\bK=\bp \bK_1&\bK_2 \ep\in\cM_{2,2}(\R)$ a\footnote{The choice of $(\bX_1,\bX_2)$, or equivalently of $(\bK_1,\bK_2)$, is not canonical, but it is locally unique.} matrix of wave vectors, $\bOm=-\transp{\bK}\bfc\in\R^2$ a temporal frequency vector, $\bfc$ the wave speed and $\bU$ an associated (scaled) wave profile normalized to satisfy
\[
\bU(\,\cdot+\beD_j\,)=\bU\,,\qquad j=1,2\,,
\]
where $(\beD_1,\,\beD_2)$ is the canonical basis of $\R^2$. When $\bcW$ is given by \eqref{e:wave-intro}, it solves \eqref{e:rd-intro} if and only if 
\be\label{stand_eq}
0\,=\,\transp{(\bK\nabla)}(\bK\nabla)\bU+\transp{(\bK\nabla)}\bG(\bU) + (\transp{\bK}\bfc \cdot \nabla) \bU + \bff(\bU)\,.
\ee

In our analysis, we shall take as an assumption the existence of one specific wave, spectrally stable in a suitable sense. Yet the reader may wonder what is the relevance of this kind of objects and whether there is a robust universal mechanism supporting the existence of such objects. We claim that this is indeed the case and that such objects are somehow ubiquitous. To support the claim, we briefly recall, in words of \cite{R}, one of the prominent paradigms of the general field including studies in pattern formation, coherent structure, nonlinear waves and hydrodynamic instabilities. 

Quite often transition to instability of a certain form of solutions often gives rise to a new family of patterns whose stability may in turn also be investigated. Hence the classical strategy --- for equations involving some set of parameters --- consisting in carrying a parametric study of stability/instability. Starting from a simple family of solutions, explicit or even trivial, known to be stable for a certain range of parameters, one varies these parameters up to a transition to instability. At this threshold emerges a new family of special solutions, whose stability is also tracked when varying parameters and that can also yield yet another family of solutions, and so on and so forth. The patterns emerging from the first transition are usually called primary instabilities, those coming next secondary instabilities. Although one may artificially build systems exhibiting an infinite number of such transitions, it seems that in most of classical physical problems the instability of secondary patterns leads rather to chaos then turbulence. An argument supporting this phenomenological rule of thumb is that the emergence of new patterns often goes with a symmetry breaking increasing the dimensional complexity: trivial solutions are zero dimensional, primary instabilities one-dimensional, secondary ones two-dimensional, then comes chaos. Two-dimensional periodic waves studied here typically arise as secondary instabilities, the role of primary instabilities being played by plane periodic traveling waves, but may also emerge directly as primary instabilities from constant states. From this point of view the forthcoming \cite{BdR-R}, that studies the bifurcation of two-dimensional periodic waves from plane waves, appears as a companion paper. We refer the reader interested in supporting examples and further developments of the foregoing notions to \cite{CH93,Manneville-instabilities,Cross-Greenside,Charru-hydrodynamic-stability} and to \cite{UW14} the reader interested in an example of a numerical parametric study in a context close to our nonlinear analysis. 

From the point of view of mathematical analysis the near-constant study is in many ways more tractable than the near plane-wave one. Correspondingly, the mathematical literature devoted to proofs of existence of two-dimensional periodic waves is overwhelmingly focused on their arising from constant states, in particular through Turing bifurcations. As illustrated by \cite{Knobloch}, a large part of this literature outgrows from the trailblazing of equivariant bifurcation in \cite{Sattinger_bifurcation,Golubitsky-Stewart-Schaeffer_II}. We refer to \cite[Section~2]{DSSS03} for a short review of this kind of analysis. To help the reader navigate through this primary-instability literature we add a few general comments on the bifurcation of small-amplitude two-dimensional periodic patterns.
\begin{enumerate}
\item Many of the studied systems are invariant under the symmetry $\bfx\mapsto -\bfx$. As a result all small traveling waves built in this case are \emph{standing} waves, that is, $\bOm$ is identically zero along the family of waves. Such a symmetry happens for \eqref{e:rd-intro} when $\bG$ is assumed to be identically zero.
\item A significant part of the literature focuses on an even smaller class, the one of isotropic systems, that is, the one of systems invariant under the action of any linear rotation. Again this happens for \eqref{e:rd-intro} when $\bG$ is zero. As a consequence, in this case, built small-amplitude periodic waves arise with limiting wavevectors $(\bK_1^{(0)},\bK_2^{(0)})$ sharing the same norm. Generically, then, the only vectors of $\bK_1^{(0)}\Z+\bK_2^{(0)}\Z$ with this norm are $\bK_1^{(0)}$, $\bK_2^{(0)}$, $-\bK_1^{(0)}$, $-\bK_2^{(0)}$, those forming a rectangle so that corresponding standing waves are often referred to as (generalized) squares. Elementary geometry shows that the exceptions to the latter happen exactly when the angle between $\bK_1^{(0)}$ and $\bK_2^{(0)}$ is $\pi/N$ (or $-\pi/N$, $\pi+\pi/N$, $\pi-\pi/N$) for some integer $N\geq 3$, so that the vectors of the lattice on the prescribed circle form a regular polygon with $2N$-vertices. Standing waves arising from the case $N=3$ are often referred to as hexagon patterns. The exceptionally regular case is associated with a higher-dimensional limiting kernel but this dimension is reduced by enforcing extra symmetries on the sought pattern. We emphasize that all cases are equally captured by our nonlinear analysis (when suitable spectral stability is met).
\item The small-amplitude existence studies are often completed by a stability diagram. We warn the reader that the stability analyzed in the literature is restricted to perturbations with the symmetries of the pattern, in particular with the same periodicity. This is by far a much simpler task that the one we tackle here, and may be deduced from the computation of a normal form on a suitable center manifold. In contrast the secondary-instability analysis of \cite{BdR-R} provides exactly the notion of spectral stability needed here. 
\end{enumerate}
Let us stress again that despite the absence, until now, of a complete mathematical treatment of dynamics on extended domains, multidimensional periodic patterns are currently observed in numerous real life contexts where spatial domains are far from resembling fundamental domains of the patterns. Jointly with \cite{BdR-R} the goal of the present contribution is to help bridging this gap. 

One reason to restrict our analysis to dimension two is that we believe that technical gaps concerning tools available to study stability issues lie on one hand between constant solutions and non-constant solutions --- zero-dimensional objects to one-dimensional objects --- and on the other hand between plane waves and genuinely multi-dimensional waves --- one-dimensional objects to two-dimensional objects ---. Therefore the present analysis is expected to be representative of other multi-dimensional analyses. As a first sign of this gap, we stress that a large part of the technical tools classically used in the analysis of plane waves hinges on spatial dynamics building from ODE interpretations of both profile equations (\eqref{stand_eq} here) and spectral problems, none of them being available in the present context (at least in an obvious way). \emph{A priori} this rules out techniques as common in the field as phase portrait analysis, Evans functions, accompanying pointwise bounds on Green functions, \emph{etc.}

\subsection{Stability}

From now on we pick a specific periodic wave solution and use underlining to denote wave quantities related to this specific wave, including $\ubU$, $\ubK$, $\ubc$, \emph{etc.} To analyze the dynamics near this specific wave it is convenient to work in an adapted co-moving frame. Introducing $\bW$ through 
\begin{align}\label{e:co-moving}
\bcW(t,\bfx)\,=\,\bW\left(t,\transp{\ubK}\,(\bfx-t\,\ubc)\right)
\end{align}
turns \eqref{e:rd-intro} into
\be\label{rd}
\bW_t=\transp{(\ubK\nabla)}(\ubK\nabla)\bW +\transp{(\ubK\nabla)} \bG(\bW) + (\transp{\ubK}\ubc \cdot \nabla) \bW +\bff(\bW).
\ee
By design, $(t,\bfx)\mapsto\ubU(\bfx)$ is a stationary $(\beD_1,\beD_2)$-periodic solution to \eqref{rd}. Linearizing \eqref{rd} about $\ubU$ yields the periodic-coefficient
equation $(\d_t-L)\bV=0$ with $L$ given by
\be\label{e:lin-op}
L\bV:= \transp{(\ubK\nabla)}(\ubK\nabla) \bV + \transp{(\ubK\nabla)} \dD\bG(\ubU)(\bV) + (\transp{\ubK}\ubc \cdot \nabla) \bV + \dD\bff(\ubU)(\bV)\,.
\ee
From a functional-analytic point of view, we shall consider $L$ as an operator on $L^2(\R^2;\R^n)$ with domain $H^2(\R^2;\R^n)$. 

By (variations on) classical arguments --- detailed in Appendix~\ref{s:spec} ---, based on a suitable integrable transform --- the Bloch transform ---, the analysis of the action of $L$ on functions over $\R^2$ is reduced to the study of a Bloch symbol $\bfxi\mapsto L_\bfxi$, that with each $\bfxi\in[-\pi,\pi]^2$ %\cong \R^2/(2\pi\Z)^2$ 
associates an operator on $(\beD_1,\beD_2)$-periodic functions. More explicitly, for each $\bfxi\in[-\pi,\pi]^2$, $L_\bfxi$ acts on $L^2([0,1]^2;\C^n)\cong L^2(\R^2/\Z^2;\C^n)$ with domain $H^2_{\rm per}(\R^2;\C^n)\cong H^2(\R^2/\Z^2;\C^n)$ through 
\be\label{e:symbol}
L_{\bfxi}\bV=
\transp{(\ubK(\nabla+\iD\bfxi))}\ubK(\nabla+\iD\bfxi)\bV +\transp{(\ubK(\nabla+\iD\bfxi))} \dD\bG(\ubU)(\bV) + (\transp{\ubK}\ubc \cdot (\nabla+\iD\bfxi))\bV + \dD\bff(\ubU)(\bV)\,,
\ee
and, as such, has compact resolvents hence discrete spectrum, reduced to eigenvalues of finite multiplicity. As a consequence of the Bloch-wave representation and the continuity of $\sigma(L_{\bfxi})$ with respect to variations in $\bfxi$, in particular, the following spectral decomposition holds
\be\label{e:spec-decomp}
\sigma(L)\,=\,\bigcup_{\bfxi\in[-\pi,\pi]^2} \sigma(L_\bfxi)\,.
\ee
See Appendix \ref{ss:pert} for a proof.

To motivate the definition of the relevant notion of spectral stability, we point out that it follows from translational invariance of \eqref{e:rd-intro} that for any $\bfvarphi_0\in\R^2$, $\ubU(\cdot+\bfvarphi_0)$ is also a periodic traveling-wave profile associated with $(\ubK,\ubc)$. As a consequence, differentiating the corresponding profile equations with respect to $\bfvarphi_0$ shows that $\d_1\ubU$ and $\d_2\ubU$ lie in the kernel of $L_{\b0}$. Note moreover that the independence of $\d_1\ubU$ and $\d_2\ubU$ is precisely the condition that the wave under consideration is a genuinely multi-dimensional wave (and not a plane-wave in disguise). Besides, the real symmetry of spectra (stemming from the fact that \eqref{e:rd-intro} has real coefficients) implies that the real part of the eigenvalues of $L_{\bfxi}$ arising from the zero eigenvalue of $L_{\b0}$ when $\bfxi$ is small cannot be of order $\|\bfxi\|$ unless the background is unstable. Therefore the best one may expect is \emph{diffusive spectral stability} in the sense of the following conditions 
\begin{enumerate}
  \item[\cond1] There exist $\theta>0$ and $C>0$ such that for any $\bfxi\in[-\pi,\pi]^2$ and any $t\geq0$
\[
\lnor \eD^{t\,L_{\bfxi}} \rnor\,\leq\,C e^{-\theta t\|\bfxi\|^2}\,,
\]
where $(\eD^{t\,L_{\bfxi}})_{t\geq 0}$ denotes the semigroup on $L^2([0,1]^2;\R^n)$ generated by $L_{\bfxi}$ and $\lnor\cdot\rnor$ stands for operator norms.
  \item[\cond2]
 The spectrum of $L_{\b0}$ intersects $\iD\R$ only at $\lambda=0$ and $\lambda=0$ is an eigenvalue of $L_{\b0}$ of algebraic multiplicity $2$, its generalized eigenspace $\Sigma_{\b0}$ being spanned by $\d_1\ubU$ and $\d_2\ubU$.
\end{enumerate}

Assumption~\cond2 encodes that the criticality of the co-periodic spectrum is minimal. In turn, we think assumption~\cond1 as at least two-fold. The fact that the bound holds when $\|\bfxi\|\geq \xi_0$ for some fixed $\xi_0>0$ (and $(\theta,C)$ depending on $\xi_0$) is equivalent to the fact that for any $\bfxi$ such that $\|\bfxi\|\geq \xi_0$, $\sigma(L_{\bfxi})\subset\left\{\ \lambda\ ;\ \Re(\lambda)<0\ \right\}$. Once this is known to hold for any $\xi_0>0$ and \cond2 is also enforced, \cond1 is equivalent to the modulation system, introduced below, being hyperbolic-parabolic in a suitable Kawashima sense. We detail the latter in Appendix~\ref{s:spec}. We have chosen to summarize all these aspects in the form \cond1 mostly because it is particularly convenient for our linear and nonlinear stability analysis but, in Appendix~\ref{s:spec}, we provide more concrete equivalent characterizations and even simpler sufficient conditions. It is important to note that, unlike what happens for one-dimensional waves \cite{JNRZ-RD1,JNRZ-RD2}, condition \cond1 is in general stronger than
\begin{enumerate}
  \item[\cond{0}] There exists $\theta>0$ such that for any $\bfxi\in[-\pi,\pi]^2$ we have
\[
\sigma(L_{\bfxi})\subset\left\{\ \lambda\ ;\ \Re(\lambda)\leq-\theta\|\bfxi\|^2\ \right\}\,.
\]
\end{enumerate}
The difference between \cond0 and \cond1 lies in uniform control on diagonalization/symmetrization near $(\lambda,\bfxi)=(0,\b0)$. To stress that something is at stake, we point out that in general one cannot hope for a consistent diagonalization, smooth in $\bfxi$, near $\bfxi=\b0$ for the eigenvalues arising from the double eigenvalue $\lambda=0$ at $\bfxi=\b0$, but instead, at best in general, one expects the diagonalization to be smooth in $(\|\bfxi\|,\frac{\bfxi}{\|\bfxi\|})$. Incidentally we observe that this lack of smoothness on the symbolic side is intrinsically tied to dispersive effects, potential enhancing decay in poorly-localized topologies but deteriorating it in localized topologies. For more details, we refer the reader to the distinction between Cases~\case{a} and~\case{b} below and the discussion and results surrounding it. 

Our first result converts spectral stability in the diffusive sense of \cond1-\cond2 into nonlinear asymptotic stability in a suitable sense. Though our statement contains a more detailed description, it should be thought as providing stability in the space-modulated sense of \cite{JNRZ-conservation} (see also \cite{R,R_Roscoff}). This consists in measuring in the classical definition of stability the proximity to $\ubU$ of both initial data and solutions at later times by 
\[
\inf_{\bPhi\textrm{ invertible}}\|\bU\circ\bPhi-\ubU\|_{\bX} + \|\nabla(\bPhi-\Id)\|_{\bY}
\]
where $(\bX,\bY)$ are some functional spaces (possibly different for initial data and the solution). For comparison, note that naive stability requires control on $\|\bU-\ubU\|_{\bX}$ whereas orbital stability requires control on 
\[
\inf_{\bPhi\textrm{ uniform translation}}\|\bU\circ\bPhi-\ubU\|_{\bX}\,.
\]
Note that it is possible to choose $\bY$ as a space of curl-free vector-fields, on which it is natural to choose $\|\,\cdot\,\|_\bY$ as $\|\bV\|_{\bY}:=\|\Div(\bV)\|_{\tbY}$ for some functional space $\tbY$.

\bt[Stability]\label{th:nonlinear_stab}
Let $\ubU$ be a stationary $(\beD_1,\beD_2)$-periodic solution to \eqref{rd} associated to the matrix of wave vectors $\ubK$ and speed $\ubc$. Assume \cond1-\cond2. There exists $\eps_0>0$ and $C>0$ such that if for some sublinear\footnote{By this, we mean that $\bfphi_0$ may differ from $\Delta^{-1}(\Delta \bfphi_0)$ by a constant function but not by a non-constant affine function. This is for instance the case if one enforces $\nabla\bfphi_0\in\Span\left(\bigcup_{1\leq p<\infty}L^p(\R^2;\cM_2(\R))\right)$.} $\bfphi_0$
\[
E_0:=\|\bW_{0}(\cdot-\bfphi_{0})-\ubU \|_{(H^{2} \cap W^{2,4})(\R^2;\R^n)}+\|\Delta \bfphi_{0}\|_{(L^{1} \cap W^{1,4})(\R^2;\R^2)}\,\leq\,\eps_0
\]
then, there exist a unique global solution to \eqref{rd} with initial datum $\bW_{0}$ and a phase shift $\bfphi$ with $\bfphi(0,\cdot)= \bfphi_{0}$ such that, for any $t\geq0$,
\begin{align*}
\|\bW(t,\cdot-\bfphi(t,\cdot))-\ubU\|_{W^{2,4}(\R^2;\R^n)}
+\|\nabla\bfphi(t,\cdot)\|_{W^{2,4}(\R^2;\cM_2(\R))}
+\|\d_t\bfphi(t,\cdot)\|_{W^{2,4}(\R^2;\R^2)}
&\leq \frac{C\,E_{0}}{(1+t)^{\frac{1}{4}}}\,.
\end{align*}
Furthermore, with constants independent of $(\bW_0,\bfphi_0)$ and no further restriction on $E_0$, 
\begin{enumerate}
\item for any $t\geq0$, 
\begin{align*}
\|\bW(t,\cdot-\bfphi(t,\cdot))-\ubU\|_{L^\infty(\R^2;\R^n)}
+\|\nabla\bfphi(t,\cdot)\|_{L^\infty(\R^2;\cM_2(\R))}
+\|\d_t\bfphi(t,\cdot)\|_{L^\infty(\R^2;\R^2)}
&\leq C\,E_0\,\frac{\ln(2+t)}{(1+t)^{\frac{1}{2}}}\,;
\end{align*}
\item for any $2<p_{0}<q_{0}<\infty$, there exists a constant $C_{p_{0},q_{0}}>0$, such that for any $p \in [p_{0},q_{0}]$, and any $t\geq0$
\[
\|\bW(t,\cdot-\bfphi(t,\cdot))-\ubU\|_{L^p(\R^2;\R^n)}
+\|\nabla\bfphi(t,\cdot)\|_{L^p(\R^2;\cM_2(\R))}
+\|\d_t\bfphi(t,\cdot)\|_{L^p(\R^2;\R^2)}
\,\leq \frac{C_{p_{0},q_{0}}\,E_{0}}{(1+t)^{\frac{1}{2}-\frac{1}{p}}}\,;
\]
\item for any $\ell\in\N$, $\ell\geq2$, there exists $C_\ell$ such that if moreover 
\[
E_{0,\ell}:=\|\bW_{0}(\cdot-\bfphi_{0})-\ubU \|_{(H^{2} \cap W^{\ell,4})(\R^2;\R^n)}+\|\Delta \bfphi_{0}\|_{(L^{1} \cap W^{\ell-1,4})(\R^2;\R^2)}\,<+\infty
\]
then for any $t\geq0$,
\begin{align*}
\|\bW(t,\cdot-\bfphi(t,\cdot))-\ubU\|_{W^{\ell,4}(\R^2;\R^n)}
+\|\nabla\bfphi(t,\cdot)\|_{W^{\ell,4}(\R^2;\cM_2(\R))}
+\|\d_t\bfphi(t,\cdot)\|_{W^{\ell,4}(\R^2;\R^2)}
&\leq \frac{C_\ell\,E_{0,\ell}}{(1+t)^{\frac{1}{4}}}\,.
\end{align*}
\end{enumerate}
\et

Note that adding an affine function to $\bfphi_{0}$ would alter the background matrix of wave vectors so that the arising solution should be compared with another periodic solution of \eqref{rd}. That is the main reason why we only deal with sublinear initial phases.

\begin{figure}[htb]
\begin{subfigure}{.4\textwidth}
\includegraphics[width=\textwidth]{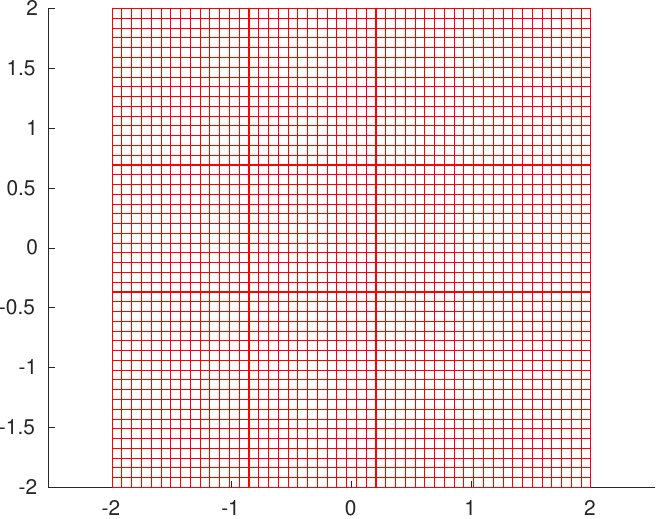}
\caption{Original grid}
\end{subfigure}%
\begin{subfigure}{.4\textwidth}
\includegraphics[width=\textwidth]{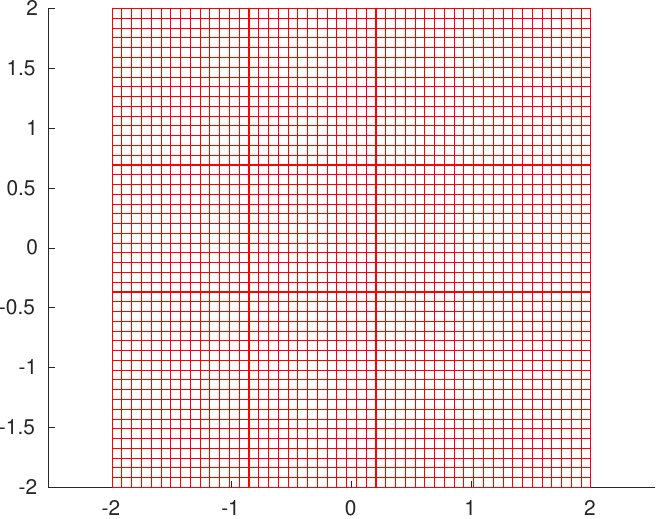}
\caption{Original grid}
\end{subfigure}
\begin{subfigure}{.4\textwidth}
\includegraphics[width=\textwidth]{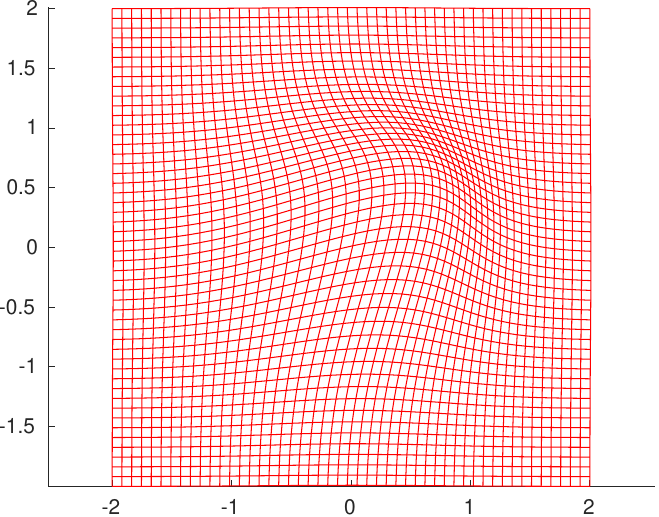}
\caption{$\bfphi$ of Gaussian type}
\end{subfigure}%
\begin{subfigure}{.4\textwidth}
\includegraphics[width=\textwidth]{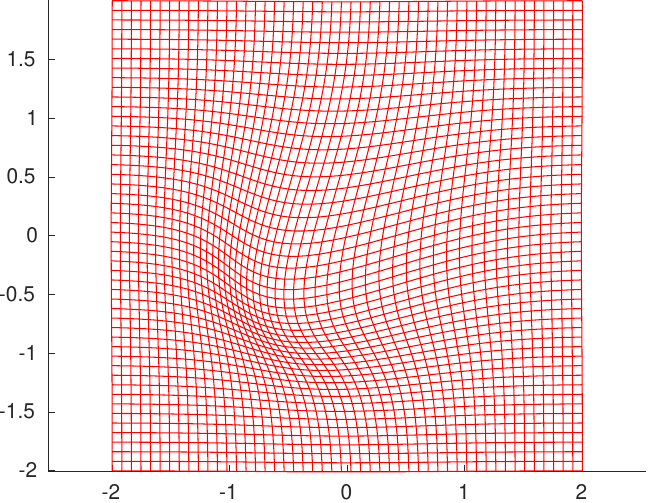}
\caption{$-\bfphi$ of Gaussian type}
\end{subfigure}
\begin{subfigure}{.4\textwidth}
\includegraphics[width=\textwidth]{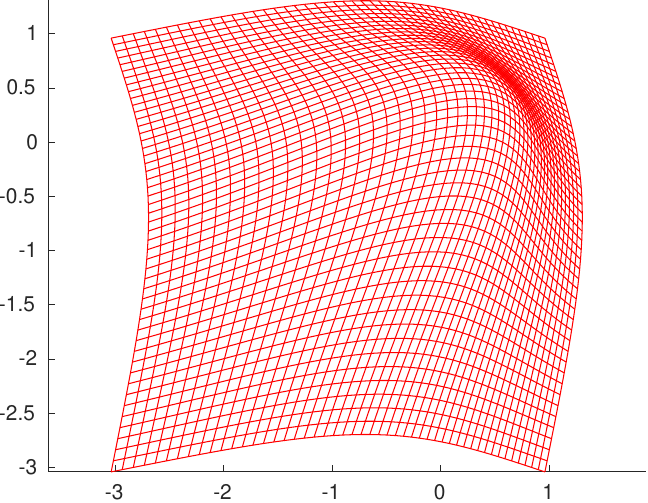}
\caption{$\Delta \bfphi$ of Gaussian type (source)}
\end{subfigure}%
\begin{subfigure}{.4\textwidth}
\includegraphics[width=\textwidth]{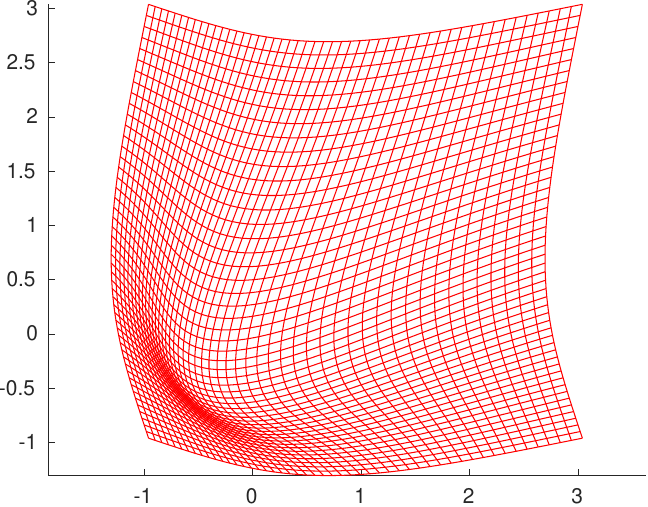}
\caption{$-\Delta \bfphi$ of Gaussian type (sink)}
\end{subfigure}
\caption{Illustration of the various levels of localizations. We plot the image of a square reference grid through the function $\text{Id}_{\R^{2}}-\bfphi$.}
\label{source-sink}
\end{figure}

What drives decay rates is the initial localization. From this point of view, the key part of the assumption is $\bW_{0}(\cdot-\bfphi_{0})-\ubU\in L^2$ and $\Delta \bfphi_{0}\in L^1$ (with small norms). The assumption $\Delta \bfphi_{0}\in L^1$ should be thought as a relaxation of $\nabla \bfphi_{0}\in L^2$. Indeed,  $\Delta \bfphi_{0}\in L^1$ implies that $\nabla \bfphi_{0}$ belongs to $L^{2,\infty}$, the weak-$L^2$ space. Moreover, if $(1+\| \cdot \|)\,\Delta \bfphi_{0}\in L^1$ and $\Delta \bfphi_{0}\in L^p$ for some $p>1$, then $\nabla \bfphi_{0}$ belongs to $L^2$ if and only if $\int_{\R^2} \Delta \bfphi_{0}\,=\,\b0$. In particular, the relaxed assumption allows to prescribe for $\Delta \bfphi_{0}$ any small regularized version of a multiple of a Dirac mass, hence the term phase source/sink (depending on the sign of the Dirac mass) in the title. See Figure \ref{source-sink}.

Likewise, $\Delta \bfphi_{0}\in L^1$ ensures that $\bfphi_{0}$ belongs to $BMO$, the space of functions with bounded mean oscillation, which may be thought as a relaxed version of $\bfphi_{0}\in L^\infty$, whereas, when, say, $(1+\| \cdot \|^2)\,\Delta \bfphi_{0}\in L^2$, one shows that $\bfphi_{0}$ belongs to $L^\infty$ if and only if $\int_{\R^2} \Delta \bfphi_{0}\,=\,\b0$. For proofs and further comments on the recovering of $\bfphi_{0}$ and $\nabla \bfphi_{0}$ from $\Delta \bfphi_{0}$, we refer to Appendix~\ref{s:phase-like}. 

The regime of localization chosen here is critical from the point of view of large-time asymptotic decay, in the sense that nonlinear terms are not asymptotically irrelevant. Moreover, our further study of asymptotic behavior identifies a leading-order nonlinear asymptotic description, allowing us to analyze the sharpness\footnote{In this direction, let us anticipate that in the end we remove essentially all $\log(2+t)$ factors at the cost of making stronger regularity assumptions. See for instance Remark~\ref{rk:log}.} of decay estimates stated above. 

In contrast, we make no claim on optimality of our regularity assumptions, encoded by the choice of the space $W^{2,4}$. We regard the regularity question as largely irrelevant for the problem at hand and, correspondingly, in parts of the proof where this plays a role, we have decided to apply simpler and/or more robust arguments instead of sharper ones. At a technical level, the choice of an $L^4$-based space is designed to ensure that quadratic terms lie in $L^2$ so that their contribution to large-time decay may be analyzed through Hausdorff-Young inequalities. Then, among $L^4$-based spaces, the choice of $W^{2,4}$ enforces embedding in $W^{1,\infty}$. Though we expect our regularity framework to be suboptimal, we warn the reader that one should not lower the regularity on $\nabla\bfphi$ below the threshold ensuring that $\Id-\bfphi$ is invertible.

For comparison and comprehensiveness we also provide a stability result under more localized perturbations. 

\bt[Subcritical perturbations]\label{th:nonlinear-localized}
Let $\ubU$ be a stationary $(\beD_1,\beD_2)$-periodic solution to \eqref{rd} associated to the matrix of wave vectors $\ubK$ and speed $\ubc$. Assume \cond1-\cond2. There exists $\eps_0>0$ and $C>0$ such that if for some $\bfphi_0$
\[
\cE_{0}:=\|\bW_{0}(\cdot-\bfphi_{0})-\ubU \|_{(L^{1}\cap H^2\cap W^{2,4})(\R^2;\R^n)}+\|\nabla \bfphi_{0}\|_{(L^{1} \cap H^2\cap W^{2,4})(\R^2;\cM_2(\R))}\,\leq\,\eps_0
\]
then, there exist a unique global solution to \eqref{rd} with initial datum $\bW_{0}$ and a phase shift $\bfphi$ with $\bfphi(0,\cdot)= \bfphi_{0}$ such that, for any $t\geq0$, for any $2\leq p\leq 4$
\begin{align*}
\|\bW(t,\cdot-\bfphi(t,\cdot))-\ubU\|_{W^{2,p}(\R^2;\R^n)}
+\|\nabla\bfphi(t,\cdot)\|_{W^{2,p}(\R^2;\cM_2(\R))}
+\|\d_t\bfphi(t,\cdot)\|_{W^{2,p}(\R^2;\R^2)}
&\leq \frac{C\,\cE_{0}}{(1+t)^{1-\frac{1}{p}}}\,.
\end{align*}
Furthermore, with constants independent of $(\bW_0,\bfphi_0)$ and no further restriction on $\cE_0$, 
\begin{enumerate}
\item there exists a constant $C>0$, such that for any $2 \leq p \leq \infty$ and any $t\geq0$
\[
\|\bW(t,\cdot-\bfphi(t,\cdot))-\ubU\|_{L^p(\R^2;\R^n)}
+\|\nabla\bfphi(t,\cdot)\|_{L^p(\R^2;\cM_2(\R))}
+\|\d_t\bfphi(t,\cdot)\|_{L^p(\R^2;\R^2)}
\,\leq \frac{C\,\cE_{0}}{(1+t)^{1-\frac{1}{p}}}\,,
\]
\item for any $\ell\in\N$, $\ell\geq2$, there exists $C_\ell$ such that if moreover 
\[
\cE_{0,\ell}:=\|\bW_{0}(\cdot-\bfphi_{0})-\ubU \|_{(L^{1}\cap H^2 \cap W^{\ell,4})(\R^2;\R^n)}+\|\nabla \bfphi_{0}\|_{(L^{1} \cap H^2 \cap W^{\ell,4})(\R^2;\cM_2(\R))}\,<+\infty \,,
\]
then for any $t\geq0$, and $2\leq p\leq 4$
\begin{align*}
\|\bW(t,\cdot-\bfphi(t,\cdot))-\ubU\|_{W^{\ell,p}(\R^2;\R^n)}
+\|\nabla\bfphi(t,\cdot)\|_{W^{\ell,p}(\R^2;\cM_2(\R))}
+\|\d_t\bfphi(t,\cdot)\|_{W^{\ell,p}(\R^2;\R^2)}
&\leq \frac{C_\ell\,\cE_{0,\ell}}{(1+t)^{1-\frac{1}{p}}}\,.
\end{align*}
\end{enumerate}
Finally, there exists a constant $\bfphi_{\infty}$ depending only on $\bfphi_{0}$ such that for any $t \geq 0$
\[
\|\bW(t,\cdot)-\ubU \|_{L^{p}(\R^2;\R^n)} + \| \bfphi(t,\cdot)- \bfphi_{\infty} \|_{L^{p}(\R^2;\R^2)} \leq \frac{C \cE_{0}}{(1+t)^{\frac{1}{2}-\frac{1}{p}}}  \,, \qquad 2  \leq p \leq +\infty,
\]
and if $\bfphi_{0} \in L^{1}(\R^{2})$ and $\cE_{0}$ is small enough, for any $t \geq 0$
\[
\|\bW(t,\cdot)-\ubU \|_{L^{p}(\R^2;\R^n)} \leq C (\cE_{0} + \| \bfphi_{0}  \|_{L^{1}(\R^{2})}) \frac{\ln(2+t)}{(1+t)^{1-\frac{1}{p}}}  \,, \qquad 2  \leq p \leq +\infty.
\]
\et

We stress that assuming more localization on initial data, including enforcing $\bfphi_0\equiv 0$, would not bring extra decay. Moreover, though the extra localization assumed here does bring some minor simplifications, the scheme of proof of Theorem~\ref{th:nonlinear-localized} is not significantly different from the one for Theorem~\ref{th:nonlinear_stab}. The reason for that is that the corresponding regime of decay is barely subcritical and, thus, some care is needed to carry out the argument. With this respect, it is instructive to compare the two-dimensional analysis in \cite{Johnson-Zumbrun-cons-1D-2D} with the three-dimensional analysis in \cite{Oh-Zumbrun-cons-nonlinear}.

\subsection{Modulational behavior}

Theorem~\ref{th:nonlinear_stab} contains that, up to a remainder of size $t^{-\left(\tfrac12-\tfrac1p\right)}$ in $L^p$, $2<p<\infty$, the solution $\bcW(t,\cdot)$ to the original \eqref{e:rd-intro} is well-described by $\ubU\circ \bPsi(t,\cdot)$ for some $\bPsi$ such that $\nabla\bPsi-\ubK$ and $\d_t\bPsi-\ubOm$ are also decaying at the same remainder rate in $L^p$.

To go further, one would like to capture the leading-order part of the near-constant dynamics of $(\nabla \bPsi,\d_t\bPsi)$. It turns out that this is closely related to the obtention of a refined description of $\bcW(t,\cdot)$, that is, of a description up to a faster-decaying remainder. The latter requires not only a space-time modulation of the position of the wave profile $\ubU$ but also of its shape, hence, as a preliminary, an understanding of nearby waves. 

As we prove in Proposition~\ref{p:structure}, Assumption~\cond2 is sufficient to elucidate the structure of nearby two-dimensional periodic waves. The upshot is that corresponding profiles, wavenumbers and speeds may be smoothly parametrized as $(\bK,\bfvarphi) \mapsto(\bU^{\bK}(\cdot+\bfvarphi),\bK,\bfc^{\bK})$. The following result shows that by modulating also in wavenumber, besides the modulation in position, one does improve the asymptiotic description of solutions.

\bt[Modulational behavior]\label{th:mod-behavior}
Assume \cond1-\cond2 and consider a wave parametrization as in Proposition~\ref{p:structure}.\\
With notational conventions of Theorem~\ref{th:nonlinear_stab}, one may also ensure that the solution $\bcW$ to \eqref{e:rd-intro} obtained from $\bW$ through \eqref{e:co-moving} satisfies, 
\begin{align*}
\|\bcW(t,\cdot)-\bU^{\bcK(t,\cdot)}(\bPsi(t,\cdot))\|_{L^{p}(\R^2;\R^n)}
&\leq C_{p_0,q_0}\,E_{0}\,\frac{\ln(2+t)}{(1+t)^{\frac{1}{2}-\frac{1}{p}+\frac{1}{2}\left(\frac12+\frac1p\right)}}\,,&
2<p_0\leq p\leq q_0<\infty\,,\\
\|\bcW(t,\cdot)-\bU^{\bcK(t,\cdot)}(\bPsi(t,\cdot))\|_{L^{\infty}(\R^2;\R^n)}
&\leq C\,E_{0}\,\frac{(\ln(2+t))^2}{(1+t)^{\frac{1}{2}+\frac14}}\,,&
\end{align*}
with $\bcK:=\nabla\bPsi$ and $\bPsi$ defined from $\bfphi$ as
\[
\bPsi(t,\bfx):=(\Id-\bfphi(t,\cdot))^{-1}\left(\transp{\ubK}\,(\bfx-t\,\ubc)\right)\,,
\]
thus satisfying
\begin{align*}
\|\bcK(t,\cdot)-\ubK\|_{L^p(\R^2;\cM_2(\R))}
+\|\d_t\bPsi(t,\cdot)-\ubOm\|_{L^p(\R^2;\R^2)}
&\,\leq \frac{C_{p_{0},q_{0}}\,E_{0}}{(1+t)^{\frac{1}{2}-\frac{1}{p}}}\,,&
2<p_0\leq p\leq q_0<\infty\,,\\
\|\bcK(t,\cdot)-\ubK\|_{L^\infty(\R^2;\cM_2(\R))}
+\|\d_t\bPsi(t,\cdot)-\ubOm\|_{L^\infty(\R^2;\R^2)}
&\leq C\,E_0\,\frac{\ln(2+t)}{(1+t)^{\frac{1}{2}}}\,.
\end{align*}
\et

Note that Theorem~\ref{th:mod-behavior} encodes through $\bcK=\nabla\bPsi$ that modulation in wavenumbers result from spatial variations of modulation in positions. 

We believe that the proof of Theorem~\ref{th:mod-behavior} is both robust and representative of what could be expected in much more general situations. Yet, as we show below, in the present two-dimensional case, the estimates of Theorem~\ref{th:mod-behavior} are deceptively pessimistic. Indeed, Theorem~\ref{th:mod-behavior} is the result of the combination of various worst-case bounds, whereas, in the present two-dimensional case under study, it turns out that not all the difficulties may be present simultaneously.

To elucidate this, let us anticipate on the description of the near-constant dynamics of $\bcK(t,\cdot)$. As we prove below, at leading order, this dynamics obeys an hyperbolic-parabolic system. In full generality, the hyperbolic part of such systems (linearized about the reference constant state) contains both scalar-type components and dispersive-type components. For the full linearized hyperbolic-parabolic system, one expects the former to decay exactly as solutions of the heat equation, that is, in dimension $2$, as $t^{-\left(\tfrac1q-\tfrac1p\right)}$, in $L^p$ for initial data in $L^q$, $1\leq q\leq p\leq+\infty$. In contrast, for the latter, one expects wave-diffusion type decay; for instance, in $L^p$, starting from $L^1$ data, the dispersion enhances the decay when $p>2$ but slows it down when $p<2$. For general systems, the decay of the full solution is thus prescribed by the worst rates between heat-like and wave-diffusion-like decay rates. We refer the reader to the detailed analysis of the isentropic compressible Navier-Stokes system in \cite{Hoff_Zumbrun-NS_compressible_ponctuel,Hoff_Zumbrun-NS_compressible_pres_de_zero,Kobayashi-Shibata,Rodrigues-compressible} (and further comments in \cite{Rodrigues-these}) for a worked-out example supporting such intuition.

However, as we prove in\footnote{We expect this fact to be well-known by experts of hyperbolic systems as it is related to the Strang analysis \cite{Strang_hyp} of systems of two equations in arbitrary dimension. Yet we have not found it in the literature.} Lemma~\ref{l:hyp}, in dimension two, constant-coefficient hyperbolic systems of two equations are either composed of two uncoupled scalar equations or they are strictly hyperbolic and fully of dispersive type. When the hyperbolic part of the effective dynamics for $\bcK(t,\cdot)$ is of dispersive type (Case~\case{a} below), it turns out that the decay are faster than proved in Theorems~\ref{th:nonlinear_stab} and~\ref{th:mod-behavior}, and the dynamics is actually asymptotically linear. When it is of scalar-like type (Case~\case{b} below), the estimates of Theorem~\ref{th:nonlinear_stab} are sharp\footnote{Up to logarithmic factors.} but those of Theorem~\ref{th:mod-behavior} are not because they rely on some linear estimates that may be improved in the scalar-type case but seemingly\footnote{Actually, finer comparisons with \cite{Hoff_Zumbrun-NS_compressible_ponctuel,Kobayashi-Shibata} suggest that there is some room for improvement here also but at a high technical price and in a way essentially useless at the nonlinear level.} not in the dispersive-type case.

To make this discussion more concrete, let us point out that the relevant linearized first-order dynamics is 
\be\label{e:linear-group-intro}
\d_t \bfpsi-\dD_{\bK}\bOm(\ubK)\,(\nabla \bfpsi)\,=\,\b0\,,
\ee
where $\nabla \bfpsi$ plays the role of a linear approximation of $\bcK-\ubK$. Note that $\bcK$ is $\cM_2(\R)$-valued but curl-free, hence satisfying two constraints, so that, as encoded in \eqref{e:linear-group-intro}, its dynamics is effectively two-dimensional. As we prove in Lemma~\ref{l:robust}, under conditions~\cond1-\cond2, System~\ref{e:linear-group-intro} is hyperbolic and therefore one, and only one, of the following two conditions hold.
\begin{enumerate}
  \item[\case{a}] For any unitary $\bfxi_0\in\R^2$, $\sum_{j=1}^2\dD_{\bK_j}\bOm(\ubK)\,(\bfxi_0)\,\transp{\beD_j}$ has real distinct eigenvalues. 
  \item[\case{b}] Matrices $\sum_{j=1}^2\dD_{\bK_j}\bOm(\ubK)\,(\beD_1)\,\transp{\beD_j}$ and $\sum_{j=1}^2\dD_{\bK_j}\bOm(\ubK)\,(\beD_2)\,\transp{\beD_j}$ are simultaneously diagonalizable over $\R$.
\end{enumerate}
Let us point out that in the special case where $\bfc^{\bK}$ does not depend on $\bK$ (or more generally when $\dD_{\bK}\bfc(\ubK)$ is zero), System \eqref{e:linear-group-intro} reduces to $\d_t \bfpsi+(\ubc\cdot\nabla)\bfpsi\,=\,\b0$, a genuinely scalar system, yielding Case~\case{b} in a trivial way. We identify this subcase as
\begin{enumerate}
  \item[\case{b0}] Matrices $\sum_{j=1}^2\dD_{\bK_j}\bOm(\ubK)\,(\beD_1)\,\transp{\beD_j}$ and $\sum_{j=1}^2\dD_{\bK_j}\bOm(\ubK)\,(\beD_2)\,\transp{\beD_j}$ are scalar.
\end{enumerate}
and we shall prove for it slightly sharper estimates with significantly simpler proofs.

Case~\case{a} is arguably the hardest case to analyze but we expect that it is also the generic one in the absence of extra symmetry. In particular this is the one proved to arise at the secondary bifurcation studied in \cite{BdR-R}. We recall however that in the isotropic case small-amplitude waves fit in Subcase~\case{b0}.

The following theorems prove the above claims about decay rates specialized to either Case~\case{a} or Case~\case{b}. 

\bt[Dispersive case]\label{th:dispersive}
Assume \cond1-\cond2 and Case~\case{a} of the alternative.\\
Assume that $E_{0}$ is small enough. One may improve the estimates of Theorem~\ref{th:nonlinear_stab} into 
\begin{align*}
\|\bW(t,\cdot-\bfphi(t,\cdot))-\ubU\|_{W^{2,4}}
+\|\nabla\bfphi(t,\cdot)\|_{W^{2,4}}
+\|\d_t\bfphi(t,\cdot)\|_{W^{2,4}}
&\leq \frac{C\,E_{0}}{(1+t)^{\frac{3}{8}}}\,,\\
\|\bW(t,\cdot-\bfphi(t,\cdot))-\ubU\|_{L^p}
+\|\nabla\bfphi(t,\cdot)\|_{L^p}
+\|\d_t\bfphi(t,\cdot)\|_{L^p}
&\leq \frac{C_{p_{0},q_0}\,E_{0}}{(1+t)^{\frac{3}{4}-\frac32\frac{1}{p}}}\,,
&2<p_0\leq p\leq q_0<\infty\,,\\
\|\bW(t,\cdot-\bfphi(t,\cdot))-\ubU\|_{L^\infty}
+\|\nabla\bfphi(t,\cdot)\|_{L^\infty}
+\|\d_t\bfphi(t,\cdot)\|_{L^\infty}
&\leq C\,E_{0}\frac{\ln(2+t)}{(1+t)^{\frac{3}{4}}}\,,&\\
\|\bW(t,\cdot-\bfphi(t,\cdot))-\ubU\|_{W^{\ell,4}}
+\|\nabla\bfphi(t,\cdot)\|_{W^{\ell,4}}
+\|\d_t\bfphi(t,\cdot)\|_{W^{\ell,4}}
&\leq \frac{C_{\ell}\,E_{0,\ell}}{(1+t)^{\frac{3}{8}}}\,,
&\ell\in\N,\,\ell\geq2,
\end{align*}
and those of Theorem~\ref{th:mod-behavior} into \footnote{We refrain from stating estimates for $\|\bcW(t,\cdot)-\bU^{\bcK(t,\cdot)}(\bPsi(t,\cdot))\|_{L^{p}}$ when $p>4$ because the corresponding decay rates would be artificially limited by a lack of regularity assumption on the data; see the related Remark \ref{rk:log}.}
\begin{align*}
\|\bcW(t,\cdot)-\bU^{\bcK(t,\cdot)}(\bPsi(t,\cdot))\|_{L^{p}}
&\leq \,\frac{C_{p_0,p_{1}}\,E_{0}}{(1+t)^{\frac{3}{4}-\frac32\frac{1}{p}+\frac{1}{2}}}\,,&
2<p_0\leq p \leq p_{1} < 4\,,\\
\|\bcW(t,\cdot)-\bU^{\bcK(t,\cdot)}(\bPsi(t,\cdot))\|_{L^{4}}
&\leq \,\frac{C \,E_{0} \ln(2+t)}{(1+t)^{\frac78}}\,,
\end{align*}
with $\bcK:=\nabla\bPsi$ and $\bPsi$ defined from $\bfphi$ as
\[
\bPsi(t,\bfx):=(\Id-\bfphi(t,\cdot))^{-1}\left(\transp{\ubK}\,(\bfx-t\,\ubc)\right)\,,
\]
thus satisfying
\begin{align*}
\|\bcK(t,\cdot)-\ubK\|_{L^p}
+\|\d_t\bPsi(t,\cdot)-\ubOm\|_{L^p}
&\,\leq \frac{C_{p_{0},q_{0}}\,E_{0}}{(1+t)^{\frac{3}{4}-\frac32\frac{1}{p}}}\,,&
2<p_0\leq p\leq q_0<\infty\,,\\
\|\bcK(t,\cdot)-\ubK\|_{L^\infty}
+\|\d_t\bPsi(t,\cdot)-\ubOm\|_{L^\infty}
&\leq C\,E_0\,\frac{\ln(2+t)}{(1+t)^{\frac{3}{4}}}\,.
\end{align*}
\et

\bt[Dispersive case, subcritical perturbations]\label{th:dispersive-localized}
Assume \cond1-\cond2 and Case~\case{a} of the alternative. Assume that $\cE_{0}$ is small enough. One may improve the estimates of Theorem~\ref{th:nonlinear-localized} into 
\begin{align*}
\|\bW(t,\cdot-\bfphi(t,\cdot))-\ubU\|_{W^{2,p}}
+\|\nabla\bfphi(t,\cdot)\|_{W^{2,p}}
+\|\d_t\bfphi(t,\cdot)\|_{W^{2,p}}
&\leq \frac{C\,\cE_{0}}{(1+t)^{\frac54-\frac32\frac1p}}\,,& 2\leq p\leq 4\,,\\
\|\bW(t,\cdot-\bfphi(t,\cdot))-\ubU\|_{L^p}
+\|\nabla\bfphi(t,\cdot)\|_{L^p}
+\|\d_t\bfphi(t,\cdot)\|_{L^p}
&\leq \,\frac{C\,\cE_0}{(1+t)^{\frac54-\frac32\frac1p}}\,,&
2\leq p\leq \infty\,,\\
\|\bW(t,\cdot-\bfphi(t,\cdot))-\ubU\|_{W^{\ell,p}}
+\|\nabla\bfphi(t,\cdot)\|_{W^{\ell,p}}
+\|\d_t\bfphi(t,\cdot)\|_{W^{\ell,p}}
&\leq \frac{C_\ell\,\cE_{0,\ell}}{(1+t)^{\frac54-\frac32\frac1p}}\,,&
\ell\in\N,\,\ell\geq2,\,2\leq p\leq4\,.
\end{align*}
\et

\bt[Scalar-type case]\label{th:scalar}
Assume \cond1-\cond2 and Case~\case{b} of the alternative. One may remove the $\log$-factor in the estimates of Theorem~\ref{th:nonlinear_stab} 
\begin{align*}
\|\bW(t,\cdot-\bfphi(t,\cdot))-\ubU\|_{L^p}
+\|\nabla\bfphi(t,\cdot)\|_{L^p}
+\|\d_t\bfphi(t,\cdot)\|_{L^p}
&\leq \frac{C_{p_{0}}\,E_{0}}{(1+t)^{\frac{1}{2}-\frac1p}}\,,
&2<p_0\leq p\leq \infty\,.\\
\end{align*}
Moreover,  if
\[
E_0^{(p)}:=\|\bW_{0}(\cdot-\bfphi_{0})-\ubU \|_{(H^{2}\cap W^{2,p})(\R^2;\R^n)}+\|\Delta \bfphi_{0}\|_{(L^{1} \cap W^{1,p})(\R^2;\R^2)} < +\infty \,,
\]
one can improve the estimates of Theorem~\ref{th:mod-behavior} into
\begin{align*}
\|\bcW(t,\cdot)-\bU^{\bcK(t,\cdot)}(\bPsi(t,\cdot))\|_{L^{p}}
&&\leq \,
\begin{cases}\ds
\frac{\ln(2+t)}{(1+t)^{1-\frac{1}{p}}}\ C_{p_0}\,E_0\,,&
2<p_0\leq p\leq 4\,,\\[0.75em]\ds
\frac{\ln(2+t)}{(1+t)^{1-\frac{1}{p}}}\ C_{p_1}\,E_0^{(p_1)}\,\,,&
4\leq p\leq p_1 < \infty \,,\\[0.75em]\ds
\frac{1}{(1+t)^{1-\frac{1}{p_1}}}\ C_{p_1,p_2}\,E_0^{(p_1)}\,,&
p_1<p_2\leq p\leq \infty\,,
\end{cases}
\end{align*}
with $\bcK:=\nabla\bPsi$ and $\bPsi$ defined from $\bfphi$ as
\[
\bPsi(t,\bfx):=(\Id-\bfphi(t,\cdot))^{-1}\left(\transp{\ubK}\,(\bfx-t\,\ubc)\right)\,,
\]
thus satisfying
\begin{align*}
\|\bcK(t,\cdot)-\ubK\|_{L^p}
+\|\d_t\bPsi(t,\cdot)-\ubOm\|_{L^p}
&\,\leq \frac{C_{p_{0}}\,E_{0}}{(1+t)^{\frac{1}{2}-\frac{1}{p}}}\,,&
2<p_0\leq p\leq \infty\,.
\end{align*}
\et

As already implicitly pointed out, the decay rates in Theorem~\ref{th:scalar} should be compared with those of solutions to the heat equation whereas those of Theorems~\ref{th:dispersive} and~\ref{th:dispersive-localized} should be compared with those for the viscously damped wave equation. Roughly\footnote{This means, in particular, that at this informal level we do not bother to state admissibility conditions for $\alpha$, $\ell$, $p$, $q$.} speaking, on $\R^2$ 
\begin{enumerate}
\item when $\d_tu-\Delta u=0$, an initial control on $\Delta^{\ell}u(0,\cdot)$ in $L^q$ yields a $t^{-\left(\tfrac1q-\tfrac1p+\tfrac12\left(|\alpha|-\ell\right)\right)}$ decay for $\d^\alpha u$ in $L^p$;
\item when $\d_t^2u-\Delta u-\Delta\d_tu=0$, an initial control on $\Delta^{\ell} \partial_{t} u(0,\cdot)$ and $\Delta^{\ell} u(0,\cdot)$ in $L^1$  yield a $t^{-\left(\tfrac34- \tfrac32\tfrac1p +\tfrac12\left(|\alpha|-\ell\right)\right)}$ decay for $\d^\alpha u$ in $L^p$ for $p \geq 2$.
\end{enumerate}
The latter bound is classical but non trivial and we refer the reader to \cite{Shibata} for precise statements and proofs. A significant part of the proofs of Theorems~\ref{th:dispersive}, \ref{th:dispersive-localized} and~\ref{th:scalar} is actually, in disguise, an extension of the large-time estimates for constant-coefficient hyperbolic-parabolic operators to classes of systems with no particular structure, including on one hand systems with a scalar-like --- but not scalar --- hyperbolic part and a cross-diffusion, and on the other hand anisotropic systems with dispersion. Concerning the latter, we point out that for the most part of the literature the extra decay due to dispersion is merely overlooked, whereas the remaining body of works is restricted to isotropic systems, most often variations on the wave equation. We also stress that even in the case where the hyperbolic part and the diffusive part commute, combining the well-developed dispersive estimates for the former with dissipative estimates for the latter yields non sharp decay rates; see the related detailed discussion in \cite{Hoff_Zumbrun-NS_compressible_ponctuel}.

\br\label{rk:log}
The $\ln(2+t)$ factors in Theorem~\ref{th:dispersive} and Theorem~\ref{th:scalar} and the lack of optimal decay in Theorem~\ref{th:dispersive} are due to the limited smoothness we assume and to the way we use it. Assuming $E_{0}$ small enough and that $E_{0,3}<\infty$, we also prove that, in Case~\case{a} of Assumptions~\cond1-\cond2, 
\begin{align*}
\|\bW(t,\cdot-\bfphi(t,\cdot))
&-\ubU\|_{L^p}
+\|\nabla\bfphi(t,\cdot)\|_{L^p}
+\|\d_t\bfphi(t,\cdot)\|_{L^p}\\
&+\|\bcK(t,\cdot)-\ubK\|_{L^p}
+\|\d_t\bPsi(t,\cdot)-\ubOm\|_{L^p}
\,\leq \frac{C_{p_{0}}\,E_{0,3}}{(1+t)^{\frac{3}{4}-\frac32\frac{1}{p}}}\,,
&2<p_0\leq p\leq \infty\,,
\end{align*}
\begin{align*}
\|\bcW(t,\cdot)-\bU^{\bcK(t,\cdot)}(\bPsi(t,\cdot))\|_{L^{p}}
&\leq \,
\begin{cases}\ds
\frac{C_{p_0,p_{1}}\,E_{0,3}}{(1+t)^{\frac34-\frac{3}{2}\frac{1}{p}+\frac12}}\,,&
2 < p_0 \leq p \leq p_{1} < \infty\,,\\[1em]\ds
\frac{C \,E_{0,3}\,(1%E_{0}
+E_{0,3})}{(1+t)^{\frac34+\frac12}}\,,&
p = \infty\,,
\end{cases}
\end{align*}
and, in Case~\case{b} of Assumptions~\cond1-\cond2, for any $p_{1} \in [4,\infty)$,
\begin{align*}
\|\bcW(t,\cdot)-\bU^{\bcK(t,\cdot)}(\bPsi(t,\cdot))\|_{L^{p}}
&\leq \,
\begin{cases}\ds
\frac{C_{p_0}\,E_{0,3}}{(1+t)^{1-\frac{1}{p}}}\ \,,&
2<p_0\leq p\leq 4\,,\\[1em]\ds
\frac{C_{p_1}\,E_{0,3}}{(1+t)^{1-\frac{1}{\min(\{p,p_1\})}}}\,\,,&
4\leq p\leq\infty\,.
\end{cases}
\end{align*}
We prove this remark in Subsection \ref{ss:additional_estim_whitham}.
\er

\subsection{Averaged systems}

In the foregoing subsection, to provide educated guesses on expected decay rates, we have already largely anticipated that the dynamics of local wavevectors $\bcK(t,\cdot)$ obeys at leading-order an hyperbolic-parabolic system, in the neighborhood of the constant $\ubK$. We make such a claim precise here. Combining this with results on modulation behavior will complete the leading-order description of the dynamics near the periodic wave of profile $\ubU$.

We shall compare $\bPsi$ and $\bcK$ with respectively $\bPsiW$ and $\bcKW=\nabla\bPsiW$ solving an equation of the form
\be\label{e:W-phase-intro}
\d_t\bPsiW=\bOm(\nabla\bPsiW)+\ubLambda_0[\nabla](\nabla\bPsiW)
\ee
thus also
\be\label{e:W-wn-intro}
\d_t\bcKW=\nabla(\bOm(\bcKW))+\nabla(\ubLambda_0[\nabla](\bcKW))\,,
\ee
with $\bcKW$ curl-free. In the foregoing systems, both $\ubLambda_0[\nabla](\nabla(\cdot))$ and $\nabla(\ubLambda_0[\nabla](\cdot))$ are elliptic linear operators whose action is encoded by homogeneous second-order Fourier multipliers. These elliptic operators are actually differential operators in Case~\case{b}. Moreover our choice of $\ubLambda_0$ enforces that, in any case, $\ubLambda_0[\nabla](\nabla(\cdot))$ commutes with $\dD_{\bK}\bOm(\ubK)(\nabla(\cdot))$. We precisely define this operator in Appendix \ref{ss:artificial}.

We provide a partly formal derivation of \eqref{e:W-phase-intro}/\eqref{e:W-wn-intro} in Appendix~\ref{ss:formal}. This derivation is dramatically more involved than the corresponding one for one-dimensional reaction-diffusion systems --- for which we refer for instance to \cite[Section~4.3]{DSSS} ---, and even trickier than the general one-dimensional case analyzed in \cite{JNRZ-conservation}. As in \cite{JNRZ-conservation}, there are two separate steps in our derivation. The first one, carried out in Appendix~\ref{ss:formal}, is purely formal and is an adaptation of the strategy tailored in \cite{Noble-Rodrigues} and subsequently used in \cite{JNRZ-conservation}. We insert a suitable geometrical optics' \emph{ansatz} in \eqref{e:rd-intro}, identify a few orders of the formal expansion and group together some of the equations to obtain differential systems similar to \eqref{e:W-phase-intro}/\eqref{e:W-wn-intro}; see \eqref{e:W-phase-formal}/\eqref{e:W-wn-formal}. Alternatively, this first step could be replaced with a spectrally motivated derivation; see Remark~\ref{rk:formal-spectral}. At the level of local wavevectors, both the geometrical optics' derivation expounded here and the alternative spectral derivation hinge on low-frequency expansions so that the only piece of information on their structure, inherited from \cond{2}, is that their first-order part is hyperbolic and that they are diffusive in the low-frequency regime. Unfortunately, except in the scalar Subcase~\case{b0}, this is insufficient to deduce that a second-order system is well-posed; see in particular the concrete example given in Appendix~\ref{ss:artificial} to illustrate that this may fail even in Case~\case{b}. Our second step of the derivation finds a canonical way to replace the system obtained in the first step with a well-posed parabolic system sharing, at leading-order, the same low-frequency properties. It is in this step that, in Case~\case{a}, we need to leave the differential frame for the Fourier-multiplier class. This second part of the derivation is analytical and parallel but significantly harder to the analysis in \cite[Appendix~B.2]{JNRZ-conservation}. We stress that in this second step the question to solve has a much broader significance than the study of the averaged dynamics near periodic traveling waves ; it is a general question about the dynamics of second-order systems near constant states, relevant even when the original system is well-posed, and dissipative in some hypocoercive sense but not genuinely parabolic. In particular, our analysis in Appendix~\ref{ss:artificial} extends in various ways, including the class of systems considered and the sharpness of estimates proved, the analysis about \emph{artificial viscosity systems} in \cite{Hoff_Zumbrun-NS_compressible_pres_de_zero,Rodrigues-compressible} (discussed further in \cite{Rodrigues-these} and \cite[Appendix~A]{R}).

Though, for the sake of brevity, we shall not dwell on this line of investigation, we mention that, in the spirit of Theorem~\ref{th:mod-behavior}, and with the same kind of shortcomings, we could validate \eqref{e:W-phase-intro}/\eqref{e:W-wn-intro} without specializing to either Case~\case{a} or Case~\case{b}. 

\bt[Whitham equation]\label{th:whitham}
Assume that we are under the assumptions of Theorem~\ref{th:nonlinear_stab} with notational conventions of Theorem~\ref{th:mod-behavior}.\\ %Assume that $\bW_{0}(\cdot-\bfphi_{0})-\ubU \in L^{\frac32}(\R^{2})$ in Case~\case{a}, that $E_{0}$ is small enough and let $\eta>0$ in Case~\case{b} but when Subcase~\case{b0} fails. 
Let $\bPsiW$ and $\bcKW$ satisfy respectively \eqref{e:W-phase-intro} and \eqref{e:W-wn-intro} with initial data
\begin{align*}
\bPsiW(0,\bfx)&=(\Id-\bfphi_{0})^{-1}\left(\transp{\ubK} \bfx \right)\,,&
\bcKW(0,\cdot)&= \nabla \bPsiW(0,\cdot)\,. 
\end{align*}
Then for any $t \geq 0$, if Case~\case{a} and $\bW_{0}(\cdot-\bfphi_{0})-\ubU \in L^{\frac32}(\R^{2})$ hold
\begin{align*}\ds
\|\bcK(t,\cdot)  - \bcKW(t,\cdot) \|_{L^{p}} \leq
\frac{C\left(E_{0,3}\,(1+E_{0,3}) + \| \bW_{0}(\cdot-\bfphi_{0})-\ubU \|_{L^{\frac32}}\right)}{(1+t)^{\frac34-\frac32\,\frac1p + \frac12}}\,, \quad
2 \leq p \leq  \infty, 
&\textrm{ Case~\case{a}},
\end{align*}
\begin{align*}\ds
\|\bPsiW(t,\cdot)  - \bPsi(t,\cdot) \|_{L^{p}} \leq
\frac{C\left(E_{0,3}\,(1+E_{0,3}) + \| \bW_{0}(\cdot-\bfphi_{0})-\ubU \|_{L^{\frac32}}\right)}{(1+t)^{\frac34-\frac32\,\frac1p}}\,, \quad
2 < p_{0} \leq p \leq  \infty, 
&\textrm{ Case~\case{a}},
\end{align*}
whereas if Case~\case{b} holds
\begin{align*}
\|\bcK(t,\cdot)  - \bcKW(t,\cdot) \|_{L^{p}} \leq\begin{cases}\ds
\frac{C E_{0,3}}{(1+t)^{\frac12-\frac1p + \frac12}}\,,
&2 \leq p \leq \infty,\textrm{ Subcase~\case{b0}},\\[1em]\ds
\frac{C_{\eta} (E_{0,3} + E_{0,3}^2)}{(1+t)^{\frac12-\frac1p + \frac12 (1-\frac1p) - \eta}} \,, 
&2 \leq p \leq \infty,\,\eta>0,\textrm{ Case~\case{b} but \case{b0} fails}\,,
\end{cases}
\end{align*}
\begin{align*}
\| \bPsiW(t,\cdot)  - \bPsi(t,\cdot) \|_{L^{p}} \leq\begin{cases}\ds
\frac{C_{p_{0}} E_{0,3}}{(1+t)^{\frac12-\frac1p}} \,, 
&2 < p_{0} \leq p \leq \infty,\textrm{ Subcase~\case{b0}},\\[1em]\ds
\frac{C_{\eta,p_{0}} (E_{0,3} + E_{0,3}^2)}{(1+t)^{\frac12 -\frac1p - \frac12 \frac1p - \eta}} \,, 
&2 < p_{0} \leq p \leq \infty,\,\eta>0,\textrm{ Case~\case{b} but \case{b0} fails}.
\end{cases}
\end{align*}
\et

\br 
The extra localization of $\bW_{0}(\cdot-\bfphi_{0})-\ubU$ in Case~\case{a} is crucially used to obtain the optimal decay. Without this assumption, one only gets
\[
\|\bcK(t,\cdot)  - \bcKW(t,\cdot) \|_{L^{p}} \leq
\frac{C E_{0,3}}{(1+t)^{1 - \frac1p}}\,, \qquad
2 \leq p \leq \infty.
\]
\er

Actually in the dispersive case~\case{a}, the nonlinear terms are subcritical so that it is enough to retain from \eqref{e:W-phase-intro} and \eqref{e:W-wn-intro} their linear approximants 
\be\label{e:W-phase-linear-intro}
\d_t\bPsiW_{lin}\,=\,\ubOm+\dD_{\bK}\bOm(\ubK)(\nabla \bPsiW_{lin}-\ubK) + \ubLambda_0[\nabla](\nabla \bPsiW_{lin}-\ubK)\,,
\ee
and
\be\label{e:W-wn-linear-intro}
\d_t\bcKW_{lin}=\nabla(\dD_{\bK}\bOm(\ubK)(\bcKW_{lin}-\ubK))+\nabla(\ubLambda_0[\nabla](\bcKW_{lin}-\ubK))\,.
\ee
Obviously the same hope holds for Case~\case{b} when subcritical perturbations are considered. Indeed we prove the following results.

\bt[Dispersive case, linearized Whitham equation]\label{th:dispersive-linear-whitham}
Under the Case~\case{a} assumptions of Theorem~\ref{th:nonlinear_stab} and with notational conventions of Theorem~\ref{th:mod-behavior}, let $\bPsiW_{lin}$ and $\bcKW_{lin}$ satisfy respectively \eqref{e:W-phase-linear-intro} and \eqref{e:W-wn-linear-intro} with initial data 
\begin{align*}
\bPsiW(0,\bfx)&=(\Id-\bfphi_{0})^{-1}\left(\transp{\ubK} \bfx \right)\,,&
\bcKW(0,\cdot)&= \nabla \bPsiW(0,\cdot)\,. 
\end{align*}
Then, for any $t \geq 0$,
\begin{align*}
\| \bcK(t,\cdot) - \bcKW_{lin}(t,\cdot) \|_{L^{p}} &\ds
\leq \frac{C\,E_{0,3}}{(1+t)^{\frac34-\frac32 \frac1p + \frac12 \frac1p}} \,,\qquad
&2 \leq p \leq \infty\,,\\
\| \bPsi(t,\cdot) - \bPsiW_{lin}(t,\cdot) \|_{L^{p}} &\ds
\leq \frac{C_{p_{0}} \,E_{0,3}}{(1+t)^{\frac14-\frac32 \frac1p + \frac12 \frac1p}} \,,\qquad
&2 < p_{0} \leq p \leq \infty\,.
\end{align*}
\et

\bt[Subcritical perturbations, linearized Whitham equation]\label{th:subcritical-linear-whitham}
Under the assumptions of Theorem~\ref{th:nonlinear-localized} and in Case~\case{a} of Theorem~\ref{th:dispersive-localized} and with notational conventions of Theorem~\ref{th:mod-behavior}, let $\bPsiW_{lin}$ and $\bcKW_{lin}$ satisfy respectively \eqref{e:W-phase-linear-intro} and \eqref{e:W-wn-linear-intro} with initial data 
\begin{align*}
\bPsiW(0,\bfx)&=(\Id-\bfphi_{0})^{-1}\left(\transp{\ubK} \bfx \right)\,,&
\bcKW(0,\cdot)&= \nabla \bPsiW(0,\cdot)\,. 
\end{align*}
Then, for any $t \geq 0$,
\begin{align*}
\|\bcW(t,\cdot)-\bU^{\bcK(t,\cdot)}(\bPsi(t,\cdot))\|_{L^{p}} \leq
\begin{cases}\ds
\frac{C \cE_{0}}{(1+t)^{\frac54-\frac32\,\frac1p + \frac12}}\,,
&2 \leq p \leq \infty\textrm{, Case~\case{a}},\\[1em]\ds
\frac{C \cE_{0}}{(1+t)^{1-\frac1p + \frac12}} \,,
&2 \leq p \leq \infty \textrm{, Case~\case{b}}\,,
\end{cases}
\end{align*}
and
\begin{align*}
\|\bcK(t,\cdot)  - \bcKW_{lin}(t,\cdot) \|_{L^{p}} \leq
\begin{cases}\ds
\frac{C \cE_{0} \ln(2+t)}{(1+t)^{\frac54-\frac32\,\frac1p + \frac12}}\,,
&2 \leq p \leq \infty \textrm{, Case~\case{a}},\\[1em]\ds
\frac{C \cE_{0} \ln(2+t)}{(1+t)^{1-\frac1p + \frac12}} \,,
&2 \leq p \leq \infty \textrm{, Case~\case{b0}},\\[1em]\ds
\frac{C_{p_{1}} \cE_{0}}{(1+t)^{1-\frac1p + \frac12 (1-\frac1p)}} \,, 
&2 \leq p \leq p_{1} < \infty\textrm{, Case~\case{b} but \case{b0} fails},\\[1em]\ds
\frac{C \cE_{0} \ln(2+t)}{(1+t)^{1-\frac1p + \frac12 (1-\frac1p)}} \,,
&p = \infty \textrm{, Case~\case{b} but \case{b0} fails},
\end{cases}
\end{align*}
\begin{align*}
\|\bPsi(t,\cdot)  - \bPsiW_{lin}(t,\cdot) \|_{L^{p}} \leq
\begin{cases}\ds
\frac{C \cE_{0} \ln(2+t)}{(1+t)^{\frac34-\frac32\,\frac1p + \frac12}}\,,
&2 \leq p \leq \infty \textrm{, Case~\case{a}},\\[1em]\ds
\frac{C \cE_{0} \ln(2+t)}{(1+t)^{\frac12-\frac1p + \frac12}} \,,
&2 \leq p \leq \infty\textrm{, Case~\case{b0}},\\[1em]\ds
\frac{C_{p_{1}} \cE_{0}}{(1+t)^{\frac12-\frac1p + \frac12 (1-\frac1p)}} \,,
&2 \leq p \leq p_{1} < \infty\textrm{, Case~\case{b} but Subcase~\case{b0} fails},\\[1em]\ds
\frac{C \cE_{0} \ln(2+t)}{(1+t)^{\frac12-\frac1p + \frac12 (1-\frac1p)}} \,,
&p = \infty \textrm{, Case~\case{b} but Subcase~\case{b0} fails}.
\end{cases}
\end{align*}
\et
Again, as in Remark \ref{rk:log}, one can remove almost all the $\ln(2+t)$ in the previous bounds by assuming more regularity. There is one exception for $p=\infty$ in  Case ~\case{b} (but \case{b0} fails) since then the logarithm factor already occurs at the linear level, through Proposition~\ref{pr:ubL0-scalar}.

\subsection{Perspectives, outline and notation}\label{s:notation}

$ $

\smallskip

\noindent\emph{Perspectives.} The present contribution offers an essentially complete analysis of the nonlinear dynamics near spectrally stable periodic waves of parabolic systems without conservation laws. Most natural follow-up questions are two-fold. On one hand in order to apply the present results it is important to provide stability diagrams in all relevant bifurcation scenari, multiplying the type of analysis carried out in \cite{BdR-R}. On the other hand it is equally important to enlarge the class of systems encompassed by our nonlinear analysis so as to allow conservation laws, in the same way as \cite{JNRZ-conservation} extends \cite{JNRZ-RD1,JNRZ-RD2}. This is crucial to be able to deal with most hydrodynamic applications.

\smallskip

\noindent\emph{Outline.} The organization of the remainder of the paper, after the present introduction, reflects the plan of the introduction. The second section proves stability theorems, Theorems~\ref{th:nonlinear_stab} and~\ref{th:nonlinear-localized}. The third section proves modulational-behavior theorems, Theorems~\ref{th:mod-behavior}, \ref{th:dispersive}, \ref{th:dispersive-localized} and~\ref{th:scalar}. The fourth section proves theorems on averaged modulation systems, Theorems~\ref{th:whitham}, \ref{th:dispersive-linear-whitham} and~\ref{th:subcritical-linear-whitham}. The paper is concluded with four appendices, devoted respectively to 
\begin{itemize}
\item elements of Bloch-wave spectral analysis, 
\item the geometric structure of profile equations, 
\item phase estimates, mostly showing how to bound $\nabla\bfphi$ with $\Delta \bfphi$, 
\item geometric optics as needed in the direct derivation of averaged systems.
\end{itemize}

\smallskip

% perp, L^2_perp, espaces fonctionnels

\noindent\emph{Notation.} We conclude this introduction by collecting some elements of our notational conventions. 

When $a \in \R$, $a_{+}:=\max(a,0)$. When $\bfx=(x_1,x_2)$, $\bfx^\perp:=(-x_2,x_1)$.

When $a$ and $b$ are two elements of the same set, $\delta_{a,b}$ is $1$ if $a=b$, $0$ otherwise. 

We often identify vectors of $\R^{n}$ with column vectors, elements of $\cM_{n,1}(\R)$. When $A$ and $B$ are linear operators, $[A,B]:=AB-BA$. Our complex inner scalar products are skew-linear in their first argument, linear in their second argument. We denote the canonical basis of $\R^2$ as $(\beD_1,\beD_2)$.

For a map $\bF : \R^{p} \to \R^{n}$ and $\bfx_{0} \in \R^{p}$, we denote by $\dD \bF(\bfx_{0})(\cdot)$ the differential of $\bF$ at $\bfx_{0}$, a linear map from $\R^p$ to $\R^n$, by $\dD^2 \bF(\bfx_{0})(\cdot,\cdot)$ the second differential of $\bF$ at $\bfx_{0}$, a bilinear map from $\R^p\times\R^p$ to $\R^n$, and by $\nabla \bF (\bfx_0)$ the gradient of $\bF$ at $\bfx_0$, the transpose of the Jacobian at $\bfx_0$, an element of $\cM_{p,n}$ defined by $(\nabla \bF(\bfx_0))_{j,\ell} = \partial_{j} \bF_{\ell}(\bfx_0)$ (with standard notation for partial derivative and coordinate). We add a subscript when using partial differential operators, such as $\dD_{\bfx} \bH(\bfx_{0},\bfy_{0})(\cdot)$ to denote the partial differential of $\bH$ with respect to the $\bfx$-variable at $(\bfx_{0},\bfy_{0})$ when $\bH :\R^{p} \times \R^{p} \to\R^{n}$, $(\bfx,\bfy) \mapsto\bH(\bfx,\bfy)$.

We make two main exceptions to the previous differential notation by changing at some specific places where we mark the evaluation point $\bfx_0$. Explicitly, in Section \ref{s:mod}, we denote by $\dD_{\bfxi}\bfq^{\b0}_{j}(\cdot)$, $\dD_{\bfxi}\tbq^{\b0}_{\ell}(\cdot)$ and $\dD_{\bfxi}\bD_{\b0}(\cdot)$ the differential at $\bf0$ of respectively $\bfxi \mapsto \bfq_{j}^{\bfxi}$, $\bfxi \mapsto \tbq_{\ell}^{\bfxi}$ and $\bfxi \mapsto \bD_{\bfxi}$ and, in Section \ref{s:averaged_dyn}, we denote by $\dD^{2}_{\bfxi}\bD_{\b0}(\cdot,\cdot)$ the second differential at $\bf0$ of $\bfxi \mapsto \bD_{\bfxi}$.

Divergence $\transp{\nabla}$ and Laplacian operators $\Delta$ are always taken with respect to the spatial variable only, and we do not mark this partial restriction. We extend this omission to gradients $\nabla$ when there is no risk of confusion. With respect to the spatial variable they are defined as follows: if $\bG : \R^{2} \to \cM_{2,n}(\R)$ then $\transp{\nabla} \bG(\bfx_0)$ is the vector of $\R^n$ given by $(\transp{\nabla} \bG(\bfx_0))_{\ell} = \sum_{j=1}^2 \partial_{j} \bG_{j,\ell}(\bfx_0)$; if $\bF: \R^{2} \to \R^n$ then $\Delta \bF(\bfx_0):=\transp{\nabla}(\nabla\bF)(\bfx_0)$ is the vector of $\R^n$ given by $(\Delta \bF(\bfx_0))_{\ell} = \sum_{j=1}^2 \partial_{j}^2 \bF_{\ell}(\bfx_0)$.

In Section \ref{ss:phase-in}, when carrying out more abstract algebraic computations with too many spatial differential operators already involved and no particular functional topology in mind, we switch from differential notation $\dD$ to linearized notation $\bL$. We use $\bL$ essentially as $\dD$ with respect to where we mark evaluation points, directions of application, restrictions, \emph{etc.} 

We identify spaces of $(\beD_1,\beD_2)$-periodic functions with closed subspaces of functions over $[0,1]^2$ satisfying suitable boundary conditions. We use the subscript ${}_{\rm per}$ to distinguish those. At an abstract level they may be defined as the closure for the topology at hand of the restrictions to $[0,1]^2$ of smooth  $(\beD_1,\beD_2)$-periodic functions.

\smallskip

\noindent\emph{Acknowledgment.} The authors thank Kevin Zumbrun for his constant interest in the present work. B.M. thanks the University of Rennes for its hospitality.

\section{Stability}

\subsection{Linear estimates}\label{ss:stab-lin}

We begin our stability analysis with linear estimates. In this part we make extensive use of the background material provided in Appendix~\ref{s:spec}.

We consider $(S(t))_{t\geq0}$ the semi-group generated by $L$. The linear counterpart to Theorem~\ref{th:nonlinear_stab} is that given some $\bfg_0=(\bfphi_0\cdot \nabla) \ubU+\bV_0$ with initial estimates on $(\bV_0,\Delta\bfphi_0)$, one is able to split $S(t)[\bfg_0]$ as 
\[
S(t)[\bfg_0]\,=\,(\bfphi(t,\cdot)\cdot \nabla) \ubU+\bV(t,\cdot)
\]
so as to ensure large-time decay estimates on $(\bV(t,\cdot),\nabla\bfphi(t,\cdot),\d_t\bfphi(t,\cdot))$. Note that since $L$ is a parabolic operator, the fact that it does generate an analytic semigroup and the accompanying short-time estimates are part of the classical theory for linear PDEs, for which we refer to \cite{Pazy}.

Condition \cond2 ensures that we may decompose $S(t)$ according to
\[
S(t)[\bfg]\,=\,(s(t)[\bfg] \cdot \nabla) \ubU\,+\,S_{1}(t)[\bfg] + S_{2}(t)[\bfg]
\]
with
\[
(s(t)[\bfg])(\bfx):= \int_{[-\pi,\pi]^{2}} \chi(\bfxi) \eD^{\iD \bfx \cdot \bfxi} 
\eD^{t \bD_{\bfxi}} 
\bp \left\langle\tbq_{1}^{\bfxi};\check{\bfg}(\bfxi,\cdot) \right\rangle_{L^{2}_{per}} \\ 
\left\langle\tbq_{2}^{\bfxi};\check{\bfg}(\bfxi,\cdot) \right\rangle_{L^{2}_{per}}  \ep \dD\bfxi\,,
\]
and
\begin{align*}
(S_{1}(t)[\bfg])(\bfx)
&:= \int_{[-\pi,\pi]^{2}} (1-\chi(\bfxi)) \eD^{\iD \bfx \cdot \bfxi} \eD^{t\,L_{\bfxi}} (\check{\bfg}(\bfxi,\cdot))(\bfx) \dD\bfxi 
+\int_{[-\pi,\pi]^{2}}  \chi(\bfxi) \eD^{\iD \bfx \cdot \bfxi} \eD^{t\,L_{\bfxi}} (\I - \Pi_{\bfxi})(\check{\bfg}(\bfxi,\cdot))(\bfx) \dD\bfxi\\
(S_{2}(t)[\bfg])(\bfx)
&:= \int_{[-\pi,\pi]^{2}} \chi(\bfxi) \eD^{\iD \bfx \cdot \bfxi}  
\bp \bfq_{1}^{\bfxi}(\bfx) - \bfq_{1}^{\b0}(\bfx) & \bfq_{2}^{\bfxi}(\bfx) -  \bfq_{2}^{\b0}(\bfx) \ep 
\eD^{t\,\bD_{\bfxi}} 
\bp \left\langle \tbq_{1}^{\bfxi};\check{\bfg}(\bfxi,\cdot) \right\rangle_{L^{2}_{per}}\\ 
\left\langle \tbq_{2}^{\bfxi};\check{\bfg}(\bfxi,\cdot)\right\rangle_{L^{2}_{per}} \ep \dD\bfxi
\end{align*}
where $\check{g}$ denotes the Bloch transform of $g$,
\begin{itemize}
\item $\chi$ is a smooth function valued in $[0,1]$, compactly supported in a sufficiently small neighborhood of $\b0$ and equal to $1$ in a (smaller) neighborhood of $\b0$;
\item $\Pi_\bfxi$ is the spectral projector of $L_\bfxi$ associated with its spectrum near the origin;
\item $(\bfq_{1}^{\bfxi},\bfq_{2}^{\bfxi})$ is a basis of the range of $\Pi_\bfxi$, smoothly dependent on $\bfxi$, such that $(\bfq_{1}^{\b0},\bfq_{2}^{\b0})=(\d_1\ubU,\d_2\ubU)$;
\item $(\tbq_{1}^{\bfxi},\tbq_{2}^{\bfxi})$ is a basis of the range of $\Pi_\bfxi^*$, smoothly dependent on $\bfxi$ and in duality with $(\bfq_{1}^{\bfxi},\bfq_{2}^{\bfxi})$;
\item $\bD_\bfxi$ is the matrix of the restriction of $L_\bfxi$ to the range of $\Pi_\bfxi$ in the basis $(\bfq_{1}^{\bfxi},\bfq_{2}^{\bfxi})$.
\end{itemize}
See Appendix \ref{ss:pert} for more details about the construction of such objects. For concreteness, we introduce $\xi_0$ a sufficiently small positive number measuring the support of $\chi$ in the sense that on the support of $\chi$, $\|\bfxi\|\leq \xi_0$.

\br\label{rk:choice}
Note that $(\bfq_{1}^{\bfxi},\bfq_{2}^{\bfxi})$ is not uniquely determined by the above conditions so that there is some flexibility in the definition of $s(t)$. This flexibility is irrelevant in the stability analysis but it will be used to impose further normalization in the asymptotic behavior part. 
\er

Throughout this work, we will use a consequence of \cond{1}-\cond{2}, related to \cond{1''} in Appendix~\ref{ss:diff}, that asserts that there exists $C_{1}>0$ such that for any $\bfxi\in[-\pi,\pi]^2$  satisfying $\|\bfxi\|\leq \xi_0$ and any $t\geq0$,
\be\label{eq:decay_D}
\lnor \eD^{t\,\bD_{\bfxi}} \rnor\,\leq\,C_{1} e^{-\theta t\|\bfxi\|^2}\,.
\ee
We first provide bounds adapted to localized perturbations.

\bpr\label{pr:lin-local}
Assume \cond1-\cond2.
\begin{enumerate}
\item There exists $\theta'>0$, such that, for any $(s,s')\in(\R_+)^2$ such that $s'\leq s$, there exists $C_{s',s}$ such that for any $t> 0$
\begin{align*}
\|S_{1}(t)[\bfg]\|_{H^{s}}
&\leq \frac{C_{s',s}}{(\min(\{1,t\}))^{\frac{(s-s')}{2}}}\,\eD^{-\theta'\,t}
\|\bfg\|_{H^{s'}}\,.
\end{align*}
\item For any $s\in\R_+$ and any $\beta\in\N^2$, there exists $C_{s,\beta}$ such that for any $1\leq q \leq 2 \leq p \leq +\infty$, and any $t\geq0$
\begin{align*}
\|S_{2}(t)[\d_\bfx^\beta\bfg]\|_{W^{s,p}}
&\leq \frac{C_{s,\beta}}{(1+t)^{\frac{1}{2}+\frac{1}{q}-\frac{1}{p}}}\,\|\bfg\|_{L^{q}}\,.
\end{align*}
\item For any $\alpha\in\N^2$, any $\beta\in\N^2$ and any $\ell\in\N$, there exists $C_{\alpha,\ell,\beta}$ such that for any $1\leq q \leq 2 \leq p \leq +\infty$, and any $t\geq0$
\begin{align*}
\|\,\d_\bfx^\alpha\,\d_t^\ell\,s(t)[\d_\bfx^\beta\bfg]\|_{L^p}
&\leq \frac{C_{\alpha,\ell,\beta}}{(1+t)^{\frac{|\alpha|+\ell}{2}+\frac{1}{q}-\frac{1}{p}}}\,\|\bfg\|_{L^{q}}\,.
\end{align*}
\end{enumerate}
\epr

\begin{proof}
To prove the first point, it is sufficient to combine an $H^{s'}\to H^s$ bound for $0<t\leq 1$ with an $H^s\to H^s$ bound for $t\geq 0$. Moreover the former follows from the parabolicity of $L$ (combined with bounds on $s(t)$ and $S_2(t)$ proved below). In turn, the latter may be derived, through Parseval's identity, from
\begin{align*}
\lnor \eD^{t\,L_{\bfxi}} \rnor_{H^s_{\rm per}\to H^s_{\rm per}}&\,\leq\,C e^{-\theta'\, t}\,,&
\|\bfxi\|\geq \xi_0\,,\\ 
\lnor \eD^{t\,L_{\bfxi}}\,(\I-\Pi_\bfxi)\rnor_{H^s_{\rm per}\to H^s_{\rm per}}&
\,\leq\,C e^{-\theta'\, t}\,,&
\|\bfxi\|\leq \xi_0\,\,.
\end{align*}
As pointed out in Appendix~\ref{ss:diff}, these bounds stem from condition \cond1-\cond2.

To prove the second point, by using \eqref{eq:decay_D} and integration by parts in scalar products, from Hausdorff-Young and H\"older inequalities one derives 
\begin{align*}
\|S_{2}(t) [\bfg]  \|_{W^{s,p}} 
&\lesssim \left\|\bfxi\mapsto \|\bfxi\|\,\eD^{-\theta t \|\bfxi\|^{2}} \|\check{\bfg}(\bfxi,\cdot)\|_{L^{q}_{\rm per}}\right\|_{L^{p'}_\bfxi}\\
&\lesssim  \left\|\bfxi\mapsto \|\bfxi\|\,\eD^{-\theta t \|\bfxi\|^{2}}\right\|_{L^{r}_\bfxi}\ \times\ \|\check{\bfg} \|_{L^{q'}_{\bfxi} L^{q}_{\rm per}}\\
&\lesssim (1+t)^{-\left(\frac{1}{2}+\frac{1}{r}\right)}\,\|g\|_{L^q}
\end{align*}
with $p'$, $q'$ Lebesgue conjugate respectively to $p$ and $q$, and $1/r=1/p'-1/q'=1/q-1/p$. Hence the second bound.

To prove the third point, note that
\[
(\d_\bfx^\alpha\,\d_t^\ell s(t)[\d_\bfx^\beta\bfg])(\bfx) = \int_{[-\pi,\pi]^{2}} \chi(\bfxi) \eD^{\iD \bfx \cdot \bfxi} 
\,(\iD\bfxi)^\alpha\,(\bD_\bfxi)^\ell\,\eD^{t \bD_{\bfxi}} 
\bp \left\langle\tbq_{1}^{\bfxi};(\d_\bfx+\iD\bfxi)^\beta\check{\bfg}(\bfxi,\cdot) \right\rangle_{L^{2}_{per}} \\ 
\left\langle\tbq_{2}^{\bfxi};(\d_\bfx+\iD\bfxi)^\beta\check{\bfg}(\bfxi,\cdot) \right\rangle_{L^{2}_{per}}  \ep \dD\bfxi\,.
\]
From here the third bound is proved essentially as was the second one, the decay stemming from
\[
\left\|\bfxi\mapsto \|\bfxi\|^{|\alpha|+\ell}\,\eD^{-\theta t \|\bfxi\|^{2}}\right\|_{L^{r}_\bfxi}
\lesssim (1+t)^{-\left(\frac{|\alpha|+\ell}{2}+\frac{1}{r}\right)}\,.
\]
\end{proof}

We now focus on bounds adapted to initial data given as phase modulations. Throughout we make regular use of bounds from Appendix~\ref{s:phase-like} and we implicitly assume that $\bfphi$ has no affine component at $\infty$, in the sense that $\bfphi=\Delta^{-1}\Delta\bfphi$. Consistently, the phases built with $s(t)$ also satisfy the latter condition.

\bpr\label{pr:lin-mod}
Assume \cond1-\cond2.
\begin{enumerate}
\item There exists $\theta'>0$, such that, for any $(s,s')\in\R_+$ such that $s'+1\leq s$ and any $(p_0,p_1)$ such that $2<p_0<p_1<\infty$, there exists $C_{p_0,p_1,s,s'}$ such that for any $t>0$, and any $p_0\leq p\leq p_1$,
\begin{align*}
\|S_{1}(t)[(\bfphi \cdot \nabla) \ubU]\|_{W^{s,p}}
&\leq \frac{C_{p_0,p_1,s,s'}}{(\min(\{1,t\}))^{\frac{(s-(s'+1))}{2}}}\,\eD^{-\theta'\,t}
\,\|\Delta\bfphi\|_{L^1\cap H^{s'}}\,.
\end{align*}
\item For any $s\in\R_+$ and any $p_0>2$, there exists $C_{p_0,s}$ such that for any $p_0 \leq p \leq +\infty$, and any $t\geq0$
\begin{align*}
\|S_{2}(t)[(\bfphi \cdot \nabla) \ubU]\|_{W^{s,p}}
&\leq \frac{C_{p_0,s}}{(1+t)^{\frac{1}{2}-\frac{1}{p}}}\,\|\Delta\bfphi\|_{L^1}\,.
\end{align*}
\item For any $\alpha\in\N^2$, any $\ell\in\N$ and any $2 \leq p \leq +\infty$ such that $|\alpha|+\ell-\tfrac2p>0$, there exists $C_{p,\alpha,\ell}$ such that for any $t\geq0$
\begin{align*}
\|\,\d_\bfx^\alpha\,\d_t^\ell\,s(t)[(\bfphi \cdot \nabla) \ubU]\|_{L^p}
&\leq \frac{C_{p,\alpha,\ell}}{(1+t)^{\frac{|\alpha|+\ell}{2}-\frac{1}{p}}}\,\|\Delta\bfphi\|_{L^1}\,.
\end{align*}
\end{enumerate}
\epr

To ease comparisons with bounds of Proposition~\ref{pr:lin-local}, we recall that $\|\Delta\bfphi\|_{L^1\cap H^{s}}$ should be thought as a relaxed version of $\|\nabla\bfphi\|_{H^{s+1}}$. Note moreover that the condition $|\alpha|+\ell-\tfrac2p>0$ may be written more explicitly as $|\alpha|+\ell\geq 2$ or ($|\alpha|+\ell=1$ and $p>2$).

\begin{proof}
To establish various bounds it is convenient to single out the low-frequency part of $\bfphi$, according to
\begin{align*}
\bfphi&\,=\,\bfphi_{LF}+\bfphi_{HF}\,,&
\widehat{(\bfphi_{LF})}&=\chi \widehat{\bfphi}\,.
\end{align*} 
The contribution of $\bfphi_{HF}$ to the first bound may be deduced from the corresponding estimate in Proposition~\ref{pr:lin-local}. Indeed, since $2\leq p<\infty$,
\begin{align*}
\|S_{1}(t)[(\bfphi_{HF} \cdot \nabla) \ubU]\|_{W^{s,p}}
&\lesssim \|S_{1}(t)[(\bfphi_{HF} \cdot \nabla) \ubU]\|_{H^{s+1}}\,,\\
\|(\bfphi_{HF} \cdot \nabla) \ubU\|_{H^{s'+2}}
&\lesssim\|\bfphi_{HF}\|_{H^{s'+2}}
\lesssim\|\Delta\bfphi\|_{H^{s'}}\,.
\end{align*}
The analysis of the contribution of $\bfphi_{LF}$ requires more care. To begin with, we recall that 
\begin{align}\label{eq:low-phi}
\widecheck{((\bfphi_{LF} \cdot \nabla) \ubU)}(\bfxi,\bfx) = 
(\widehat{(\bfphi_{LF})}(\bfxi) \cdot \nabla) \ubU(\bfx)
\end{align}
and observe that this may be used to gain an extra $\|\bfxi\|$-factor in the second part of the definition of $S_1$ through
\[
(\I - \Pi_{\bfxi})(\widecheck{((\bfphi_{LF} \cdot \nabla) \ubU)}(\bfxi,\cdot))
\,=\,(\I - \Pi_{\bfxi})(\Pi_{\b0} - \Pi_{\bfxi})
(\widecheck{((\bfphi_{LF} \cdot \nabla) \ubU)}(\bfxi,\cdot))
\] 
since $\d_1\ubU$ and $\d_2\ubU$ lie in the range of $\Pi_0$ and $(\I - \Pi_{\bfxi})$ is a projector. Moreover, an extra $\|\bfxi\|$-factor in the first part of $S_1$ is readily obtained from the trivial $(1-\chi(\bfxi))\lesssim \|\bfxi\|$. With this in hands, from Hausdorff-Young inequalities and the embedding $H^{s+1}_{\rm per}\hookrightarrow W^{s,p}_{\rm per}$, one derives
\begin{align*}
\|S_{1}(t)[(\bfphi_{LF} \cdot \nabla) \ubU]\|_{W^{s,p}}
\lesssim
\eD^{-\theta'\,t}\,\,\|\Delta\bfphi_{LF}\|_{L^1}\,
\left\|\bfxi\mapsto \|\bfxi\|^{-1}\right\|_{L^{p'}_\bfxi}
\lesssim
\eD^{-\theta'\,t}\,\,\|\Delta\bfphi\|_{L^1}
\end{align*}
since $p>2$, thus $p'<2$. This achieves the proof of the first bound.

The contribution of $\bfphi_{HF}$ to the second bound may also be deduced from the corresponding estimate in Proposition~\ref{pr:lin-local}. Indeed
\begin{align}\label{estim:high-phi}
\|(\bfphi_{HF} \cdot \nabla) \ubU\|_{L^1}
&\lesssim\|\bfphi_{HF}\|_{L^1}
\lesssim\|\Delta\bfphi\|_{L^1}\,.
\end{align}
The analysis of the contribution of $\bfphi_{LF}$ to the second bound follows from \eqref{eq:low-phi}, Hausdorff-Young inequalities and, since $p>2$,  
\[
\left\|\bfxi\mapsto \|\bfxi\|^{-1}\,\,\eD^{-\theta t \|\bfxi\|^{2}}\right\|_{L^{p'}_\bfxi}
\lesssim (1+t)^{-\left(\frac12-\frac1p\right)}\,.
\]

The third bound is proved similarly. 
\end{proof}

The last set of linear estimates we need to close our nonlinear stability argument consists in short-time bounds. It shall be used to ensure that the nonlinear phase $\bfphi$ does satisfy $\bfphi(0,\cdot)=\bfphi_0$. At the linear level we just need to prove that $s(0)[(\bfphi \cdot \nabla) \ubU]$ is not too far from $\bfphi$. Note that if one relaxes Theorem~\ref{th:nonlinear_stab} by removing the condition $\bfphi(0,\cdot)=\bfphi_0$ from its statement, these short-time bounds become irrelevant. 

\bl\label{l:time-bl}
Assume \cond1-\cond2. For any $\alpha\in\N^2$, any $\ell\in\N$ and any $2 \leq p \leq +\infty$ such that $p>2$ if $|\alpha|=0$ and $p<\infty$ if $|\alpha|\geq2$, there exists $C_{p,\alpha}$ such that for any $t\geq 0$
\begin{align*}
\|\,\d_\bfx^\alpha\,\left(s(t)[(\bfphi \cdot \nabla)\ubU]-\bfphi\right)\|_{L^p}
&\leq C_{p,\alpha}\,\|\Delta\bfphi\|_{L^1\cap W^{(|\alpha|-2)_+,p}}
\,\begin{cases}
\,(1+t)^{\frac12\,\left(1-\left(|\alpha|-\frac2p\right)\right)_+}&\quad\textrm{if } |\alpha|-\tfrac2p\neq 1\\
\ln(2+t)&\quad\textrm{otherwise }
\end{cases}\,.
\end{align*}
\el

\begin{proof}
Since $|\alpha|+1-\tfrac2p>0$, from Proposition~\ref{pr:lin-mod} stems
\begin{align*}
\|\,\d_\bfx^\alpha\,(s(t)-s(0))[(\bfphi \cdot \nabla)\ubU]\|_{L^p}
&\lesssim 
\|\Delta\bfphi\|_{L^1}\,\int_0^t\frac{\dD\tau}{(1+\tau)^{\frac{|\alpha|+1}{2}-\frac{1}{p}}}\,\,,
\end{align*}
which predicts the growth time rates. An examination of the proof of Proposition~\ref{pr:lin-mod} also gives
\begin{align*}
\|\,\d_\bfx^\alpha\,s(t)[(\bfphi_{HF} \cdot \nabla)\ubU]\|_{L^p}
&\lesssim 
\frac{\|\Delta\bfphi\|_{L^1}}{(1+t)^{\frac{|\alpha|}{2}+1-\frac{1}{p}}}
\lesssim 
\,\|\Delta\bfphi\|_{L^1}\,\,.
\end{align*}
Since, moreover the conditions on $p$ ensure
\[
\|\bfphi_{HF}\|_{W^{|\alpha|,p}}
\lesssim \|\Delta\bfphi\|_{L^1\cap W^{(|\alpha|-2)_+,p}}\,,
\]
there only remains to bound $s(0)[(\bfphi_{LF} \cdot \nabla)\ubU]-\bfphi_{LF}$.

Now we observe that
\begin{align*}
\left(s(0)[(\bfphi_{LF} \cdot \nabla) \ubU] - \bfphi_{LF}\right)(\bfx)
&=\int_{[-\pi,\pi]^{2}} \eD^{\iD \bfxi \cdot\bfx}\,(\chi(\bfxi)-1)\,\widehat{(\bfphi_{LF})}(\bfxi)\,\dD\bfxi\\
&\quad+\int_{[-\pi,\pi]^{2}}\eD^{\iD \bfxi \cdot\bfx}\,\chi(\bfxi) \,
\bp \left\langle \tbq_{1}^{\bfxi} - \tbq_{1}^{\b0};(\widehat{(\bfphi_{LF})}(\bfxi) \cdot \nabla) \ubU\right\rangle_{L^{2}_{per}} \\ 
\left\langle \tbq_{2}^{\bfxi} - \tbq_{2}^{\b0} ;(\widehat{(\bfphi_{LF})}(\bfxi) \cdot \nabla) \ubU\right\rangle_{L^{2}_{per}} \ep 
\dD\,\bfxi\,.
\end{align*}
Thus
\[
\|\,\d_\bfx^\alpha\,\left(s(0)[(\bfphi \cdot \nabla)\ubU]-\bfphi\right)\|_{L^p}
\lesssim \|\Delta\bfphi\|_{L^1}
\,\times\,\left\|\bfxi\mapsto \|\bfxi\|^{|\alpha|-1}\right\|_{L^{p'}_\bfxi}\,.
\]
Hence the result (since $p'<2$ when $|\alpha|=0$).
\end{proof}

\subsection{Introducing phases}\label{ss:phase-in}

To carry out our analysis, it is convenient to write equation~\eqref{rd} directly in terms of $\bfphi$ and $\bV$ such that $\bW(t,\cdot-\bfphi(t,\cdot))=\ubU+\bV(t,\cdot)$. When doing so, we have in mind Lemma~\ref{l:var-change} and its variants.

We would like to stress here that the argument is quite robust and to spare unnecessary detailed computations we provide it in abstract form. To do so, we introduce
\begin{align*}
\cP[\bW]
&=\cP(\bW,\nabla\bW,\nabla^2\bW)
\,:=\,
\transp{\nabla}\left(\transp{\ubK}\ubK\,\nabla\bW\right)
+\transp{\nabla}\left(\transp{\ubK}\,
(\bG(\bW)+\ubc\,\transp{\bW})\right)+\bff(\bW)
\end{align*}
and consider its image under a change of variable $\bPhi$
\begin{align}\label{def:bP}
\bP[\tbW,\bPhi]
&=\bP(\tbW,\nabla\tbW,\nabla^2\tbW,\nabla\bPhi,\nabla^2\bPhi)
:=(\cP[\tbW\circ \bPhi^{-1}])\circ\bPhi\\\nn
&\,=\,
|\nabla\bPhi|^{-1}
\transp{\nabla}\left(|\nabla\bPhi|\,\transp{(\ubK\,[\nabla\bPhi]^{-1})}(\ubK\,[\nabla\bPhi]^{-1})\nabla\tbW\right)\\\nn
&\quad+\ |\nabla\bPhi|^{-1}
\transp{\nabla}\left(|\nabla\bPhi|\,\transp{(\ubK\,[\nabla\bPhi]^{-1})}
(\bG(\tbW)+\ubc\,\transp{\tbW})\right)+\bff(\tbW)\,.
\end{align}
At the linear level, the key observation is that
\[
-\bL_{\bPhi}\bP[\ubU,\Id](\bfphi)
\,=\,\bL_\bW\cP [\ubU]((\bfphi\cdot\nabla)\ubU)-(\bfphi\cdot\nabla)(\cP[\ubU])
\] 
thus, since $\cP[\ubU]\equiv 0$,
\begin{equation}\label{eq:key-lin}
\bL_{(\tbW,\bPhi)}\bP[\ubU,\Id](\bV,-\bfphi)
\,=\,\bL_\bW \cP [\ubU](\bV+(\bfphi\cdot\nabla)\ubU)
\,=\,L(\bV+(\bfphi\cdot\nabla)\ubU)\,,
\end{equation}
where $\bL_{\bPhi}$, $\bL_\bW$ and $\bL_{(\tbW,\bPhi)}$ stand for linearized operators. With this in hands, we may rephrase \eqref{rd}.

\bl\label{lem:cancellation-separation}
Let $\bW$ and $(\bfphi,\bV)$ be smooth\footnote{Since the content of the present lemma is essentially algebraic and in the end we only consider classical solutions, we do not make precise assumptions about the level of regularity needed here.} functions such that
\be\label{pertvar} 
\bW(t,\bfx-\bfphi(t,\bfx))\ =\ \ubU(\bfx)+\,\bV(t,\bfx)\,,
\ee
and for any $t$, $\|\nabla\bfphi(t,\cdot)\|_{L^\infty(\R^2)}<1$. Then $\bW$ satisfies \eqref{rd} if and only if $(\bfphi,\bV)$ satisfies
\begin{equation}\label{veq}
\left(\d_t-L\right)(\bV+(\bfphi \cdot \nabla) \ubU)=\cN[\bV,\bfphi]\,,
\end{equation}
or equivalently
\begin{equation}\label{veq0}
\d_t\bV-L_{\bfphi}\bV
\,=\,-(\bfphi_t \cdot [\I_2-\nabla\bfphi]^{-1}\nabla)(\ubU+\bV)
+\bP[\ubU,\Id-\bfphi]-\bP[\ubU,\Id]+\cN_0[\bV,\bfphi]\,,
\end{equation}
with 
\begin{align*}
L_{\bfphi}\bV
&:=\bL_{\tbW}\bP[\ubU,\Id-\bfphi](\bV)\,,\\
\cN_0[\bV,\bfphi]&=\cN_0(\nabla\bfphi,\nabla^2\bfphi,\bV,\nabla\bV)
:=\bP[\ubU+\bV,\Id-\bfphi]
-\bP[\ubU,\Id-\bfphi]-\bL_{\tbW}\bP[\ubU,\Id-\bfphi](\bV)\,,\\
\cN[\bV,\bfphi]&=\cN(\bfphi_t,\nabla\bfphi,\nabla^2\bfphi,\bV,\nabla\bV,\nabla^2\bV)\\
&:=\ \left(\bL_{\tbW}\bP[\ubU,\Id-\bfphi]-\bL_{\tbW}\bP[\ubU,\Id]\right)(\bV)
+\cN_0[\bV,\bfphi]
-(\bfphi_t \cdot \nabla\bfphi\,[\I_2-\nabla\bfphi]^{-1}\nabla)\ubU\\[0.5em]
&
\quad-(\bfphi_t\cdot [\I_2-\nabla\bfphi]^{-1}\nabla) \bV
+\bP[\ubU,\Id-\bfphi]-\bP[\ubU,\Id]-\bL_{\bPhi}\bP[\ubU,\Id](-\bfphi)\,.
\end{align*}
\el

\begin{proof}
To begin with, note that if $\bPsi$ is defined by $\bPsi(t,\cdot):=\bPhi(t,\cdot)^{-1}$, then
\begin{align*}
\nabla\bPsi(t,\bfx)&=[\nabla\bPhi(t,\bPsi(t,\bfx))]^{-1}\,,&
\d_t\bPsi(t,\bfx)&=-[\transp{(\nabla\bPhi(t,\bPsi(t,\bfx)))}]^{-1}
\,\d_t\bPhi(t,\bPsi(t,\bfx))\,.
\end{align*}
Define $\tbW$ by $\bW(t,\bPhi(t,\bfx))=\tbW(t,\bfx)$ or equivalently $\tbW(t,\bfx)=\bW(t,\bPsi(t,\bfx))$. Then $\bW$ solves \eqref{rd} if and only if $(\bPhi,\tbW)$ satisfies
\begin{align*}
\tbW_t- (\bPhi_t \cdot [\nabla\bPhi]^{-1}\nabla) \tbW
&\ =\ \bP[\tbW,\bPhi]\,.
\end{align*}
Inserting $\tbW=\ubU+\bV$ and $\bPhi=\Id-\bfphi$, we readily deduce \eqref{veq0} and derive \eqref{veq} by combining it with \eqref{eq:key-lin}.
\end{proof}

The form \eqref{veq} is adapted to the large-time analysis whereas the form \eqref{veq0} is used in nonlinear regularity estimates. Obviously, we do not need full details of the expression of $\cN$ to carry out our analysis. Yet, for concreteness' sake, let us observe that $\nabla\bfphi(t,\bfx)$ commutes with $\I_2-\nabla\bfphi(t,\bfx)$ and
\be\label{eg:nabla_bfphi}
[\I_2-\nabla\bfphi(t,\bfx)]^{-1}\ =\ \I_2+\nabla\bfphi(t,\bfx)\,[\I_2-\nabla\bfphi(t,\bfx)]^{-1} =\ \I_2+\nabla\bfphi(t,\bfx)+(\nabla\bfphi(t,\bfx))^2\,[\I_2-\nabla\bfphi(t,\bfx)]^{-1}\,,
\ee
so that we have the pointwise estimate
\be\label{e:pointwise_N}
\| \cN[\bV,\bfphi] \| \lesssim \| \nabla \bV \| \left( \| \bfphi_{t} \| + \| \nabla \bfphi \| + \| \nabla^{2} \bfphi \| \right) + \| \nabla^{2} \bV \| \| \nabla \bfphi \| +  \| \nabla \bfphi \| \left( \| \nabla^{2} \bfphi \| + \| \bfphi_{t} \| + \| \nabla \bfphi \| \right).
\ee
The main upshot of Lemma~\ref{lem:cancellation-separation} is that, as long as $\|\nabla\bfphi\|_{L^\infty}<1$, \eqref{rd} is equivalently written as
\[
\bV(t,\cdot)+(\bfphi(t,\cdot) \cdot \nabla) \ubU
\,=\,S(t)[\bV_0+(\bfphi_0 \cdot \nabla) \ubU]
+\int_0^tS(t-\tau)\,\cN[\bV(\tau,\cdot),\bfphi(\tau,\cdot)]\,\dD \tau\,.
\]
At this stage, we need to make a choice so as to split the foregoing equation. We would like to simply use the semigroup splitting of the linear analysis but we need to enforce $\bfphi(0,\cdot)=\bfphi_0$. To do so, we pick $\tchi$ a smooth function on $\R_+$ valued in $[0,1]$, compactly supported in $[0,1]$ and equal to $1$ on $[0,\tfrac12]$. Then, we consider
\begin{align}\label{bfphi}
\bfphi(t,\cdot)&=s(t)[\bV_{0} + (\bfphi_{0} \cdot \nabla) \ubU] 
+ \int_{0}^{t} s(t-\tau)\cN[\bV(\tau,\cdot),\bfphi(\tau,\cdot)]\dD \tau\\
&\ +\tchi(t)\,\left(\bfphi_{0}-s(t)[\bV_{0} + (\bfphi_{0} \cdot \nabla) \ubU]\right)
\nn\\\label{V}
\bV(t,\cdot)&=(S_{1}+S_{2})(t) \left[ \bV_{0}+(\bfphi_{0} \cdot \nabla) \ubU \right] + \int_{0}^{t} (S_{1}+S_{2})(t-\tau)\cN[\bV(\tau,\cdot),\bfphi(\tau,\cdot)]\dD\tau\\
&\ -\tchi(t)\,\Big(\left(\bfphi_{0}-s(t)[\bV_{0} + (\bfphi_{0} \cdot \nabla) \ubU]\right)\cdot\nabla\Big)\ubU
\nn
\end{align}
and observe that, as long as $\|\nabla\bfphi\|_{L^\infty}<1$, \eqref{bfphi}-\eqref{V} imply that $\bW$ defined by 
\[
\bW(t,\cdot):=(\ubU+\bV(t,\cdot))\circ (\Id-\bfphi(t,\cdot))^{-1}
\]
satisfies \eqref{rd} with $\bW(0,\cdot):=(\ubU+\bV_0)\circ (\Id-\bfphi_0)^{-1}$. 

\br\label{rk:time-layer}
Though our strategy is inspired from \cite{JNRZ-RD1,JNRZ-conservation}, we point out that we make here a small departure in the way the short-time layer argument is incorporated. The choice in \cite{JNRZ-RD1,JNRZ-conservation} enforces $\bfphi(t,\cdot)\equiv \bfphi_0$ when $0\leq t\leq\tfrac12$ but results in a slightly more cumbersome analogue of \eqref{bfphi}-\eqref{V}. For comparison, we also observe that in \cite{BJNRZ-KS} and in other previous pieces of work where the linear separation is presented in terms of Green functions rather than semigroups, the time-layer is hidden in the definition of the object playing the role of $(s(t))_{t\geq 0}$.
\er

\br\label{rk:uniqueness}
Let us stress that uniqueness in solving \eqref{bfphi}-\eqref{V} is essentially useless since \eqref{bfphi}-\eqref{V} implies but is not equivalent to \eqref{rd}. However, under our assumptions $\bW(0,\cdot)$ is a Lipschitz bounded function and classical theory for semilinear heat equations provides a local well-posedness result for \eqref{rd} for data in $BUC^0(\R^2;\R^n)$, the space of bounded uniformly continuous functions, with blow-up criterion expressed in terms of $\|\bW(t,\cdot)\|_{L^\infty}$.
\er

To some extent, the introduction of the phase $\bfphi$ has turned the semilinear parabolic system \eqref{rd} into a quasilinear parabolic equation. In particular, since $\cN$ contains terms involving $\nabla^2\bV$, using directly the Duhamel formulas \eqref{bfphi}-\eqref{V} to prove the existence of $(\bV,\bfphi)$ satisfying suitable bounds would be, if not impossible, at least extremely inconvenient\footnote{As it would require large-time maximal regularity estimates}. Instead, we shall use a fixed-point scheme, classical for quasilinear equations and involving, here $L^4$, high-frequency slaving energy estimates to close in regularity. The linear estimates used to derive the latter are provided in the following lemma.

\bl\label{l:energy-estimate}
\begin{enumerate}
\item For any $0\leq\eta_0<1$ and any $\ell\in\N$, there exist $\theta>0$ and $C\geq 0$ such that for any $\bfphi$ such that $\|\nabla\bfphi\|_{L^\infty(\R^2;\cM_2(\R))}\leq \eta_0$ and any $\bV\in W^{\ell,4}(\R^2;\R^n)$ such that $\cD_\ell(\bV)<+\infty$ and $L_{\bfphi}\bV\in W^{\ell,4}(\R^2;\R^n)$
\bas
\sum_{|\alpha|=\ell} \int_{\R^2}\,\|\d^\alpha\bV\|^2\,\d^\alpha\bV\cdot \d^\alpha(L_{\bfphi}\bV)
&\leq -\theta\,\cD_\ell(\bV)^4 \,+\,C\, \|\bV\|_{L^4}^4\,\left(1+\|\nabla^2\bfphi\|_{L^\infty}^{2\,(2\ell+1)}\right)\\
&+C \|\nabla\bfphi\|_{W^{\ell,4}}^4\,\|\bV\|_{W^{1,\infty}}^{4}\,
\eas
where
\[
\cD_\ell(\bV):=\left(\sum_{|\alpha|=\ell}\,\int_{\R^2}\,\|\d^\alpha\bV\|^2\,\|\nabla\d^\alpha\bV\|^2\,\right)^{\frac14}\,.
\]
\item For any $\ell\in\N$, there exists $C$ such that for any $\bV\in W^{\ell+1,4}(\R^2;\R^n)$ and any $j\in\N$, $j\leq \ell$,
\begin{align*}
\|\nabla^j\bV\|_{L^4}
&\leq C\,\|\bV\|_{L^4}^{1-\alpha}
\,\cD_\ell(\bV)^{\alpha}\,,&
\textrm{with }\alpha:=\frac{j}{\ell+\tfrac12}\,.
\end{align*}
\end{enumerate}
\el

\begin{proof}
We first prove the second point. Note that since for any $(j,\ell)$, $j\leq \ell$, 
\[
\|\nabla^j\bV\|_{L^4}
\,\lesssim\,\|\bV\|_{L^4}^{\frac{\ell-j}{\ell}}\,\|\nabla^\ell\bV\|_{L^4}^{\frac{j}{\ell}}
\]
it is sufficient to prove the case $j=\ell$. The subcase $j=\ell=0$ is trivial. Moreover, an integration by parts shows that for any $\alpha$ with $|\alpha|\geq1$,
\begin{align*}
\|\d^\alpha\bV\|_{L^4}^4\,
&\lesssim\,\|\d^\alpha\bV\|_{L^4}\ \|\nabla^{|\alpha|-1}\bV\|_{L^4}\ \left\|\|\d^\alpha\bV\|\,\|\nabla\,\d^\alpha\bV\|\,\right\|_{L^2}\,.
\end{align*} 
Thus for any $\ell\geq1$,
\[
\|\nabla^\ell\bV\|_{L^4}^3
\lesssim \|\nabla^{\ell-1}\bV\|_{L^4}^{\frac13}\,\cD_\ell(\bV)^{\frac{2}{3}}
\lesssim 
\|\bV\|_{L^4}^{\frac{1}{\ell}}\,
\|\nabla^\ell\bV\|_{L^4}^{\frac{\ell-1}{\ell}}
\,\cD_\ell(\bV)^{2}\,,
\]
and the result follows.

As for the first part, to underline the core of the argument we begin by proving the case when $\alpha=\b0$. We first observe that
\begin{align*}
\int_{\R^2}
&\|\bV\|^2\,\bV\cdot\transp{\nabla}\left(\transp{(\ubK\,[\I_2-\nabla\bfphi]^{-1})}(\ubK\,[\I_2-\nabla\bfphi]^{-1})\nabla\bV\right)\\
&=-\int_{\R^2}\|\bV\|^2\,\left\|\ubK\,[\I_2-\nabla\bfphi]^{-1}\nabla\bV\right\|^2
-\frac12\int_{\R^2}\,\left\|\ubK\,[\I_2-\nabla\bfphi]^{-1}\nabla(\|\bV\|^2)\right\|^2
&\lesssim -\cD_0(\bV)^4\,.
\end{align*}
From here one deduces readily
\begin{align*}
\int_{\R^2}\,\|\bV\|^2\,\bV\cdot L_{\bfphi}\bV
&\lesssim -\cD_0(\bV)^4
+\|\bV\|_{L^4}^2\,\|\nabla^2\bfphi\|_{L^\infty}\,\cD_0(\bV)^2
+\|\bV\|_{L^4}^4\,(1+\|\nabla^2\bfphi\|_{L^\infty})\\
&\lesssim -\cD_0(\bV)^4
+\|\bV\|_{L^4}^4\,(1+\|\nabla^2\bfphi\|_{L^\infty}^2)\,.
\end{align*}

We now go back to the general case. Combining commutator bounds recalled below with the arguments of the case $\alpha=\b0$ yields for any $\alpha$
\begin{align*}
\int_{\R^2}\,&\|\d^\alpha\bV\|^2\,\d^\alpha\bV\cdot \d^\alpha(L_{\bfphi}\bV)\\
&\lesssim -\cD_{0}(\d^\alpha\bV)^4
\,+\,\|\d^\alpha\bV\|_{L^4}\,\cD_0(\d^\alpha\bV)^2
\,\left(\,\|\bV\|_{W^{|\alpha|,4}}\,(1+\|\nabla^2\bfphi\|_{L^\infty})
+\|\nabla\bfphi\|_{W^{|\alpha|,4}}\,\|\bV\|_{W^{1,\infty}}\,\right)\,.
\end{align*}
Thus, for any $\ell$,
\begin{align*}
\sum_{|\alpha|=\ell}\int_{\R^2}\,\|\d^\alpha\bV\|^2\,\d^\alpha\bV&\cdot \d^\alpha(L_{\bfphi}\bV)
\lesssim -\cD_{\ell}(\bV)^4
\,+\,\|\bV\|_{L^4}^{\frac{1}{\ell+\frac12}}\,\cD_{\ell}(\bV)^{2+\frac{2\ell}{\ell+\frac12}}\,(1+\|\nabla^2\bfphi\|_{L^\infty})\\[-1em]
&\,+\,\|\bV\|_{L^4}^{\frac{\frac12}{\ell+\frac12}}\,\cD_{\ell}(\bV)^{2+\frac{\ell}{\ell+\frac12}}
\,\left(\|\bV\|_{L^{4}}\,(1+\|\nabla^2\bfphi\|_{L^\infty})
+\|\nabla\bfphi\|_{W^{\ell,4}}\,\|\bV\|_{W^{1,\infty}}\right)\,,
\end{align*}
from which the result follows through Young inequalities. We point out that in the case $\ell=1$, as in the case $\ell=0$, one may obtain the same estimate without the $\|\nabla\bfphi\|_{W^{\ell,4}}\,\|\bV\|_{W^{1,\infty}}$ term.
\end{proof}

In the foregoing proof, to carry out integrations by parts, we have used the standard facts that over $\R^d$, $d\in\N^*$, if $\Div \bfa$ is integrable and $\bfa\in \Span\left(\bigcup_{1\leq p<\frac{d}{d-1}}L^{p}\right)$ then the integral of $\Div \bfa$ is zero and that if $\nabla a$ is integrable and $a\in \Span\left(\bigcup_{1\leq p<\infty}L^{p} \right)$ then the integral of $\nabla a$ is zero.

In the foregoing proof, to bound commutator terms, we have used the following standard nonlinear bounds --- to be used intensively later on --- :
\begin{align*}
\|\nabla^\ell(a\,b)\|_{L^p}
&\lesssim \|a\|_{W^{\ell,\infty}}\,\|b\|_{W^{\ell,p}}\,,\\
\|\nabla^\ell(a\,b)\|_{L^p}
&\lesssim \|\nabla^\ell a\|_{L^p}\,\|b\|_{L^\infty}+\|a\|_{L^\infty}\,\|\nabla^{\ell} b\|_{L^p}\,,\\
\|[\d^\alpha,a]b\|_{L^p}
&\lesssim \|\nabla^{|\alpha|} a\|_{L^p}\,\|b\|_{L^\infty}+\|\nabla a\|_{L^\infty}\,\|\nabla^{|\alpha|-1} b\|_{L^p}\,,\\
\|\nabla^\ell(g\circ \bfa)\|_{L^p}
&\lesssim \|\nabla g\|_{W^{(\ell-1)_+,\infty}(\overline{B}(\b0,\|\bfa\|_{L^\infty}))}\,
\left(1+\|\bfa\|_{L^\infty}^{(\ell-1)_+}\right)\,\|\nabla^\ell \bfa\|_{L^p}\,,&
\textrm{ when }g(\b0)=0\,.
\end{align*}
We shall also use 
\begin{align*}
\|\nabla^\ell(g(\bfa)-g(\bfb))\|_{L^p}
&\lesssim
\|\nabla g\|_{L^\infty(\overline{B}_{\bfa,\bfb})}\,\|\nabla^\ell(\bfa-\bfb)\|_{L^p}\\
&\,+\,
\|\nabla g\|_{W^{\ell,\infty}(\overline{B}_{\bfa,\bfb})}\,
(1+\|\bfa\|_{L^\infty}^{(\ell-2)_+}+\|\bfb\|_{L^\infty}^{(\ell-2)_+})
\,(\|\nabla^\ell \bfa\|_{L^p}+\|\nabla^\ell \bfb\|_{L^p})
\,\|\bfa-\bfb\|_{L^\infty}
\end{align*}
where $\overline{B}_{\bfa,\bfb}:=\overline{B}(\b0,\max(\{\|\bfa\|_{L^\infty},\|\bfb\|_{L^\infty}\}))$.

We conclude this subsection devoted to the analysis of the effects of the introduction of $\bfphi$ by making explicit the affine auxiliary problems used to set up a fixed-point problem at the nonlinear level. 

To begin with, note that the cases $\ell=0$ and $\ell=1$ of the first estimate of Lemma~\ref{l:energy-estimate} (or rather the precised versions of its proof) provide the bounds necessary to check by standard\footnote{That is, by approximating $L_{\phiin}$ by piece-wise constant-in-time operators.} arguments that when $\phiin$ is such that $\nabla\phiin\in\cC^0(\R^+;W^{1,\infty}(\R^2;\cM_2(\R)))$ with $\|\nabla\phiin\|_{L^\infty}<1$, $\d_t-L_{\phiin}$ does generate an evolution system $(S_{\phiin}(t,s))_{t\geq s\geq0}$ on $L^4(\R^2;\R^n)$.

Under suitable assumptions on $\phiin$, $\gin$, we claim that we may likewise solve uniquely in $(\Vout,\phiout)$ (in relevant spaces)
\begin{align}\label{Vout}
\Vout(t,\cdot)&=
S_{\phiin}(t,0)(\bV_0)
+\int_{0}^{t} S_{\phiin}(t,\tau)[(-\left(\d_t-L\right)((\phiout \cdot \nabla) \ubU)\,+\,\gin)(\tau,\cdot)]\dD \tau\,,\\\label{phiout}
\phiout(t,\cdot)&=s(t)[\bV_{0} + (\bfphi_{0} \cdot \nabla) \ubU] 
+ \int_{0}^{t} s(t-\tau)[((L_{\phiin}-L)\Vout\,+\,\gin)(\tau,\cdot)]\dD \tau\\
&\quad+\tchi(t)\,\left(\bfphi_{0}-s(t)[\bV_{0} + (\bfphi_{0} \cdot \nabla) \ubU]\right)\,.
\nn
\end{align}
Note that the problem is designed to also ensure that for $t \geq 0$
\begin{align}\label{Vout_bis}
\Vout(t,\cdot)&=(S_{1}+S_{2})(t) \left[ \bV_{0}+(\bfphi_{0} \cdot \nabla) \ubU \right] \\\nonumber
&\quad+ \int_{0}^{t} (S_{1}+S_{2})(t-\tau)[((L_{\phiin}-L)\Vout\,+\,\gin)(\tau,\cdot)]\dD\tau\\\nonumber
&\quad-\tchi(t)\,\Big(\left(\bfphi_{0}-s(t)[\bV_{0} + (\bfphi_{0} \cdot \nabla) \ubU]\right)\cdot\nabla\Big)\ubU\,,
\end{align}
so that 
\be\label{eq:Vout_eq}
\partial_{t} \Vout - L_{\phiin} \Vout = \gin + \tchi(t)\, (\partial_{t} -L_{\phiin} ) \Big( \Big(\left(\bfphi_{0}-s(t)[\bV_{0} + (\bfphi_{0} \cdot \nabla) \ubU]\right)\cdot\nabla\Big)\ubU \Big).
\ee
Since, because of the short-time cut-off $\tchi$, the problem is not invariant by time-translations, the most convenient way to prove the claim is to observe that on any time interval $[0,T]$, there exists an iterate of the natural fixed-point map that is strictly contracting. We omit to provide more details on the statement and the proof of the claim, partly because some of these details are tedious and classical, partly because the other ones are essentially redundant with those used below to prove nonlinear stability.

\subsection{Nonlinear stability: proof of Theorem~\ref{th:nonlinear_stab}}

To prove Theorem~\ref{th:nonlinear_stab}, we set up a fixed point argument on the map that associates with a given $(\Vin,\phiin)$ the solution $(\Vout,\phiout)$ to \eqref{Vout}-\eqref{phiout} with 
\begin{align*}
\gin&:=-(\phiin_t \cdot \nabla\phiin\,[\I_2-\nabla\phiin]^{-1}\nabla)\ubU-(\phiin_t\cdot [\I_2-\nabla\phiin]^{-1}\nabla) \Vin\\[0.5em]
&\qquad+\bP[\ubU,\Id-\phiin]-\bP[\ubU,\Id]-\bL_{\bPhi} \bP [\ubU,\Id](-\phiin)+\cN_0[\Vin,\phiin]\,,
\end{align*}
where $\bP$ and $\cN_0$ are respectively defined through \eqref{def:bP} and in Lemma~\ref{lem:cancellation-separation}.

For some $C_0>0$ to be taken sufficiently large, we consider $\cX$ the space of functions $(\bV,\bfphi)$ such that $(\bV,\nabla\bfphi,\bfphi_t)\in \cC^0(\R^+;L^4)\cap L^\infty(\R^+;W^{2,4})$ and $\bfphi-\bfphi_0\in\cC^0(\R^+;L^4)$ with $\bfphi(0,\cdot)=\bfphi_0$ that satisfy $\|(\bV,\bfphi)\|_{\cX}\leq C_0\,E_0$, where 
\[
\|(\bV,\bfphi)\|_{\cX}:=\sup_{t\geq0}\,(1+t)^{\frac14}\|(\bV,\nabla\bfphi,\bfphi_t)(t,\cdot)\|_{W^{2,4}}\,.
\]
Note that 
\[
\sup_{t\geq0}\,(1+t)^{-\frac34}\|\bfphi(t,\cdot)-\bfphi_0\|_{L^{4}}
\leq \sup_{t\geq0}\,(1+t)^{\frac14}\|\bfphi_t(t,\cdot)\|_{L^{4}}
\leq \|(\bV,\bfphi)\|_{\cX}\,.
\]
We constraint $E_0$ to be sufficiently small compared to $1/C_0$ so that $(\bV,\bfphi)\in\cX$ implies 
\[
\|(\bV,\nabla\bfphi,\bfphi_t)\|_{L^\infty(\R^+;W^{1,\infty})}\,\leq \frac12\,.
\] 
The latter is sufficient to ensure that the constants introduced below do not depend on $C_0$ and $E_0$.

We first show that if $C_0$ is sufficiently large and, accordingly, $E_0$ is sufficiently small then the map $(\Vin,\phiin)\mapsto(\Vout,\phiout)$ introduced above is well-defined from $\cX$ to itself. 

Let us pick $(\Vin,\phiin)\in\cX$ and consider the associated $(\Vout,\phiout)$. From Lemma~\ref{l:energy-estimate}, product estimates and the embedding $W^{1,4}\hookrightarrow L^\infty$, we deduce that \eqref{eq:Vout_eq} implies for some $C, C'>0$ and $\theta>0$,
\begin{align*}
&\|\Vout(t,\cdot)\|_{W^{2,4}}^4\\
&\leq \|\bV_0\|_{W^{2,4}}^4\,\eD^{-\theta\,t}
\,+\,C\,\int_0^t\,\eD^{-\theta\,(t-\tau)}\,
\left(\,\|\Vout\|_{L^{4}}^4+\|\nabla\phiout\|_{W^{2,4}}^4
+\|\phiout_t\|_{W^{1,4}}^4
+\|\gin\|_{W^{1,4}}^4\,\right)(\tau)\,\dD \tau\\
&\leq \|\bV_0\|_{W^{2,4}}^4\eD^{-\theta\,t}
+C'\,\int_0^t\,\eD^{-\theta\,(t-\tau)}
\left(\left(\|\Vout\|_{L^{4}}^4+\|\nabla\phiout\|_{W^{2,4}}^4
+\|\phiout_t\|_{W^{1,4}}^4\right)(\tau)
+\frac{\|(\Vin,\phiin)\|_{\cX}^8}{(1+\tau)}\right)\dD \tau\,.
\end{align*}
Besides, from Propositions~\ref{pr:lin-local} and~\ref{pr:lin-mod} and Lemma~\ref{l:time-bl}, product estimates and the embedding $H^1\hookrightarrow L^4$, we deduce that \eqref{phiout}-\eqref{Vout_bis} imply for some $C, C'>0$ and $\theta>0$
\begin{align*}
&\|\Vout(t,\cdot)\|_{L^{4}}+\|\nabla\phiout(t,\cdot)\|_{W^{2,4}}\\
&\leq \frac{C\,E_0}{(1+t)^{\frac14}}
\,+\,C\,\int_0^t\,\left(\frac{\eD^{-\theta\,(t-\tau)}}{(\min(\{1,t-\tau\}))^{\frac{1}{2}}}+\frac{1}{(1+t-\tau)^{\frac34}}\right)\,
\left(\,\|\nabla\phiin\|_{W^{1,4}}\|\Vout\|_{W^{2,4}}
+\|\gin\|_{L^2}\,\right)(\tau)\,\dD \tau\\
&\leq \frac{C\,E_0}{(1+t)^{\frac14}}
+C'\,\int_0^t\,\frac{1}{(t-\tau)^{\frac{1}{2}}\,(1+t-\tau)^{\frac14}}\,
\left(\frac{\|(\Vin,\phiin)\|_{\cX}}{(1+\tau)^{\frac14}}\|\Vout(\tau,\cdot)\|_{W^{2,4}}
+\frac{\|(\Vin,\phiin)\|_{\cX}^2}{(1+\tau)^{\frac12}}\right)\dD \tau\,.
\end{align*}
Likewise for some $C, C'>0$ 
\begin{align*}
\|\phiout_t(t,\cdot)\|_{W^{2,4}}&\\
\leq \frac{C\,E_0}{(1+t)^{\frac14}}
&+C\,\left(\frac{\|(\Vin,\phiin)\|_{\cX}}{(1+t)^{\frac14}}\|\Vout(t,\cdot)\|_{W^{2,4}}
+\frac{\|(\Vin,\phiin)\|_{\cX}^2}{(1+t)^{\frac12}}\right)\\
&+C'\,\int_0^t\,\frac{1}{(t-\tau)^{\frac{1}{2}}\,(1+t-\tau)^{\frac14}}\,
\left(\frac{\|(\Vin,\phiin)\|_{\cX}}{(1+\tau)^{\frac14}}\|\Vout(\tau,\cdot)\|_{W^{2,4}}
+\frac{\|(\Vin,\phiin)\|_{\cX}^2}{(1+\tau)^{\frac12}}\right)\dD \tau\,.
\end{align*}
As a result, for 
\[
\zeta(t):=\sup_{0\leq \tau \leq t}(1+\tau)^{\frac14}\|(\Vout,\nabla\phiout,\phiout_t)(t,\cdot)\|_{W^{2,4}}\,,
\]
by combining the foregoing inequalities with Gr\"onwall like arguments and direct integrations, one obtains, for some constant $K>0$, for any $t\geq0$ 
\[
\zeta(t)\,\leq K\,(E_0+C_0\,E_0\,\zeta(t)+(C_0\,E_0)^2)
\]
that implies, when $K\,C_0\,E_0\leq\tfrac12$,
\[
\zeta(t)\,\leq 2\,K\,(E_0+(C_0\,E_0)^2)\,.
\]
Therefore, if $C_0>2\,K$ and $E_0$ is sufficiently small (depending on the choice of $C_0$), one concludes that $\cX$ is indeed left invariant by the map.

To conclude the proof, we first point out that $\cX$ is a complete space for the distance
\[
d_\cX((\bV_1,\bfphi_1),(\bV_2,\bfphi_2))
:=\sup_{t\geq0}\,(1+t)^{\frac14}\|(\bV_1-\bV_2,\nabla(\bfphi_1-\bfphi_2),(\bfphi_1-\bfphi_2)_t)(t,\cdot)\|_{L^{4}}
\]
and we leave to the reader to check that, when $E_0$ is sufficiently small, estimates similar to the ones expounded above prove that the map $(\Vin,\phiin)\mapsto(\Vout,\phiout)$ is strictly contracting for the distance $d_\cX$. This achieves, by the Banach fixed-point theorem, the proof of the first part of Theorem~\ref{th:nonlinear_stab}.

There only remains to prove the further bounds. On one hand, using the embedding $H^{\frac32 }\hookrightarrow L^p$, for any $2<p\leq\infty$, for some $K_p$, for any $t\geq0$, 
\begin{align*}
\|(\bV,\nabla\bfphi,\bfphi_t)(t,\cdot)\|_{L^{p}}&\\
\leq \frac{K_p\,E_0}{(1+t)^{\frac12-\frac1p}}
&+\frac{K_p\,(C_0\,E_0)^2}{(1+t)^{\frac12}}
+K_p\,(C_0\,E_0)^2\,\int_0^t\,\frac{1}{(t-\tau)^{\frac{3}{4}}\,(1+t-\tau)^{\frac14-\frac1p}}\,
\frac{1}{(1+\tau)^{\frac12}}\dD \tau\,.
\end{align*}
This yields the claimed $L^p$-bounds by integration. On the other hand, when $\ell\in\N$, $\ell\geq2$, by using the smoothing effects of $(s(t))_{t\geq0}$, one derives for some $K_\ell$ and any $t\geq0$
\[
\|(\nabla\bfphi,\bfphi_t)(t,\cdot)\|_{W^{\ell,4}}
\leq \frac{K_\ell\,(E_{0,\ell}+(C_0\,E_0)^2)}{(1+t)^{\frac14}}
\]
and, this may be used to show that for some $K_\ell$ and $\theta_\ell>0$ and any $t\geq0$
\begin{align*}
\|\bV(t,\cdot)\|_{W^{\ell,4}}^4
&\leq \|\bV_0\|_{W^{\ell,4}}^4\,\eD^{-\theta_\ell\,t}
\,+\,K_\ell\,\int_0^t\,\eD^{-\theta_\ell\,(t-\tau)}\,
\left(\,\frac{(C_0\,E_{0,\ell})^4}{(1+\tau)}
+\frac{(C_0\,E_0)^4}{(1+\tau)}\|\bV(\tau,\cdot)\|_{W^{\ell,4}}^4\,\right)\,\dD \tau\,,
\end{align*}
which provides the missing bound by a Gr\"onwall-type argument. Note that in the last estimate, we have crucially used the tame character of product estimates.

\subsection{More localized perturbations: proof of Theorem~\ref{th:nonlinear-localized}}

We now sketch the proof of Theorem~\ref{th:nonlinear-localized}. The arguments being quite similar to the ones for Theorem~\ref{th:nonlinear_stab}, we only stress departures from the foregoing detailed proof.

At the linear level, the main variation is that Proposition~\ref{pr:lin-mod} and Lemma~\ref{l:time-bl} should be replaced with the following proposition and lemma whose proofs are nearly identical since $\| \phi_{HF} \|_{L^{p}} \lesssim \| \nabla \phi \|_{L^{p}}$ and
\[
\widehat{\phi_{LF}} (\bfxi) = - \iD \sum_{j} \frac{\chi(\bfxi) \xi_{j}}{\| \bfxi \|^{2}} \widehat{\partial_{j} \phi} (\bfxi).
\]

\bpr\label{pr:lin-localized}
Assume \cond1-\cond2.
\begin{enumerate}
\item There exists $\theta'>0$, such that, for any $(s,s')\in\R_+$ such that $s'\leq s$ and any $q_0$ such that $2\leq q_0<\infty$, there exists $C_{q_0,s,s'}$ such that for any $t>0$, and any $2\leq p\leq q_0$,
\begin{align*}
\|S_{1}(t)[(\bfphi \cdot \nabla) \ubU]\|_{W^{s,p}}
&\leq \frac{C_{q_0,s,s'}}{(\min(\{1,t\}))^{\frac{(s-s')}{2}}}\,\eD^{-\theta'\,t}
\,\|\nabla\bfphi\|_{L^1\cap H^{s'}}\,.
\end{align*}
\item For any $s\in\R_+$, there exists $C_{s}$ such that for any $2\leq p \leq +\infty$, and any $t\geq0$
\begin{align*}
\|S_{2}(t)[(\bfphi \cdot \nabla) \ubU]\|_{W^{s,p}}
&\leq \frac{C_{s}}{(1+t)^{1-\frac{1}{p}}}\,\|\nabla\bfphi\|_{L^1}\,.
\end{align*}
\item For any $\alpha\in\N^2$, any $\ell\in\N$ and any $2 \leq p \leq +\infty$ such that $|\alpha|+\ell+1-\tfrac2p>0$, there exists $C_{p,\alpha,\ell}$ such that for any $t\geq0$
\begin{align*}
\|\,\d_\bfx^\alpha\,\d_t^\ell\,s(t)[(\bfphi \cdot \nabla) \ubU]\|_{L^p}
&\leq \frac{C_{p,\alpha,\ell}}{(1+t)^{\frac{|\alpha|+\ell+1}{2}-\frac{1}{p}}}\,\|\nabla\bfphi\|_{L^1}\,.
\end{align*}
\end{enumerate}
\epr

Note that here the constraint $|\alpha|+\ell+1-\tfrac2p>0$ is reduced to $(\alpha,\ell)\neq(\b0,0)$ or $p>2$.

\bl\label{l:time-bl-localized}
Assume \cond1-\cond2. For any $\alpha\in\N^2$, any $\ell\in\N$ and any $2 \leq p \leq +\infty$ such that $|\alpha|+1-\tfrac2p>0$, there exists $C_{p,\alpha}$ such that for any $t\geq 0$
\begin{align*}
\|\,\d_\bfx^\alpha\,\left(s(t)[(\bfphi \cdot \nabla)\ubU]-\bfphi\right)\|_{L^p}
&\leq C_{p,\alpha}\,\|\nabla\bfphi\|_{L^1\cap W^{(|\alpha|-1)_+,p}}
\,\begin{cases}
\,(1+t)^{\frac12\,\left(|\alpha|-\frac2p\right)_-}&\quad\textrm{if } |\alpha|-\tfrac2p\neq 0\\
\ln(2+t)&\quad\textrm{otherwise }
\end{cases}\,.
\end{align*}
\el

There are only three more elements that require some change.
\begin{enumerate}
\item To gain more localization on $\cN[\bV(\tau,\cdot),\bfphi(\tau,\cdot)]$, we complement with more standard $L^2$ energy estimates the less usual $L^4$ estimates of Lemma~\ref{l:energy-estimate}. This brings some control on the $L^1$-norm of $\cN[\bV(\tau,\cdot),\bfphi(\tau,\cdot)]$.
\item In the Duhamel formula part of the argument, the contribution of nonlinear terms through $S_2$, or $s$, is analyzed by breaking the integral in two parts. The $\int_0^{t/2}$ part is estimated with $L^1\to L^p$ bounds, whereas the $\int_{t/2}^t$ part is bounded using $L^2\to L^p$ estimates.
\item The last estimate follows from Proposition \ref{pr:lin-local}-(3) together with Remark \ref{r:sobolev_gradient}. Note that $\bfphi_{\infty}$ is the only constant such that $\bfphi_{0} - \bfphi_{\infty}$ belongs to $L^{2}(\R^{2})$ and that $s(t)[(\bfphi_{\infty} \cdot \nabla) \ubU] = \bfphi_{\infty}$.
\end{enumerate}

Incidentally, we point out that the arguments sketched above do prove that nonlinear terms are indeed asymptotically irrelevant in large-time.

\br
We stress that actually one may remove the $\log(2+t)$ factor of the $L^\infty$ estimate. One way to prove this is to use $L^q\to L^\infty$ bounds, with $q>2$, to estimate nonlinear contributions in the Duhamel formulation. This requires techniques beyond those of the present section, expounded in the following one. Yet, since our focus is mostly on critical decay, we shall not provide details for those extra arguments. We point out however that the $L^q\to L^\infty$ bounds, $q>2$, mentioned here scale badly in large-time and therefore are not sufficient to remove similar $\log(2+t)$ factors in the critical case.
\er

\section{Modulational behavior}\label{s:mod}

\subsection{Linear estimates}\label{ss:linear_estim_mod}

In the stability part, the starting point of the phase separation was the normalization $\bfq^{\b0}_j=\d_j\ubU$, $j=1,2$. The following lemma provides a similar, higher-order, spectral normalization to set the frame for a wavenumber identification.

\bl\label{l:spec-wn}
Assume \cond2 and consider a wave parametrization as in Proposition~\ref{p:structure}. 
\begin{enumerate}
\item For any $\bfeta \in \R^2$,
\[
\dD_{\bfxi}\bD_{\b0}(\bfeta)\,=\,-\iD \transp{\ubK}\bp \dD_{\bK_1}\bfc(\ubK)(\ubK\bfeta) & \dD_{\bK_2}\bfc(\ubK)(\ubK\bfeta)\ep\,.
\]
\item For any particular choice of the wave profile parametrization, one may normalize $((\bfq^{\bfxi}_1,\bfq^{\bfxi}_2))_{\bfxi}$, $((\tbq^{\bfxi}_1,\tbq^{\bfxi}_2))_{\bfxi}$ to ensure that, for any $\bfeta \in \R^2$ and any $j\in\{1,2\}$,
\begin{align}\label{e:normalization-wn}
\dD_{\bfxi}\bfq^{\b0}_{j}(\bfeta)\,=\,
\iD \,\dD_{\bK_j}\bU(\ubK)(\ubK \bfeta)\,.
\end{align}
\end{enumerate}
\el

In the following we will sometimes use the notation $\bA( \iD \bfxi) :=\dD_{\bfxi}\bD_{\b0}(\bfxi)$ in agreement with Appendix \ref{ss:diff}. Note that
\be\label{e:dD_dOmega}
\dD_{\bfxi}\bD_{\b0}(\bfxi)\,=\iD \left( \sum_{j=1}^2\dD_{\bK_j}\bOm(\ubK)\,(\ubK \bfxi)\,\transp{\beD_j} \right) + \iD \transp{(\ubK \bfxi)} \bfc(\ubK) \I_2
\ee
so that Assumption \case{a} or Assumption \case{b} give the behavior of $\dD_{\bfxi}\bD_{\b0}(\bfxi)$ with respect to $\bfxi$.

\begin{proof}
Along the proof, we use notation from Appendix~\ref{s:profile}. To begin with, by differentiating the definition $\bD_{\bfxi}=(\langle \tbq^\bfxi_j;L_\bfxi \bfq^\bfxi_\ell\rangle_{L^2_{\rm per}})_{j,\ell}$ and using \eqref{e:dk}, one derives for any $\bfeta \in \R^2$
\begin{align*}
\dD_{\bfxi}\bD_{\b0}(\bfeta)
&\,=\,\iD\,\bp\langle \tbq^{\b0}_j;\transp{(L^{(1)}[\d_{\ell} \ubU])}\ubK\bfeta\rangle_{L^2_{\rm per}}\ep_{j,\ell}
\,=\,-\iD\,\bp\langle \tbq^{\b0}_j;
(\transp{\ubK}\dD_{\bK_\ell}\bfc(\ubK)(\ubK\,\bfeta) \cdot \nabla)\ubU\rangle_{L^2_{\rm per}}\ep_{j,\ell}
\end{align*}
which yields the first part of the lemma.

Now, by differentiating $(\I-\Pi_{\bfxi})\,L_{\bfxi}[\bfq_{j}^{\bfxi}]\,=\,0$, $j=1,2$, (encoding that $\Pi_\bfxi$ commutes with $L_\bfxi$), one obtains, for any $\bfeta$ and $j$,
\begin{align}\label{e:wn-aux}
(\I-\Pi_{\b0}) \left( L_{\b0}[\dD_{\bfxi} \bfq^{\b0}_{j}(\bfeta)] 
+\iD \transp{(L^{(1)}[\d_{j} \ubU])}\ubK\bfeta\right)\,=\,0
\end{align}
which, combined with \eqref{e:dk}, implies
\begin{align*}
L_{\b0}[\dD_{\bfxi} \bfq^{\b0}_{j}(\bfeta)
-\iD\,\dD_{\bK_j}\bU(\ubK)(\ubK\bfeta)]
&\,=\,(\I-\Pi_{\b0})\,L_{\b0}[\dD_{\bfxi} \bfq^{\b0}_{j}(\bfeta)
-\iD\,\dD_{\bK_j}\bU(\ubK)(\ubK\bfeta)]\\
&\,=\,\iD\,(\I-\Pi_{\b0})\,((\transp{\ubK} \dD_{\bK_j}\bfc (\ubK)(\ubK\bfeta)\cdot \nabla) \ubU)\,=\,0\,.
\end{align*}
Therefore, for any $\bfeta$ and $j$, $\dD_{\bfxi} \bfq^{\b0}_{j}(\bfeta)
-\iD\,\dD_{\bK_j}\bU(\ubK)(\ubK\bfeta)$ belongs to the kernel of $L_{\b0}$.

To conclude the proof, we only need to prove that by replacing $((\bfq^{\bfxi}_1,\bfq^{\bfxi}_2))_{\bfxi}$, $((\tbq^{\bfxi}_1,\tbq^{\bfxi}_2))_{\bfxi}$ with some $((\bfp^{\bfxi}_1,\bfp^{\bfxi}_2))_{\bfxi}$, $((\tbp^{\bfxi}_1,\tbp^{\bfxi}_2))_{\bfxi}$ satisfying the same spectral conditions one may also achieve the extra normalization condition: for any $\bfeta$, $j$ and $\ell$, 
\[
\langle \tbp_\ell^{\b0};
\dD_{\bfxi} \bfp^{\b0}_{j}(\bfeta)\rangle_{L^2_{\rm per}}
\,=\,
\iD\,\langle \tbp_\ell^{\b0};\dD_{\bK_j}\bU(\ubK)(\ubK\bfeta)\rangle_{L^2_{\rm per}}\,.
\]
This may be achieved, for $\bfxi$ sufficiently small, through
\begin{align*}
\bp \bfp^{\bfxi}_1&\bfp^{\bfxi}_2\ep
&:=\bp \bfq^{\bfxi}_1&\bfq^{\bfxi}_2\ep\,
\left(\I_2+\bp \bfa^{(1,1)}\cdot\bfxi&\bfa^{(1,2)}\cdot\bfxi\\\bfa^{(2,1)}\cdot\bfxi&\bfa^{(2,2)}\cdot\bfxi\ep\right)\,,\\
\bp \tbp^{\bfxi}_1&\tbp^{\bfxi}_2\ep
&:=\bp \tbq^{\bfxi}_1&\tbq^{\bfxi}_2\ep\,
\left(\I_2+\bp \bfa^{(1,1)}\cdot\bfxi&\bfa^{(1,2)}\cdot\bfxi\\\bfa^{(2,1)}\cdot\bfxi&\bfa^{(2,2)}\cdot\bfxi\ep^*\right)^{-1}\,,
\end{align*}
with vectors of $\CC^2$, $\bfa^{(\ell,j)}$, determined by: for any $\bfeta$,
\[
\bfa^{(\ell,j)}\cdot\bfeta
\,=\,-\langle \tbq_\ell^{\b0};
\dD_{\bfxi} \bfq^{\b0}_{j}(\bfeta)\rangle_{L^2_{\rm per}}
\,+\,\iD\,\langle \tbq_\ell^{\b0};\dD_{\bK_j}\bU(\ubK)(\ubK\bfeta)\rangle_{L^2_{\rm per}}\,.
\]
\end{proof}

Under normalization \eqref{e:normalization-wn}, we may refine the decomposition of $S(t)$ into
\[
S(t)[\bfg]\,=\,\Smod(t)[\bfg]\,+\,S_{1}(t)[\bfg] +\tS_{2}(t)[\bfg]
\]
with
\begin{align*}
\Smod(t)[\bfg]&\,:=\,
(s(t)[\bfg] \cdot \nabla) \ubU + \sum_{j \in \{1,2\}} \dD_{\bK_{j}}\bU(\ubK)(\ubK \nabla s_{j}(t)[\bfg])\,,\\
(\tS_{2}(t)[\bfg])(\bfx)
&:= \int_{[-\pi,\pi]^{2}} \chi(\bfxi) \eD^{\iD \bfx \cdot \bfxi}  
\bp \bfr_{1}^{\bfxi}(\bfx)& \bfr_{2}^{\bfxi}(\bfx)\ep 
\eD^{t\,\bD_{\bfxi}} 
\bp \left\langle \tbq_{1}^{\bfxi};\check{\bfg}(\bfxi,\cdot) \right\rangle_{L^{2}_{per}}\\ 
\left\langle \tbq_{2}^{\bfxi};\check{\bfg}(\bfxi,\cdot)\right\rangle_{L^{2}_{per}} \ep \dD\bfxi
\end{align*}
where, for $j=1,2$,
\[
\bfr_{j}^{\bfxi}:=\bfq_{j}^{\bfxi}-\bfq_{j}^{\b0}
-\dD_\bfxi\bfq_{j}^{\b0}(\bfxi)\,.
\]

At this stage, one could just mimick the analysis of Section~\ref{ss:stab-lin} and derive for $\tS_2(t)$ estimates similar to those for $S_2(t)$ but with an extra $(1+t)^{-\frac12}$ decay factor. Yet because of the estimation of 
\[
\int_0^t S_1(t-\tau)\,\cN[\bV(\tau,\cdot),\bfphi(\tau,\cdot)]\,\dD\tau
\]
this would limit the $L^p$-decay of $\bcW(t,\cdot)-\bU^{\bcK(t,\cdot)}(\bPsi(t,\cdot))$ in Theorem~\ref{th:mod-behavior} to the $L^2$-decay of $\cN[\bV(t,\cdot),\bfphi(t,\cdot)]$, that is, $(1+t)^{-\frac12}$. To bypass this limitation, we extend the analysis of Section~\ref{ss:stab-lin} to incorporate $L^q\to L^p$ decay with $q>2$. The price to pay is that the analysis is more involved and does not follow readily from Hausdorff-Young estimates.

To replace Hausdorff-Young estimates, in regimes where they are not available, we shall use Green function representations of finite-rank Bloch multipliers,
\begin{align*}
\int_{[-\pi,\pi]^{2}} \,\eD^{\iD\,\bfx\cdot\bfxi}
\bfm(\bfxi,\bfx)\,\left\langle\tbp^{\bfxi}(\cdot)\,;\check{\bfg}(\bfxi,\cdot) \right\rangle_{L^{2}_{per}} \dD\bfxi
&\,=\, \frac{1}{(2\pi)^2} \int_{\R^2}\,\Gamma(\bfx,\bfy)\,\bfg(\bfy)\,\dD\bfy\,,\\
\textrm{where }\qquad
\Gamma(\bfx,\bfy)
&:=\int_{[-\pi,\pi]^{2}} \,\eD^{\iD\,(\bfx-\bfy)\cdot\bfxi}
\bfm(\bfxi,\bfx)\,\overline{\tbp^{\bfxi}(\bfy)} \dD\bfxi
\end{align*}
that are obtained from the explicit expression of the Bloch transform. For our purposes, it it sufficient to prove $L^\infty\to L^\infty$ on this type of operators, thus to bound $\|\Gamma\|_{L^\infty_\bfx(L^1_\bfy)}$. The following lemma provides convenient ways to obtain this kind of bound (with $\bfz$ playing the role of $\bfx-\bfy$). 

\bl\label{l:Lp}
\begin{enumerate}
\item There exists a constant $C$ such that if $\tilde \bfm$ is a smooth function on $\R^2/(2\pi\,\Z)^2\times \R^2$, then
\[
\Big\|\,\bfz\mapsto
\int_{[-\pi,\pi]^{2}} \,\eD^{\iD\,\bfz\cdot\bfxi}
\tilde \bfm(\bfxi,\bfz)\,\dD\bfxi\Big\|_{L^1(\R^2)}
\leq C\,\| \tilde \bfm\|_{L^\infty_\bfz(L^1_\bfxi)}^{\frac12}\,\|\tilde \bfm\|_{L^\infty_\bfz(W^{4,1}_\bfxi)}^{\frac12}\,.
\]
\item There exists a constant $C$ such that if $\tilde \bfm$ is a smooth function on $\R^2/(2\pi\,\Z)^2$, then
\[
\Big\|\,\bfz\mapsto
\int_{[-\pi,\pi]^{2}} \,\eD^{\iD\,\bfz\cdot\bfxi}
\tilde \bfm(\bfxi)\,\dD\bfxi\Big\|_{L^1(\R^2)}
\leq C\,\| \tilde \bfm\|_{L^2}^{\frac12}\,\| \tilde \bfm\|_{H^{2}}^{\frac12}\,.
\]
\end{enumerate}
\el

The form of the periodicity in $\bfxi$ used here is consistent with Remark~\ref{rk:per-xi}.

\begin{proof}
Both estimates hinge on the fact that for any $1\leq p_0\leq \infty$, if $r>2/p_0'$ (with $p_0'$ Lebesgue-conjugate to $p_0$), for any function $\Gamma$
\begin{align}\label{e:mom-L1}
\|\Gamma\|_{L^1}\,\lesssim \|\Gamma\|_{L^{p_0}}^{1-\frac{2}{r\,p_0'}}\,
\|\|\cdot\|^r\,\Gamma\|_{L^{p_0}}^{\frac{2}{r\,p_0'}}\,.
\end{align}
To prove the latter it is sufficient to optizime in $z_0>0$ the bound obtained by splitting the $L^1$-norm between contributions from $\|\bfz\|\leq z_0$ and those from $\|\bfz\|\geq z_0$.

The second bound is then obtained by applying the foregoing with $p_0=2$, $r=2$ and concluding with Parseval identities. The proof of the first bound starts from the foregoing with $p_0=\infty$, $r=4$, and the observation that $\|\bfz\|^4\eD^{\iD\,\bfz\cdot\bfxi}=\Delta_\bfxi^2(\eD^{\iD\,\bfz\cdot(\cdot)})(\bfxi)$, and is concluded by the integration by parts of $\Delta_\bfxi^2$ and triangle inequalities.
\end{proof}
 
With this tools in hands, we may now turn to linear estimates.

\bpr\label{pr:lin-behavior}
Assume \cond1-\cond2.
\begin{enumerate}
\item There exists $\theta'>0$, such that, for any $(s,s')\in(\R_+)^2$ such that $s'\leq s$ and any $2\leq q<\infty$, there exists $C_{s',s,q}$ such that for any $t> 0$
\begin{align*}
\|S_{1}(t)[\bfg]\|_{W^{s,q}}
&\leq \frac{C_{s',s,q}}{(\min(\{1,t\}))^{\frac{(s-s')}{2}}}\,\eD^{-\theta'\,t}
\|\bfg\|_{W^{s',q}}\,.
\end{align*}
\item For any $s\in\R_+$ and any $\beta\in\N^2$, there exists $C_{s,\beta}$ such that for any $2 \leq p \leq +\infty$, $1\leq q \leq p$, and any $t\geq0$
\begin{align*}
\|\tS_{2}(t)[\d_\bfx^\beta\bfg]\|_{W^{s,p}}
&\leq \frac{C_{s,\beta}}{(1+t)^{\frac{1}{2}+\frac{1}{q}-\frac{1}{p}+\min\left(\left\{\frac12,\frac1q\right\}\right)}}\,\|\bfg\|_{L^{q}}\,.
\end{align*}
\item For any $\alpha\in\N^2$, any $\beta\in\N^2$ and any $\ell\in\N$, there exists $C_{\alpha,\ell,\beta}$ such that for any  $2 \leq p \leq +\infty$, $1\leq q \leq p$, and any $t\geq0$
\begin{align*}
\|\,\d_\bfx^\alpha\,\d_t^\ell\,s(t)[\d_\bfx^\beta\bfg]\|_{L^p}
&\leq \frac{C_{\alpha,\ell,\beta}}{(1+t)^{\frac{|\alpha|+\ell}{2}+\frac{1}{q}-\frac{1}{p}+\min\left(\left\{\frac12,\frac1q\right\}\right)-\frac12}}\,\|\bfg\|_{L^{q}}\,.
\end{align*}
\end{enumerate}
\epr

\begin{proof}
We begin by proving the second estimate. We skip the proof of the cases when $q\leq 2$ as nearly identical to the proof of the corresponding estimates in Proposition~\ref{pr:lin-local}. By interpolation, it is sufficient to consider the cases when $s\in\N$ and $p=q=\infty$. Now to bound $\|\d_\bfx^\alpha\tS_{2}(t)[\d_\bfx^\beta\bfg]\|_{L^{\infty}}$, we expand both $\d_\bfx^{\alpha'}\bfr_j^\bfxi(\cdot)$, $j=1,2$, $\alpha'\leq \alpha$, and $\d_\bfx^{\beta'}\tbq_\ell^\bfxi(\cdot)$, $\ell=1,2$, $\beta'\leq \beta$, with respect to $\bfxi$ respectively up to third and first order. The part containing remainders is obtained through integration against a Green function of the form
\[
\int_{[-\pi,\pi]^{2}} \,\eD^{\iD\,(\bfx-\bfy)\cdot\bfxi}
\bfm(t,\bfxi,\bfx,\bfy)\,\dD\bfxi
\]
with $\bfm$ smooth, compactly supported in $\bfxi$ near $\b0$ and such that
\begin{align*}
\|\nabla_\bfxi^\ell \bfm(t,\bfxi,\bfx,\bfy)\|
&\,\lesssim\,
\left(\|\bfxi\|^{(3-\ell)_+}+\|\bfxi\|^3\,(1+t)^\ell\right)\,\,\eD^{-\theta\,t\,\|\bfxi\|^2}\,,&
\ell\geq0\,.
\end{align*}
Applying the first part of Lemma~\ref{l:Lp} provides the required $(1+t)^{-\frac12}$ bound for this part. The part containing only coefficients of the expansions takes the form of a sum of terms obtained by integration against a Green function of the form
\[
\bfa(\bfx)\,\bfb(\bfy)\,\int_{[-\pi,\pi]^{2}} \,\eD^{\iD\,(\bfx-\bfy)\cdot\bfxi}
\bfm(t,\bfxi)\,\dD\bfxi
\]
with $\bfa$, $\bfb$ bounded and $\bfm$ smooth, compactly supported near $\b0$ and such that 
\begin{align*}
\|\nabla_\bfxi^\ell \bfm(t,\bfxi)\|
&\,\lesssim\,
\left(\|\bfxi\|^{(2-\ell)_+}+\|\bfxi\|^2\,(1+t)^\ell\right)\,\,\eD^{-\theta\,t\,\|\bfxi\|^2}\,,&
\ell\geq0\,.
\end{align*}
Applying the second part of Lemma~\ref{l:Lp} provides the required $(1+t)^{-\frac12}$ bound for these terms. This achieves the proof of the second bound. 

The proof of the third bound is omitted since it follows readily from a combination of arguments of the proof of the corresponding bound in Proposition~\ref{pr:lin-local} and of arguments expounded hereabove to prove the second bound. 

We now focus on the first bound. By using short-time parabolic estimates, one may reduce the analysis to the case $s=s'$. Then, we observe that for any $\omega_0>0$, the contribution to $S_1(t)[\bfg]$ of the part of the spectrum with real part larger than $-2\,\omega_0$ takes the form of a finite sum of terms given as
\begin{align*}
\int_{[-\pi,\pi]^{2}} \chi_0(\bfxi) \eD^{\iD \bfx \cdot \bfxi} 
\bp \bfp_1(\bfxi,\bfx)&\cdots&\bfp_{m_0}(\bfxi,\bfx)\ep\,
\eD^{t\,\bE_\bfxi}
\,\bp \langle\tbp_1(\bfxi,\cdot);\check{\bfg}(\bfxi,\cdot)\rangle_{L^{2}_{per}}\\
\vdots\\
\langle\tbp_{m_0}(\bfxi,\cdot);\check{\bfg}(\bfxi,\cdot)\rangle_{L^{2}_{per}}\ep
\dD\bfxi
\end{align*}
with $\chi_0$ a cutt-off function and $\eD^{t\,\bE_\bfxi}$ exponentially decaying in time at a rate independent of $\omega_0$. Applying the first part of Lemma~\ref{l:Lp} to the corresponding Green functions we deduce that for some $\theta'>0$ independent of $\omega_0$ and $q$, its contribution to an $W^{s,q}\to W^{s,q}$ is bounded by a multiple\footnote{Possibly depending on $\omega_0$.} of $\eD^{-\theta'\,t}$. For clarity, let us temparily denote as $S_{1,\omega_0}(t)[\bfg]$ the remaining part. On one hand, arguing as in the proof of Proposition~\ref{pr:lin-local}, one derives that the $H^{s}\to H^{s}$ norm of $S_{1,\omega_0}(t)$ is bounded by a multiple of $\eD^{-\omega_0\,t}$. On the other hand, let us pick some $q<q_0<\infty$ and observe that for some $\omega_1>0$, one deduces from parabolic estimates that the $W^{s,q_0}\to W^{s,q_0}$ norm of $S(t)$ is bounded by a multiple of $\eD^{\omega_1\,t}$, thus, as a consequence of bounds proved so far, so is the $W^{s,q_0}\to W^{s,q_0}$ norm of $S_{1,\omega_0}(t)$. The last bound to prove then follows by chosing $\omega_0$ sufficiently large to ensure 
\[
-\omega_0\,\frac{\frac1q-\frac{1}{q_0}}{\frac12-\frac{1}{q_0}}
+\omega_1\,\frac{\frac12-\frac1q}{\frac12-\frac{1}{q_0}}\,\leq\,-\theta'
\]
and interpolating.
\end{proof}

The only linear bound left to establish before turning to the proof of Theorem~\ref{th:mod-behavior} is a small variation of Proposition~\ref{pr:lin-mod}.
 
\bl\label{l:lin-mod}
Assume \cond1-\cond2. 
\begin{enumerate}
\item There exists $\theta'>0$, such that, for any $(s,s')\in\R_+$ such that $s'\leq s$ and any $(p_0,p_1)$ such that $2<p_0<p_1<\infty$, there exists $C_{p_0,p_1,s,s'}$ such that for any $t>0$, and any $p_0\leq p\leq p_1$,
\begin{align*}
\|S_{1}(t)[(\bfphi \cdot \nabla) \ubU]\|_{W^{s,p}}
&\leq \frac{C_{p_0,p_1,s,s'}}{(\min(\{1,t\}))^{\frac{(s-s')}{2}}}\,\eD^{-\theta'\,t}
\,\|\Delta\bfphi\|_{L^1\cap W^{(s'-2)_+,p}}\,.
\end{align*}
\item For any $s\in\R_+$, there exists $C_{s}$ such that for any $2\leq p \leq +\infty$, and any $t\geq0$
\begin{align*}
\|\tS_{2}(t)[(\bfphi \cdot \nabla) \ubU]\|_{W^{s,p}}
&\leq \frac{C_{s}}{(1+t)^{1-\frac{1}{p}}}\,\|\Delta\bfphi\|_{L^1}\,.
\end{align*}
\end{enumerate}
\el

\subsection{Proof of Theorem~\ref{th:mod-behavior}}\label{proof.th:mod-behavior}

With $(\bV,\bfphi)$ as in the proof of Theorem~\ref{th:nonlinear_stab}, let us consider 
\be\label{e:defZ}
\bZ(t,\cdot)\ :=\ \bV(t,\cdot)\ - \sum_{j} \dD_{\bK_{j}} \bU(\ubK)(\ubK\nabla\bfphi_{j}(t,\cdot))\,.
\ee
Since we are now enforcing normalization~\eqref{e:normalization-wn}, $\bZ$ is equivalently written as 
\begin{align*}
\bZ(t,\cdot)= &(S_{1}+\tS_{2})(t) \left[ \bV_{0}+(\bfphi_{0} \cdot \nabla) \ubU \right] + \int_{0}^{t} (S_{1}+\tS_{2})(t-\tau)\cN[\bV(\tau,\cdot),\bfphi(\tau,\cdot)]\dD\tau\\
&+ \widetilde{\chi}(t) \left( \Smod(t)[\bV_{0}+ (\bfphi_{0} \cdot \nabla) \ubU] - (\bfphi_{0} \cdot \nabla) \ubU -  \sum_{j} \dD_{\bK_{j}} \bU(\ubK)(\ubK \nabla (\bfphi_{0})_j)\right).
\end{align*}
Let us point out that the estimates used to prove Theorem \ref{th:nonlinear_stab} also yield
\begin{align*}
\| \nabla^{2} \bfphi(t,\cdot) \|_{L^{p}} &\lesssim \frac{E_0}{(1+t)^{\frac12-\frac1p}}\,,\ 2 < p< \infty\,,&
\| \nabla^{2} \bfphi(t,\cdot) \|_{L^{\infty}} &\lesssim E_0 \frac{\ln(2+t)}{(1+t)^{\frac12}}\,,
\end{align*}
so that together with the estimates of Theorem \ref{th:nonlinear_stab} and standard product estimates on the pointwise estimate \eqref{e:pointwise_N} we obtain
\begin{align}\label{e:estim_N_q<2}
\|\cN[\bV(t,\cdot),\bfphi(t,\cdot)]\|_{L^q}&\lesssim \frac{E_0^{2}}{(1+t)^{1-\frac1q}}\,,\ \frac43 < q<4\,,&
\|\cN[\bV(t,\cdot),\bfphi(t,\cdot)]\|_{L^4}&\lesssim E_0^{2} \frac{\ln(2+t)}{(1+t)^{\frac34}}\,.
\end{align}
To bound the contribution of nonlinear terms to $\bZ$, we use an $L^{\min(\{p,4\})}\to W^{1,\min(\{p,4\})}$ estimate for $S_1$ and a $L^{q_p}\to L^p$ estimate for $\tS_2$, with $\frac{1}{q_p}:=\tfrac12\,\left(\frac12+\frac1p\right)$. This results in the following bounds
\begin{align*}
\|\bZ(t,\cdot)\|_{L^{p}}
&\lesssim \, E_0 \frac{\ln(2+t)}{(1+t)^{\frac{1}{2}-\frac{1}{p}+\frac{1}{2}\left(\frac12+\frac1p\right)}}\,,\ 2<p<\infty\,,&
\|\bZ(t,\cdot)\|_{L^{\infty}}
&\lesssim\,E_0 \frac{(\ln(2+t))^2}{(1+t)^{\frac{1}{2}+\frac14}}\,.&
\end{align*}

Then by changes of variable, Lemma~\ref{l:var-change} and a quadratic approximation in wavenumbers, we note that
\bas
&\|\bcW(t,\cdot)-\bU^{\bcK(t,\cdot)}(\bPsi(t,\cdot))\|_{L^{p}(\R^2;\R^n)} \lesssim \|\bZ(t,\cdot)\|_{L^{p}} + \Big\|\lnor\nabla\bfphi(t,\cdot)\rnor^2\Big\|_{L^p},\\
&\|\bcK(t,\cdot)-\ubK\|_{L^p(\R^2;\cM_2(\R))} +\|\d_t\bPsi-\ubOm\|_{L^p(\R^2;\R^2)} \lesssim  \| \nabla \bfphi(t,\cdot) \|_{L^p} + \| \d_t \bfphi(t,\cdot) \|_{L^p}.
\eas
Indeed, since $\bcK:=\nabla\bPsi$ and $\bPsi$ is defined from $\bfphi$ as $\bPsi(t,\bfx):=(\Id-\bfphi(t,\cdot))^{-1}\left(\transp{\ubK}\,(\bfx-t\,\ubc)\right)\,$, we notice that
\begin{align}\label{K_dev}
\bcK(t,\bfx)&=\,
\left(\ubK+\ubK\nabla\bfphi(t,\cdot)+
\ubK(\nabla\bfphi(t,\cdot))^2\,\left(\I_2-\nabla\bfphi(t,\cdot)\right)^{-1}\right)\left(\bPsi(t,\bfx)\right)\,,\\\label{dt_Psi_dev}
\d_t\bPsi(t,\bfx)&=\,
-\transp{\bcK(t,\bfx)}\ubc\,+\,\transp{(\ubK^{-1}\bcK(t,\bfx))}\,\d_t\bfphi(t,\bPsi(t,\bfx))\,.
\end{align}
The proof of Theorem~\ref{th:mod-behavior} is then achieved since
\[
\Big\|\lnor\nabla\bfphi(t,\cdot)\rnor^2\Big\|_{L^p}\lesssim
\|\nabla\bfphi(t,\cdot)\|_{L^{2p}}^2
\lesssim \begin{cases}\frac{1}{(1+t)^{\frac12-\frac1p+\frac12}}\,,& 1<p<\infty\,,\\[1em]
\frac{(\ln(2+t))^2}{(1+t)^{\frac12+\frac12}}\,,& p=\infty\,.
\end{cases}
\]

We now specialize the discussion to either Case~\case{a} or Case~\case{b} and refine the estimates correspondingly. 

\subsection{Scalar case}\label{ss:scalar0}

We begin our refined analysis with Subcase~\case{b0} since it requires less changes and is significantly simpler. Since $\bOm(\bK) = - \transp{\bK} \bfc(\bK)$, using \eqref{e:dD_dOmega}, we note that for some $\bfell_{0} \in\R^2$, $\dD_{\bfxi}\bD_{\b0}(\bfxi)  = \bA (\iD \bfxi) = \iD\,\bfell_{0} \cdot\bfxi\,\I_2$.

The main task is to improve Proposition~\ref{pr:lin-behavior}.

\bpr\label{pr:lin-scalar0}
Assume \cond1-\cond2 and Subcase~\case{b0}.
\begin{enumerate}
\item For any $s\in\R_+$ and any $\beta\in\N^2$, there exists $C_{s,\beta}$ such that for any $2 \leq p \leq +\infty$, $1\leq q \leq p$, and any $t\geq0$
\begin{align*}
\|\tS_{2}(t)[\d_\bfx^\beta\bfg]\|_{W^{s,p}}
&\leq \frac{C_{s,\beta}}{(1+t)^{1+\frac{1}{q}-\frac{1}{p}}}\,\|\bfg\|_{L^{q}}\,.
\end{align*}
\item For any $\alpha\in\N^2$, any $\beta\in\N^2$ and any $\ell\in\N$, there exists $C_{\alpha,\ell,\beta}$ such that for any  $2 \leq p \leq +\infty$, $1\leq q \leq p$, and any $t\geq0$
\begin{align*}
\|\,\d_\bfx^\alpha\,\d_t^\ell\,s(t)[\d_\bfx^\beta\bfg]\|_{L^p}
&\leq \frac{C_{\alpha,\ell,\beta}}{(1+t)^{\frac{|\alpha|+\ell}{2}+\frac{1}{q}-\frac{1}{p}}}\,\|\bfg\|_{L^{q}}\,.
\end{align*}
\end{enumerate}
\epr

\begin{proof}
We only indicate departures from the proof of Proposition~\ref{pr:lin-behavior}. We carry out the expansions of $\d_\bfx^{\alpha'}\bfr_j^\bfxi(\cdot)$, $j=1,2$, $\alpha'\leq \alpha$, and $\d_\bfx^{\beta'}\tbq_\ell^\bfxi(\cdot)$, $\ell=1,2$, $\beta'\leq \beta$, with respect to $\bfxi$ respectively up to fourth and second order so that the remainder part does provide the required $(1+t)^{-1}$ extra decay. Then, when we write the part containing only coefficients of the expansions, it involves integrals of the form
\be\label{int:scalarcase}
\int_{[-\pi,\pi]^{2}} \,\eD^{\iD\,(\bfx-\bfy+\bfell_{0}\,t)\cdot\bfxi}
\bfm(t,\bfxi)\,\dD\bfxi
\ee
with $\bfell_{0}$ the common speed introduced above and $\bfm$ such that 
\begin{align*}
\|\nabla_\bfxi^\ell \bfm(t,\bfxi)\|
&\,\lesssim\,
\left(\|\bfxi\|^{(2-\ell)_+}+\|\bfxi\|^{2+\ell}\,(1+t)^\ell\right)\,\,\eD^{-\theta\,t\,\|\bfxi\|^2}\,,&
\ell\geq0\,.
\end{align*}
Applying the second part of Lemma~\ref{l:Lp} provides the required $(1+t)^{-1}$ bound for these terms. This achieves the proof. 
\end{proof}

We now prove Theorem~\ref{th:scalar} restricted to Subcase~\case{b0}.

To remove the $\log$-factor in the estimates of Theorem~\ref{th:nonlinear_stab}, we notice that Proposition~\ref{pr:lin-scalar0} implies
\begin{align*}
\|S_2(t)[\d_\bfx^\beta\bfg]\|_{W^{s,p}}
&\leq \frac{C_{s,\beta}}{(1+t)^{\frac12+\frac{1}{q}-\frac{1}{p}}}\,\|\bfg\|_{L^{q}}\,,&
2\leq p\leq\infty\,,&\ 1\leq q\leq p\,.
\end{align*}
One may use this in $L^\infty$-bounds to replace $L^2\to L^\infty$ bounds with an $L^q\to L^\infty$ bound for some $2<q<\infty$ together with \eqref{e:estim_N_q<2} when estimating the $\int_{t/2}^t$ part of the nonlinear contribution of $S_2$. Similar estimates also give
\[
\|\cN[\bV(t,\cdot),\bfphi(t,\cdot)]\|_{L^q} \lesssim \frac{E_0^{2}}{(1+t)^{1-\frac1q}}\,,\ \frac43 < q \leq 4 \,,\,\, \|\bZ(t,\cdot)\|_{L^{p}} \lesssim \, \frac{E_0 \ln(2+t) }{(1+t)^{1-\frac{1}{p}}} \,, 2< p \leq 4.
\]
To optimize improvements in $L^p$ estimates of Theorem~\ref{th:mod-behavior} when $p>4$, we first need to derive sharp estimates for $\|\cN[\bV(t,\cdot),\bfphi(t,\cdot)]\|_{L^p}$. This follows readily from the following $L^p$-version of Lemma~\ref{l:energy-estimate}, whose proof is essentially identical to the one of the original lemma, hence omitted.

\bl\label{l:energy-estimate-Lp}
\begin{enumerate}
\item For any $0\leq\eta_0<1$, any $2\leq p<\infty$ and any $\ell\in\N$, there exist $\theta>0$ and $C\geq 0$ such that for any $\bfphi$ such that $\|\nabla\bfphi\|_{L^\infty(\R^2;\cM_2(\R))}\leq \eta_0$ and any $\bV\in W^{\ell,p}(\R^2;\R^n)$ such that $\cD_{\ell,p}(\bV)<+\infty$ and $L_{\bfphi}\bV\in W^{\ell,p}(\R^2;\R^n)$
\begin{align*}
\sum_{|\alpha|=\ell}&\int_{\R^2}\,\|\d^\alpha\bV\|^{p-2}\,\d^\alpha\bV\cdot \d^\alpha(L_{\bfphi}\bV)\\
&\leq -\theta\,\cD_{\ell,p}(\bV)^p
\,+\,C\,\left(\,\|\bV\|_{L^p}^p\,\left(1+\|\nabla^2\bfphi\|_{L^\infty}^{2\,(2\ell+1)}\right)
+\|\nabla\bfphi\|_{W^{\ell,p}}^p\,\|\bV\|_{W^{1,\infty}}^{p}\,\right)
\end{align*}
where
\[
\cD_{\ell,p}(\bV):=\left(\sum_{|\alpha|=\ell}\,\int_{\R^2}\,\|\d^\alpha\bV\|^{p-2}\,\|\nabla\d^\alpha\bV\|^2\,\right)^{\frac1p}\,.
\]
\item For any $\ell\in\N$ and any $2\leq p<\infty$, there exists $C$ such that for any $\bV\in W^{\ell+1,p}(\R^2;\R^n)$ and any $j\in\N$, $j\leq \ell$,
\begin{align*}
\|\nabla^j\bV\|_{L^p}
&\leq C\,\|\bV\|_{L^p}^{1-\alpha}
\,\cD_{\ell,p}(\bV)^{\alpha}\,,&
\textrm{with }\alpha:=\frac{j}{\ell+\tfrac{2}{p}}\,.
\end{align*}
\end{enumerate}
\el

With this in hands one can show that, if $4\leq p_1<+\infty$ and $E^{(p_{1})}_{0}<\infty$, the $W^{2,p_{1}}$-norm of $(\bV, \nabla \bfphi, \bfphi_{t})$ globally exists in time with
\begin{align*}
\|(\bV,\nabla\bfphi,\bfphi_t)(t,\cdot)\|_{W^{2,p}}&
\leq \frac{K_{p_1}\,E_0^{(p_1)}}{(1+t)^{\frac12-\frac1p}}\,,&
4\leq p \leq p_1\,.
\end{align*}
Combining them with, for any $2< p_0\leq 4$,
\begin{align*}
\|(\bV,\nabla\bfphi,\nabla^{2} \bfphi,\bfphi_t)(t,\cdot)\|_{L^p}&
\leq \frac{K_{p_0}\,E_0}{(1+t)^{\frac12-\frac1p}}\,,&
p_0\leq p \leq \infty\,.
\end{align*}
this yields, for any $\tfrac43<q_0\leq 4\leq p_1<+\infty$,
\begin{align*}
\|\cN[\bV(t,\cdot),\bfphi(t,\cdot)]\|_{L^q}&\leq 
\frac{K_{q_0,p_1}\,E_0\,E_0^{(p_1)}}{(1+t)^{1-\frac1q}}\,,&
q_0\leq q \leq p_1\,.
\end{align*}

Then, to bound the contribution of nonlinear terms to $\bZ$, we use an $L^{\min(\{p,p_1\})}\to W^{1,\min(\{p,p_1\})}$ estimate for $S_1$ and a $L^{\min(\{p,p_1\})}\to L^p$ estimate for $\tS_2$. This results in the claimed bounds for $\|\bZ\|_{L^{p}}$.

\subsection{Scalar-like case: proof of Theorem~\ref{th:scalar}}\label{ss:scalar}

We now study Case~\case{b} in general. The main difference with Subcase~\case{b0} is that the required version of Proposition~\ref{pr:lin-scalar0} is significantly harder to prove.

To begin with, note that the assumption ensures that for some $\bfell_0$, $\bfell\in\R^2$, and some invertible $\bP\in\cM_2(\R)$, 
\[
\bP \bA(\iD \bfxi) \bP^{-1} = \iD\,\bfell_0\cdot\bfeta\,\I_2+\diag(\iD\,\bfell\cdot\bfeta,-\iD\,\bfell\cdot\bfeta)\,.
\]
The already analyzed Subcase~\case{b0} corresponds to $\bfell=\b0$ and, thus, we assume here $\bfell\neq\b0$. As in the proof of Lemma~\ref{l:special}, we point out that with
\begin{align*}
\bP\left(\bD_\bfxi-\bA(\iD \bfxi) \right)\bP^{-1}&=:\bp q_{j,k}[\bfxi] \ep_{1\leq j,k\leq2}\,,&
\|\bfxi\|\leq \xi_0\,,
\end{align*}
one has for some $\theta>0$ and any $\bfxi\in\R^2$ such that $\|\bfxi\|\leq \xi_0$,
\[
\Re(q_{j,j}[\bfxi])\leq -\theta\,\|\bfxi\|^2\,,\qquad j=1,2\,.
\]
Introducing 
\begin{align*}
\bSigma_\bfxi(t)&:=
\eD^{-\iD\,\bfell_0\cdot\bfxi\,t}\,\diag(\eD^{-\iD\,\bfell\cdot\bfxi\,t},\eD^{\iD\,\bfell\cdot\bfxi\,t})
\,\bP\eD^{t\,\bD_\bfxi}\bP^{-1}\,,&\|\bfxi\|\leq \xi_0\,,
\end{align*}
one derives for $(k,j)\in \{(1,2);(2,1)\}$, 
\begin{align*}
(\bSigma_\bfxi)_{k,k}(t)&=
\eD^{q_{k,k}[\bfxi]\,t}
+\int_0^t\,\eD^{q_{k,k}[\bfxi]\,(t-\sigma)}\,
\eD^{(-1)^k\,2\,\iD\,\bfell\cdot\bfxi\,\sigma}
\,q_{k,j}[\bfxi]\,(\bSigma_\bfxi)_{j,k}(\sigma)\,\dD\sigma\,,\\
(\bSigma_\bfxi)_{j,k}(t)&=
\int_0^t\,\eD^{q_{j,j}[\bfxi]\,(t-\sigma)}\,
\eD^{(-1)^j\,2\,\iD\,\bfell\cdot\bfxi\,\sigma}
\,q_{j,k}[\bfxi]\,(\bSigma_\bfxi)_{k,k}(\sigma)\,\dD\sigma\,.
\end{align*}
Note that the fact that $\bfell\neq0$ introduces a significant anisotropy in the way the solution spreads. To measure this, with $\bfell^\perp:=(-\ell_2,\ell_1)$ (where $\bfell:=(\ell_1,\ell_2)$), let us introduce adapted coordinates $\xi_\mypar:=\bfell\cdot\bfxi$, $\xi_\perp:=\bfell^\perp\cdot\bfxi$, and denote as $\d_{\mypar}$, $\d_\perp$ corresponding partial derivatives. Then, from an integration by parts in the above time integrals one derives that for any $a\in\N$, $b\in\N$, there exists $C$ and $\theta>0$ such that for any $\bfxi\in\R^2$ such that $\|\bfxi\|\leq \xi_0$, for any $t\geq0$, and $j\neq k$
\begin{align}\label{e:Sig-aux}
\lnor\d_\perp^a\,\d_\mypar^b\,(\bSigma_{j,k})_\bfxi(t)\rnor
&+\lnor\d_\perp^a\,\d_\mypar^b\,(\bSigma_{k,k}-\eD^{q_{k,k}[\cdot]\,t})_\bfxi(t)\rnor\\
&\,\leq\,
C\,\eD^{-\theta\,t\,\|\bfxi\|^2}\,(1+t)^{\frac{a+2b}{2}}\left(|\xi_\mypar|+\min\left(\left\{1,\frac{|\xi_\perp|^2}{|\xi_\mypar|}\right\}\right)\right)\,.\nn
\end{align}

\bpr\label{pr:lin-scalar}
Assume \cond1-\cond2 and Case~\case{b}. Then estimates of Proposition~\ref{pr:lin-scalar0} still hold.
\epr

\begin{proof}
We only show how to bound the contributions to the $L^\infty\to L^\infty$ bounds of the new type of terms arising from the integral terms discussed above. We have to estimate in $L_{\bfz}^{\infty}(\R^{2})$ integrals of the form
\[
\int_{[-\pi,\pi]^{2}} \,\eD^{\iD\,(\bfz + \bfell_0  \pm \bfell )\cdot\bfxi} \eD^{q_{k,k}[\bfxi]\,t} \chi(\bfxi) \bfxi^{\bfalpha} \,\dD\bfxi,
\]
and
\be\label{int:scalarlike}
\int_{[-\pi,\pi]^{2}} \,\eD^{\iD\,(\bfz + \bfell_0  \pm \bfell )\cdot\bfxi} (\bSigma_{j,k})_{\bfxi}(t) \chi(\bfxi) \bfxi^{\bfalpha} \,\dD\bfxi \text{ , } \int_{[-\pi,\pi]^{2}} \,\eD^{\iD\,(\bfz + \bfell_0  \pm \bfell )\cdot\bfxi} ((\bSigma_{k,k})_{\bfxi}(t) - \eD^{q_{k,k}[\bfxi]\,t}) \chi(\bfxi) \bfxi^{\bfalpha} \,\dD\bfxi,\\
\ee
where $j \neq k$ and $\bfxi^{\bfalpha} = \xi_{1}^{\alpha_{1}} \xi^{\alpha_{2}}$ is a monomial with $\bfalpha \in \N^{2}$. The first term can easily be estimated from Lemma \ref{l:Lp}-(2) so we focus our attention on the last two terms. We proceed in two steps.

To begin with, we replace estimate~\eqref{e:mom-L1} in the proof of Lemma \ref{l:Lp} with a suitable anisotropic estimate. Namely, if $a_1$, $b_1$, $a_2$, $b_2$ are real numbers such that 
\begin{align}\label{e:aj-bj}
b_1&<\frac12<a_1\,,&
a_2&<\frac12<b_2\,,&
\left(a_1-\frac12\right)\,\left(b_2-\frac12\right)&>\left(\frac12-b_1\right)\,\left(\frac12-a_2\right)\,,
\end{align}
there exists a constant $C$ such that for any function $\Gamma$ of $\bfz=(z_\perp,z_\mypar)\in\R^2$, 
\begin{align*}
\|\Gamma\|_{L^1}\,\leq\,&C\,\|\Gamma\|_{L^2}^{\frac{2a_1 b_2- 2b_1 a_2+a_2-b_2+b_1-a_1}{2a_1 b_2- b_1 a_2}}\ds\\
&\times\||z_\perp|^{a_1}\,|z_\mypar|^{b_1}\,\Gamma\|_{L^2}^{\frac{b_2-a_2}{2a_1 b_2- 2b_1 a_2}}\ds\\
&\times\||z_\perp|^{a_2}\,|z_\mypar|^{b_2}\,\Gamma\|_{L^2}^{\frac{a_1-b_1}{2a_1 b_2- 2b_1 a_2}}\ds\,,
\end{align*}
though it is sufficient for this work to consider the case $a_{1}=b_{2}=1$ and $a_{2}=b_{1}=\frac14$. To prove the foregoing claim, we set 
\[
\beta:=\sqrt{\frac{(1-2\,a_2)(2\,a_1-1)}{(1-2\,b_1)(2\,b_2-1)}}
\]
so that 
\begin{align*}
2\,a_1-(1-2\,b_1)\,\beta&>1\,,&\textrm{and}&&2\,b_2-\frac{(1-2\,a_2)}{\beta}&>1\,,
\end{align*}
and, for any $r>0$, $R>0$, we split $\R^2$ into three zones defined respectively by 
\begin{align*}
\left(\frac{|z_\perp|}{R}\leq 1\quad\textrm{and}\quad\frac{|z_\mypar|}{r}\leq 1\right),&&
\left(\frac{|z_\perp|}{R}\geq 1\quad\textrm{and}\quad\frac{|z_\mypar|}{r}\leq \frac{|z_\perp|^\beta}{R^\beta}\right),&&
\left(\frac{|z_\mypar|}{r}\geq 1\quad\textrm{and}\quad\frac{|z_\mypar|}{r}\geq \frac{|z_\perp|^\beta}{R^\beta}\right),&&
\end{align*}
to derive
\begin{align*}
\|\Gamma\|_{L^1}&
\lesssim \|\Gamma\|_{L^2}\,\sqrt{r\,R}
+\||z_\perp|^{a_1}\,|z_\mypar|^{b_1}\,\Gamma\|_{L^2}\,\frac{r^{\frac12(1-2\,b_1)}}{R^{\frac12(2\,a_1-1)}}
+\||z_\perp|^{a_2}\,|z_\mypar|^{b_2}\,\Gamma\|_{L^2}\,\frac{R^{\frac12(1-2\,a_2)}}{r^{\frac12(2\,b_2-1)}}\\
&\qquad=\|\Gamma\|_{L^2}\,\sqrt{r\,R}
+\||z_\perp|^{a_1}\,|z_\mypar|^{b_1}\,\Gamma\|_{L^2}\,\frac{(\sqrt{r\,R})^{(1-2\,b_1)}}{R^{(a_1-b_1)}}
+\||z_\perp|^{a_2}\,|z_\mypar|^{b_2}\,\Gamma\|_{L^2}\,\frac{R^{(b_2-a_2)}}{(\sqrt{r\,R})^{(2\,b_2-1)}}\,.
\end{align*}
Optimizing the latter in $(r,R)$ (or equivalently in $(\sqrt{r\,R},R)$) achieves the proof of the claim.

The second step consists in proving that if
\[
\Gamma^{(\alpha)}(t,\bfz)=\int_{[-\pi,\pi]^{2}} \,\eD^{\iD\,\bfz\cdot\bfxi}
\bfm^{(\alpha)}(t,\bfxi)\,\dD\bfxi
\]
where $\alpha \in \N$, $\bfm^{(\alpha)}$ smooth, compactly supported near $\b0$ and if there exists  $N_{0} \in \N$ such that for any $a,b \in \llbracket 0,N_{0} \rrbracket$, we have
\begin{align*}
&|\d_\perp^a\,\d_\mypar^b\,\bfm^{(\alpha)}(t,\bfxi)|\\
&\,\leq\,
C\,\eD^{-\theta\,t\,\|\bfxi\|^2}\,
\left(\|\bfxi\|^{(\alpha-(a+b))_+}+\|\bfxi\|^\alpha \,(1+t)^{\frac{a+2b}{2}}\right)\,
\left(\,|\xi_\mypar|+
\min\left(\left\{1,\frac{|\xi_\perp|^2}{|\xi_\mypar|}\right\}\right)\right)\,,
\end{align*}
then for any  $a,b \in [0,N_{0}]$,
\begin{align}\label{e:Sig-aux2}
\||z_\perp|^a\,|z_\mypar|^b\,\Gamma^{(\alpha)}(t,\cdot)\|_{L^2}
\lesssim \frac{1}{(1+t)^{\frac34+\frac{\alpha}{2}-\frac{a+2b}{2}}}.
\end{align}
By interpolation, it is sufficient to analyze the case when $a$ and $b$ are integers and then one may use
\begin{align*}
\||z_\perp|^a\,|z_\mypar|^b\,\Gamma^{(\alpha)}(t,\cdot)\|_{L^2}
=\|(z_\perp)^a\,(z_\mypar)^b\,\Gamma^{(\alpha)} (t,\cdot)\|_{L^2}
&\lesssim
\|\d_\perp^a\,\d_\mypar^b\,\bfm^{(\alpha)} (t,\cdot)\|_{L^2}\,,
\end{align*}
from which the result follows. Let us point out that in the foregoing, the contribution of the $\min$ term to the required $L^2$ bound is conveniently estimated by splitting the integral in two zones corresponding to $|\xi_\mypar|\leq |\xi_\perp|^2$ and  $|\xi_\perp|^2\leq|\xi_\mypar|$, bounded the factor $e^{-\theta t |\xi_{\mypar}|^2}$ by $1$ when it is necessary.
Combining the two first steps one obtains for any such $\Gamma^{(\alpha)}$ 
\[
\|\Gamma^{(\alpha)}(t,\cdot)\|_{L^1}
\lesssim \frac{1}{(1+t)^{\frac{\alpha}{2}}} \,.
\]
This provides the missing ingredient to complete the proof along the lines of Proposition~\ref{pr:lin-scalar0}.
\end{proof}

We omit the details of the end of the proof of Theorem~\ref{th:scalar} as identical to those for the Subcase~\case{b0}.

\subsection{Dispersive case: proof of Theorems~\ref{th:dispersive} and~\ref{th:dispersive-localized}}\label{ss:dispersive}

We now turn to Case~\case{a}. Once again the main task is to improve estimates on $S_2(t)$, $\tS_2(t)$ and $s(t)$ by relying on the special structure of $\eD^{t\,\bD_\bfxi}$. 

To describe this structure in the present case, let us introduce notation $\eDr(\omega):=(\cos(\omega),\sin(\omega))$ and use the terminology that a function $\bfxi \mapsto A(\bfxi)$ is \emph{smooth in polar coordinates} if it is defined on $\{\,\bfxi\,;\,0<\|\bfxi\|\leq \xi_0\}$, for some $\xi_0>0$, and the map $(r,\omega)\mapsto A(r\,\eDr(\omega))$ extends smoothly from $(0,\xi_0]\times\R/(2\,\pi\Z)$ to $[0,\xi_0]\times\R/(2\,\pi\Z)$. In the following we shall not distinguish between maps $\bfxi\mapsto A(\bfxi)$ and $(r,\omega)\mapsto A(r\,\eDr(\omega))$, thus we identify $\d_\omega A$ with $\bfxi^\perp\cdot\nabla_\bfxi A$ and $\d_r A$ with $\tfrac{\bfxi}{\|\bfxi\|}\cdot\nabla_\bfxi A$. Note that if $A$ is smooth in polar coordinates, it follows that for any $\alpha$, $|\d_\bfxi^\alpha A(\bfxi)|\lesssim \|\bfxi\|^{-|\alpha|}$.

After these preliminary definitions, we observe that in Case~\case{a}, there exist complex-valued maps $\lambda_1$, $\lambda_2$ and complementary projector-valued maps $\pi_1$, $\pi_2$, all smooth in polar coordinates, such that
\[
\eD^{t\,\bD_\bfxi}\,=\,\sum_{j=1}^2\eD^{t\,\lambda_j(\bfxi)}\,\pi_j(\bfxi)
\]
and, for $j=1,2$, $\lambda_j$ is continuous at $\bfxi=\b0$ with value $0$ and for some $\theta>0$,
\[
\Re(\lambda_j(\bfxi))\,\leq -\theta\,\|\bfxi\|^2\,.
\]
Dispersive effects arise from the fact that the dependence of $\d_r\lambda_j$ on $\omega$ is non trivial in the following sense.

\bl\label{l:dispersive-lambda}
By lessening $\xi_0$ if necessary, one may enforce that $\d_r\lambda_j+\d^2_\omega\d_r\lambda_j$ is nowhere vanishing. 
\el 
Note that, in contrast, if $\lambda_j$ were linear in $\bfxi$, then the quantity under study would be identically zero since $\d^2_\omega\eDr=-\eDr$.

\begin{proof}
We first observe that each $\lambda_j$ satisfies
\[
\lambda_j(\bfxi)\stackrel{\bfxi\to\b0}{=}
\iD\,\bfell_0\cdot \bfxi +\epsilon\,\iD\,\sqrt{\bQ(\bfxi,\bfxi)}+\cO(\|\bfxi\|^2)\,,
\]
where $\epsilon\in\{1,-1\}$, $\bfell_0\in\R^2$ and $\bQ$ is a positive definite quadratic form.

It is sufficient to prove that $\bLambda:\,\omega\mapsto \sqrt{\bQ(\beD(\omega),\beD(\omega))}$ is such that $\d^2_\omega\bLambda+\bLambda$ is nowhere vanishing. Direct computations provide
\begin{align*}
%\d_\omega\bLambda(\omega)&=\frac{\bQ(\beD(\omega),\beD(\omega)^\perp)}{\sqrt{\bQ(\beD(\omega),\beD(\omega))}}\\
%\d_\omega^2\bLambda(\omega)&=\frac{\bQ(\beD(\omega)^\perp,\beD(\omega)^\perp)-\bQ(\beD(\omega),\beD(\omega))}{\sqrt{\bQ(\beD(\omega),\beD(\omega))}}-\frac{\bQ(\beD(\omega),\beD(\omega)^\perp)^2}{(\bQ(\beD(\omega),\beD(\omega)))^{\frac32}}\\
(\d_\omega^2\bLambda+\bLambda)(\omega)&=
\frac{\bQ(\beD(\omega),\beD(\omega))\bQ(\beD(\omega)^\perp,\beD(\omega)^\perp)-\bQ(\beD(\omega),\beD(\omega)^\perp)^2}{(\bQ(\beD(\omega),\beD(\omega)))^{\frac32}}
\end{align*}
so that the result stems from the Cauchy-Schwarz inequality.
\end{proof}

\bpr\label{pr:lin-dispersive}
Assume \cond1-\cond2 and Case~\case{a}.
\begin{enumerate}
\item For any $s\in\R_+$ and any $\beta\in\N^2$, there exists $C_{s,\beta}$ such that for any $2 \leq p \leq +\infty$, $1\leq q \leq 2$, and any $t\geq0$
\begin{align*}
\|\tS_{2}(t)[\d_\bfx^\beta\bfg]\|_{W^{s,p}}
&\leq \frac{C_{s,\beta}}{(1+t)^{1 +\frac{1}{q}-\frac{1}{p}+\,\frac12 \min \left( \left\{ \frac12 -\frac{1}{p} , \frac{1}{q} - \frac12 \right\} \right) }} \,\|\bfg\|_{L^{q}}\,.
\end{align*}
\item For any $\alpha\in\N^2$, any $\beta\in\N^2$ and any $\ell\in\N$, there exists $C_{\alpha,\ell,\beta}$ such that for any  $2 \leq p \leq +\infty$, $1\leq q \leq 2$, and any $t\geq0$
\begin{align*}
\|\,\d_\bfx^\alpha\,\d_t^\ell\,s(t)[\d_\bfx^\beta\bfg]\|_{L^p}
&\leq \frac{C_{\alpha,\ell,\beta}}{(1+t)^{\frac{|\alpha|+\ell}{2}+\frac{1}{q}-\frac{1}{p}+\,\frac12 \min \left( \left\{ \frac12 -\frac{1}{p} , \frac{1}{q} - \frac12 \right\} \right)}}\,\|\bfg\|_{L^{q}}\,.
\end{align*}
\item For any $s\in\R_+$, there exists $C_{s}$ such that for any $2\leq p \leq +\infty$, and any $t\geq0$
\begin{align*}
\|\tS_{2}(t)[(\bfphi \cdot \nabla) \ubU]\|_{W^{s,p}}
&\leq \frac{C_{s}}{(1+t)^{\frac54-\frac32\frac{1}{p}}}\,\|\Delta\bfphi\|_{L^1}\,.
\end{align*}
\item For any $\alpha\in\N^2$, any $\ell\in\N$ and any $2 \leq p \leq +\infty$ such that $|\alpha|+\ell-\tfrac{2}{p}>0$, there exists $C_{p,\alpha,\ell}$ such that for any $t\geq0$
\begin{align*}
\|\,\d_\bfx^\alpha\,\d_t^\ell\,s(t)[(\bfphi \cdot \nabla) \ubU]\|_{L^p}
&\leq \frac{C_{p,\alpha,\ell}}{(1+t)^{\frac{|\alpha|+\ell}{2}+\frac14-\frac32\frac{1}{p}}}\,\|\Delta\bfphi\|_{L^1}\,.
\end{align*}
\end{enumerate}
\epr

\begin{proof}
Let us first observe that it is sufficient to prove new $L^1\to L^\infty$ bounds since then one may interpolate with the already known bounds $L^q\to L^2$ and $L^2\to L^p$. This single argument covers all the cases except for $L^p$ bounds on  $\nabla_{\bfx,t}s(t)[(\bfphi \cdot \nabla) \ubU]$ when $2<p<\infty$ since the needed $L^1\to L^2$ bound does not hold. Here, instead, one observes that an $L^1\to L^{2,\infty}$ bound does hold and that it is sufficient so as to apply an interpolation argument. The former claim about the $L^1\to L^{2,\infty}$ bound is essentially shown by studying the most singular part of the associated Green kernel and stems from the fact that $\nabla \Delta^{-1}$ sends $L^{1}$ to $L^{2,\infty}$. We have implicitly used here part of the classical theory on Lorentz spaces and we refer the reader to \cite[Chapter 2]{LM_NS} for the necessary background. 

The first elements of the strategy to prove $L^\infty$ bounds on $\tS_2(t)[\d_\bfx^\beta\bfg]$ or $\d_\bfx^\alpha\,\d_t^\ell\,s(t)[\d_\bfx^\beta\bfg]$ from an $L^1$ bound on $\bfg$ are similar to those of previous subsections: kernel representation, expansions with respect to $\bfxi$ of left and right bases up to a stage where remainders are trivially bounded, identification of terms of the expansions as products of a smooth periodic function of $\bfx$ times a smooth periodic function of $\bfy$ times an integral of the type
\[
\int_{[-\pi,\pi]^{2}} \,\eD^{\iD\,(\bfx-\bfy)\cdot\bfxi+\lambda_j(\bfxi)\,t}
\bfm(\bfxi)\,\dD\bfxi
\]
with $\bfm$ smooth in polar coordinates, compactly supported near $\b0$ and vanishing at least at second order at $\bfxi=\b0$. The other bounds to prove may be obtained along the same lines, the only significant difference lying in the degree of vanishing of $\bfm$ at $\b0$ that may even have a first-order singularity in the worst estimate under consideration.

By using polar coordinates, one deduces that it is thus sufficient to prove bounds uniform with respect to $(r_0,\omega_0)\in (0,+\infty)\times\R/(2\,\pi\Z)$ on integrals of the form
\[
\int_0^{\xi_0}\int_{\R/(2\,\pi\Z)} 
\,\eD^{\iD\,\left(r_0\,\cos(\omega-\omega_0)+\Lambda(\omega)\right)\,r\,t}
\,\eD^{-r^2\,t\,\tLambda(r,\omega)}
\tilde \bfm(r,\omega)\,r^{a+1}\,\dD \omega\,\dD r
\]
with $a\geq-1$, $\tilde \bfm$ compactly supported and smooth in $(r,\omega)$, $\Lambda$ real-valued, smooth and such that $|\Lambda+\Lambda''|$ is positively lower bounded and $\tLambda$ smooth and such that $\Re(\tLambda)$ is positively lower bounded. The control on $\Lambda+\Lambda''$ is provided by Lemma~\ref{l:dispersive-lambda}. It follows from this control (and an examination of the regime $r_0\to+\infty$) that one may split the above integral into a finite number --- controlled uniformly with respect to $(r_0,\omega_0)$ --- of integrals over $(0,\xi_0)\times(\alpha(r_0,\omega_0),\beta(r_0,\omega_0))$ with $|\beta(r_0,\omega_0)-\alpha(r_0,\omega_0)|\leq 2\pi$ and either $|r_0\,\cos(\cdot-\omega_0)+\Lambda|$ or $|-r_0\,\cos(\cdot-\omega_0)+\Lambda''|$ positively lower bounded on $(\alpha(r_0,\omega_0),\beta(r_0,\omega_0))$.

As for pieces on which a control on $|r_0\,\cos(\cdot-\omega_0)+\Lambda|$ is available one may use 
\[
\eD^{\iD\,\left(r_0\,\cos(\omega-\omega_0)+\Lambda(\omega)\right)\,r\,t}
\,=\, \frac{1}{\iD\,\left(r_0\,\cos(\omega-\omega_0)+\Lambda(\omega)\right)\,t}
\d_r\left(\eD^{\iD\,\left(r_0\,\cos(\omega-\omega_0)+\Lambda(\omega)\right)\,\cdot\,t}\right)(r)
\]
and integrate by parts in the $r$ variable. This provides a bound by $t^{-\tfrac32-\tfrac{a}{2}}$, which is $t^{-\tfrac14}$ better than required. Concerning pieces where $|-r_0\,\cos(\cdot-\omega_0)+\Lambda''|$ is under control, we apply the classical van der Corput Lemma to the integral in $\omega$. This yields a total bound by a multiple of
\[
\int_0^{\xi_0}
\,\eD^{-\theta\,r^2\,t}\,\frac{1}{\sqrt{r\,t}}\,\,r^{a+1}\,\dD r
\,\leq\,
\frac{1}{t^{\frac54+\frac{a}{2}}}\,
\int_0^{+\infty}
\,\eD^{-\theta\,r^2}\,r^{a+\frac12}\,\dD r
\]
for some $\theta>0$ and achieves the proof. As for a statement and a proof of the van der Corput Lemma used in the final argument, we refer the reader to either \cite[Lemma~3.3]{R_linKdV} or \cite[Corollary~1.1]{Linares-Ponce}.

Finally, to prove the third and fourth points, we first proceed as in the proof of the second point of Proposition \ref{pr:lin-mod} by decomposing into low and high frequencies. The high-frequency part can be deduced from the previous point and \eqref{estim:high-phi} whereas the low-frequency part follows from \eqref{eq:low-phi} and the strategy used in the previous point.
\end{proof}

\begin{comment}
\br\label{pr:lin-dispersive-q>2}
Interpolating with $L^{p} \to L^p$ bounds for $p \geq 2$, one can get the following estimates:
\begin{enumerate}
\item For any $s\in\R_+$ and any $\beta\in\N^2$, there exists $C_{s,\beta}$ such that for any $2 \leq q \leq p \leq +\infty$, and any $t\geq0$
\begin{align*}
\|\tS_{2}(t)[\d_\bfx^\beta\bfg]\|_{W^{s,p}}
&\leq \frac{C_{s,\beta}}{(1+t)^{\frac12 + \frac74 \frac{1}{q} - \frac34 \frac{1}{p}}} \,\|\bfg\|_{L^{q}}\,.
\end{align*}
\item for any $\alpha\in\N^2$, any $\beta\in\N^2$ and any $\ell\in\N$, there exists $C_{\alpha,\ell,\beta}$ such that for any  $2\leq q \leq p \leq +\infty$, and any $t\geq0$
\begin{align*}
\|\,\d_\bfx^\alpha\,\d_t^\ell\,s(t)[\d_\bfx^\beta\bfg]\|_{L^p}
&\leq \frac{C_{\alpha,\ell,\beta}}{(1+t)^{\frac{|\alpha|+\ell}{2} + \frac74 \frac{1}{q} - \frac34 \frac{1}{p} - \frac12 }} \,\|\bfg\|_{L^{q}}\,.
\end{align*}
\end{enumerate}
\er
\end{comment}

With this in hands, the proof of Theorem~\ref{th:dispersive} is achieved in a straightforward way and we omit the involved details. Likewise, since no new idea nor significantly new estimates are needed to derive Theorem~\ref{th:dispersive-localized}, we skip its proof. 

Let us observe that under the assumptions of Theorem~\ref{th:dispersive} the arguments of the present section also yield estimates for $\|\bcW(t,\cdot)-\bU^{\bcK(t,\cdot)}(\bPsi(t,\cdot))\|_{L^{p}}$ when $p>4$. Yet the decay we would get in a straightforward way is expected to be suboptimal. As described in Remark \ref{rk:log} this spurious limitation is due to a lack of initial regularity.

%Notice that one can easily adapt Proposition \ref{pr:lin-dispersive} with bounds on $\nabla \bfphi$.

\section{Averaged dynamics}\label{s:averaged_dyn}

The last point to be elucidated is the leading-order large-time dynamics of $\bPsi$ and $\nabla\bPsi$ introduced in Theorem~\ref{th:mod-behavior} and its refinements.  In this part we use extensively elements provided in Appendix~\ref{s:WKB}.

\subsection{Linear estimates}

At the linear level, typically we would like to compare $s(t)[(\bfphi_{0} \cdot \nabla) \ubU]$ with $\SigO(t)[\bfphi_{0}]$ where $\SigO$ denotes the evolution operator for
\[
\d_t\bfphi\,=\,
-\transp{\ubK}\dD_{\bK}\bfc(\ubK)(\ubK\nabla_\bfx\bfphi)
+\ubLambda_0[\ubK\nabla_\bfx](\ubK\nabla_\bfx\,\bfphi)\,.
\]
Since the former system has been designed, in Appendix~\ref{ss:artificial}, to match the large-time low-frequency behavior from the system derived in Appendix~\ref{ss:formal}
\[
\d_t\bfphi\,=\,
-\transp{\ubK}\dD_{\bK}\bfc(\ubK)(\ubK\nabla_\bfx\bfphi)
+\bLambda^{\ubK}[\ubK\nabla_\bfx](\ubK\nabla_\bfx\,\bfphi)\,,
\]
at the spectral level it is sufficient to prove that the latter system matches the low-frequency expansion of $\bD_{\bfxi}$, or, with notation from Appendix~\ref{ss:diff}, that it is equivalently written as $\d_t\bfphi\,=\,\bDW(\nabla_\bfx)\bfphi$. In combination with Lemma~\ref{l:spec-wn}, this is the content of the following lemma. 

\br\label{rk:formal-spectral}
Let us stress that if one is not interested in getting an independent formal derivation of the latter system, one could have defined $\bLambda^{\ubK}$ by the condition $\d_t\bfphi\,=\,\bDW(\nabla_\bfx)\bfphi$ and correspondingly skipped the discussion in Appendix~\ref{ss:formal}. In disguise, this is the intermediate choice made in \cite{R_linKdV}.
\er

\bl\label{l:spec-2nd}
Assume \cond2 and consider a wave parametrization as in Proposition~\ref{p:structure}, jointly with normalization~\eqref{e:normalization-wn}. For any $\bfeta \in \R^2$, and $\bfphi \in \R^2$
\begin{align*}
\frac12\dD^{2}_{\bfxi}\bD_{\b0}(\bfeta,\bfeta)\,\bfphi\,=\,\bLambda^{\ubK}[\iD\ubK\bfeta](\iD\ubK\bfeta\,\transp{\bfphi})\,,
\end{align*}
with $\bLambda$ as in \eqref{def:blambda-1}-\eqref{def:blambda-2}.
\el

\begin{proof}
The starting point is, through \eqref{e:normalization-wn} and \eqref{def:Bloch-expand},
\begin{align*}
\frac12\,(\dD^2_{\bfxi}\bD_{\b0}(\bfeta,\bfeta))_{\ell,m}
&\,=\,-\delta_{\ell,m}\,\|\ubK\bfeta\|^2
\,+\,\langle\tbq_\ell^{\b0};
\transp{(L^{(1)}[\d_{\bK_m}\bU(\ubK)(\iD\ubK\bfeta)])}\iD\ubK\bfeta\rangle_{L^2_{per}}\\
&\quad
+\langle \dD_{\bfxi}(\tbq_\ell)^{\b0}(\bfeta);
\transp{(L^{(1)}[\d_m\ubU])}(\iD\ubK\bfeta)
+L_{\b0}[\d_{\bK_m}\bU(\ubK)(\iD\ubK\bfeta)]\rangle_{L^2_{per}}\,. 
\end{align*}
Now, by using \eqref{e:dk} and duality relations, one derives
\begin{align*}
\langle \dD_{\bfxi}(\tbq_\ell)^{\b0}(\bfeta);
&\transp{(L^{(1)}[\d_m\ubU])}(\iD\ubK\bfeta)
+L_{\b0}[\d_{\bK_m}\bU(\ubK)(\iD\ubK\bfeta)]\rangle_{L^2_{per}}\\
&=\,-\sum_{r=1}^2\langle \dD_{\bfxi}(\tbq_\ell)^{\b0}(\bfeta);\d_r\ubU 
\rangle_{L^2_{per}}\,\transp{\ubK}\dD_{\bK_m}\bfc_r(\ubK)(\iD\ubK\bfeta)\\
&=\,\sum_{r=1}^2\langle \tbq_\ell^{\b0};\dD_{\bK_r}\bU(\ubK)(\iD\ubK\bfeta)
\rangle_{L^2_{per}}\,\transp{\ubK}\dD_{\bK_m}\bfc_r(\ubK)(\iD\ubK\bfeta)\,.
\end{align*}
Hence the result.
\end{proof}

With this in hands, the general machinery developed in Appendix~\ref{ss:artificial} provides the relevant comparisons for $\SigO-\SigLF$ where $\SigLF$ is defined as 
\[
\widehat{(\SigLF(t)\bfphi)}(\bfxi)\,=\,\chi(\bfxi)\,\beD^{t\,\bD_\bfxi}\,\widehat{\bfphi}(\bfxi)\,.
\]
Therefore, at the linear level, the remaining task now is to be able to reduce, at leading order, each $s(t)[\bW]$ to a $\SigLF(t)[\bfphi_\bW]$ for a suitable $\bfphi_\bW$. The latter reduction simply arises from the first-order expansion of $\tbq_\ell^{\bfxi}$, $\ell=1,2$, thus bounds on the approximation error share many similarities with bounds on $S_2$, that is arising from the approximation error of the first-order expansion of $\bfq_j^{\bfxi}$. The precise statements are as follows and one proves them similarly to Propositions \ref{pr:lin-scalar} and \ref{pr:lin-dispersive}.

\bpr\label{pr:Wlin-dispersive}
Assume \cond1-\cond2 and Case~\case{a}. 
\begin{enumerate}
\item Let $\bF$ be a smooth $(\beD_1,\beD_2)$-periodic function. For any $\alpha\in\N^2$, any $k\in\N$ such that $k\leq|\alpha|+1$, and any $\ell\in\N$, there exists $C_{\alpha,\ell}$ such that for any  $2 \leq p \leq +\infty$, $1\leq q \leq 2$, and any $t\geq0$
\begin{align*}
\left\|\,\d_\bfx^\alpha\,\d_t^\ell\,s(t)[\varphi\,\bF]
-\d_\bfx^\alpha\,\d_t^\ell\,\SigLF(t)\bp\varphi\,\langle\tbq_1^{\b0};\bF\rangle\\[0.5em]
\varphi\,\langle\tbq_2^{\b0};\bF\rangle\ep\right\|_{L^p}
&\leq 
\frac{C_{\alpha,\ell}}{(1+t)^{\frac{|\alpha|-k+\ell}{2}+\frac{1}{q}-\frac{1}{p}+\frac12+\,\frac12 \min \left( \left\{ \frac12 -\frac{1}{p} , \frac{1}{q} - \frac12 \right\} \right)}}
\,\|\nabla^k\varphi\|_{L^{q}}\,.
\end{align*}
\item For any $\alpha\in\N^2$, any $\ell\in\N$ and any $2 \leq p \leq +\infty$ such that $|\alpha|+\ell+1-\tfrac{2}{p}>0$, there exists $C_{p,\alpha,\ell}$ such that for any $t\geq0$
\begin{align*}
\left\|\,\d_\bfx^\alpha\,\d_t^\ell\,s(t)[(\bfphi\cdot\nabla)\ubU]
-\d_\bfx^\alpha\,\d_t^\ell\,\SigLF(t)[\bfphi]\right\|_{L^p}
&\leq \frac{C_{p,\alpha,\ell}}{(1+t)^{\frac{|\alpha|+\ell}{2}+\frac34-\frac32\frac{1}{p}}}\,\|\Delta\bfphi\|_{L^1}\,.
\end{align*}
\end{enumerate}
\epr

\bpr\label{pr:Wlin-scalar}
Assume \cond1-\cond2 and Case~\case{b}. 
\begin{enumerate}
\item Let $\bF$ be a smooth $(\beD_1,\beD_2)$-periodic function. For any $\alpha\in\N^2$, any $k\in\N$ such that $k\leq|\alpha|+1$, and any $\ell\in\N$, there exists $C_{\alpha,\ell}$ such that for any  $2 \leq p \leq +\infty$, $1\leq q \leq p$, and any $t\geq0$
\begin{align*}
\left\|\,\d_\bfx^\alpha\,\d_t^\ell\,s(t)[\varphi\,\bF]
-\d_\bfx^\alpha\,\d_t^\ell\,\SigLF(t)\bp\varphi\,\langle\tbq_1^{\b0};\bF\rangle\\[0.5em]
\varphi\,\langle\tbq_2^{\b0};\bF\rangle\ep\right\|_{L^p}
&\leq \frac{C_{\alpha,\ell}}{(1+t)^{\frac{|\alpha|-k+\ell}{2}+\frac{1}{q}-\frac{1}{p}+\frac12}}\,\|\nabla^k\varphi\|_{L^{q}}\,.
\end{align*}
\item For any $\alpha\in\N^2$, any $\ell\in\N$ and any $2 \leq p \leq +\infty$ such that $|\alpha|+\ell+1-\frac2p>0$, there exists $C_{p,\alpha,\ell}$ such that for any $t\geq0$
\begin{align*}
\left\|\,\d_\bfx^\alpha\,\d_t^\ell\,s(t)[(\bfphi\cdot\nabla)\ubU]
-\d_\bfx^\alpha\,\d_t^\ell\,\SigLF(t)[\bfphi]\right\|_{L^p}
&\leq \frac{C_{p,\alpha,\ell}}{(1+t)^{\frac{|\alpha|+\ell}{2}-\frac{1}{p}+\frac{1}{2}}}\,\|\Delta\bfphi\|_{L^1}\,.
\end{align*}
\end{enumerate}
\epr

We observe that the scalar products involved in the foregoing propositions shall be ultimately computed by relying on the facts that $\langle\tbq_\ell^{\b0};\d_j\ubU \rangle\,=\,\delta_{j,\ell}$, $1\leq j,\,\ell\leq 2$, and that, when $\bF$ belongs to the range of $L_{\b0}$, $\langle\tbq_\ell^{\b0};\bF\rangle=0$, $\ell=1,2$.

\br\label{lin_per}
We shall use the possibility to trade some time decay against spatial derivatives on $\varphi$ to distribute integrability constraints when estimating contributions of nonlinear terms through Duhamel formula. We point out that to a lesser extent this is also possible when estimating $S_2(t)$ and $\tS_{2}(t)$. This leads to bounds, for $k \leq 1$,
\begin{align*}
\left\| \, S_{2}(t)[\varphi\,\bF] \right\|_{W^{s,p}} \leq
\begin{cases}
\frac{C_{s} \|\nabla^{k} \varphi\|_{L^{q}}}{(1+t)^{\frac{1-k}{2} + \frac{1}{q}-\frac{1}{p} +\,\frac12 \min \left( \left\{ \frac12 -\frac{1}{p} , \frac{1}{q} - \frac12 \right\} \right)}} 
\,,\qquad 1\leq q\leq2\leq p\leq\infty\,, 
&\textrm{ in Case~\case{a}}\,,\\
\frac{C_{s} \|\nabla^{k} \varphi\|_{L^{q}}}{(1+t)^{\frac{1-k}{2} + \frac{1}{q}-\frac{1}{p}}},
\hspace{3.75cm}
 1\leq q\leq p\leq\infty\,,
&\textrm{ in Case~\case{b}}\,,
\end{cases}
\end{align*}
to similar bounds for $(1+t)^{\frac12}\left\| \tS_{2}(t)[\varphi\,\bF] \right\|_{W^{s,p}}$ when $k \leq 2$ and for $(1+t)^{\frac{|\alpha|+\ell-1}{2}} \left\|\partial_\bfx^{\alpha} \partial_{t}^{\ell} \,s(t)[\varphi\,\bF] \right\|_{L^{p}}$ when $|\alpha|+\ell \geq 1$ and $k \leq |\alpha|+\ell+1$.
\er

\subsection{Additional preliminary estimates}\label{ss:additional_estim_whitham}

To prepare the final comparison with averaged equations, we need to transfer a few more properties of the geometrical expansions, measured in powers of $\eps$ in Appendix~\ref{ss:formal}, into large-time asymptotics, measured in powers of $t^{-1}$. The key points to reproduce are that
\begin{enumerate}
\item the leading-order description is of modulation type;
\item the evolution of local parameters is slow;
\item at leading-order one may express time derivatives of $\bfphi$ as a combination of its space derivatives.
\end{enumerate}

The estimates of Section~\ref{s:mod} already prove a version of the first point. However, to analyze nonlinear terms we need a version with higher-order derivatives. Let us observe that the arguments of Section~\ref{s:mod} do yield that, with notation from Theorem \ref{th:nonlinear_stab},
\begin{align*}
\|(\bV, \nabla \bfphi,\bfphi_{t})(t,\cdot)\|_{W^{1,r}} \leq \, C_{p_{0}} \,E_{0} 
\begin{cases}
\frac{\ln(2+t)}{(1+t)^{\frac34-\frac32 \frac1r}}, \qquad\,
2 < p_{0} \leq r \leq \infty
&\textrm{ in Case~\case{a}}\,,\\
\frac{1}{(1+t)^{\frac12-\frac1r}},
\hspace{1.05cm}
2 < p_{0} \leq r \leq \infty
&\textrm{ in Case~\case{b}}\,,
\end{cases}
\end{align*}
\begin{align*}
\|\bZ(t,\cdot)\|_{W^{2,4}} \leq \, C \,E_{0,3}
\begin{cases}
\frac{\ln(2+t)}{(1+t)^{\frac{7}{8}}}&\textrm{ in Case~\case{a}}\,,\\
\frac{\ln(2+t)}{(1+t)^{\frac{3}{4}}}&\textrm{ in Case~\case{b}}\,.
\end{cases}
\end{align*}
These bounds are sufficient to show that $\cN[\bV,\bfphi] = \cN[(\bV-\bZ)+\bZ,\bfphi] $ is a sum of terms of the form $\varphi\,\bF$ with $\bF$ smooth and periodic and $\varphi$ quadratic in $(\bfphi_t,\nabla\bfphi,\nabla^2\bfphi)$, and of a faster decaying remainder (given as a sum of terms that are at least cubic and of quadratic terms involving $\Z$). 

Concerning the second point, we need to prove that placing extra derivatives on $\bfphi$ brings extra decay. All the linear estimates on $s(t)$ contain a version of the latter. We only need to observe that this may be transferred at the nonlinear level into the form
\begin{align*}
\|\nabla^2\bfphi\|_{W^{2,4}} + \|\nabla\d_t\bfphi \|_{W^{2,4}} + \| \bZ(t,\cdot) \|_{W^{2,4}} \leq \, C \,E_{0,3}
\begin{cases}
\frac{1}{(1+t)^{\frac{7}{8}}}&\textrm{ in Case~\case{a}}\,,\\
\frac{1}{(1+t)^{\frac{3}{4}}}&\textrm{ in Case~\case{b}}\,.
\end{cases}
\end{align*}
Note that the initial control follows from the embedding $L^{1} \cap W^{3,4} \hookrightarrow H^{2} \cap W^{2,\infty}$, whereas the propagation is proved by a continuity argument on $\|(\nabla^{2} \bfphi,\nabla \bfphi_{t})(t,\cdot) \|_{W^{2,4}}$, the $\bZ$-bound then following. The main difference with our analysis so far is that we use the above observation about the form of $\cN[\bV,\bfphi]$ and, for the parts of the decomposition of the form $\varphi\,\bF$, with $\bF$ periodic and $\varphi$ quadratic in $(\bfphi_t,\nabla\bfphi)$, we bound their $\int_{t/2}^t$-contributions through Duhamel formula transferring an extra derivative on $\varphi$ thanks to Remark \ref{lin_per}. The same kind of argument also allows to remove almost all the $\ln(2+t)$ in previous bounds so that one can prove the bounds of Remark~\ref{rk:log} except for the bound on $\|\bcW(t,\cdot)-\bU^{\bcK(t,\cdot)}(\bPsi(t,\cdot))\|_{L^{\infty}}$ in Case~\case{a} (that we discuss separately below).

Then, to solve the third point about trading time derivatives for space derivatives, we only need to obtain linear bounds. However $\left(\d_t+\transp{\ubK}\bp \dD_{\bK_1}\bfc(\ubK)(\ubK\nabla) & \dD_{\bK_2}\bfc(\ubK)(\ubK\nabla)\ep\right)s(t)$ satisfies the same bounds as the ones proved for $\nabla^2_{\bfx}s(t)$ thanks to Lemma~\ref{l:spec-wn}, so that $\d_t\bfphi+\sum_j\transp{\ubK}\dD_{\bK_j}\bfc(\ubK)(\ubK\nabla\bfphi_j)$ decays following the bounds proved for $\nabla^2\bfphi$. 

With this in hands, we may now achieve the analysis of nonlinear terms preliminary to comparisons with averaged equations.

\bpr\label{pr:nonlin-W}
Assume \cond1-\cond2 (and normalization \eqref{e:normalization-wn}). Then the pair $(\bV,\bfphi)$ given by Theorem~\ref{th:nonlinear_stab} may be chosen to ensure
\begin{align*}
\cN[\bV,\bfphi]
&:=
\left(\left(
\frac12\dD^{2}_{\bK} \bOm(\bK)(\ubK\nabla\bfphi,\ubK\nabla\bfphi)
+\transp{(\ubK\,\nabla\bfphi)}\dD_{\bK}\bfc(\ubK)(\ubK\nabla\bfphi)
-\transp{\ubK}\dD_{\bK}\bfc(\ubK)(\ubK(\nabla\bfphi)^2)\right) \cdot \nabla\right)\ubU\\
&\quad
-\left(L_{\b0}\left(\frac12\dD_{\bK}^2\bU(\ubK)(\ubK\bfzeta,\ubK\bfzeta)+\dD_{\bK}\bU(\ubK)(\ubK\,\bfzeta^2) \right)\right)_{\bfzeta=\nabla\bfphi}+\bfr\,,
\end{align*}
with
\begin{align*}
\|\bfr\|_{L^{p}} \leq
\begin{cases}\ds
\frac{C_{p_{0}}\,E_{0}\,E_{0,3}}{(1+t)^{\frac34-\frac32\frac{1}{p}+\frac54}}\,,\qquad
&\frac43 <p_0 \leq p \leq 4\,,
\textrm{ in Case~\case{a}}\,,\\[1em]\ds
\frac{C_{p_{0}} \,E_0\,E_{0,3}}{(1+t)^{\frac12-\frac{1}{p}+1}}\ \,,\qquad
&\frac43 <p_0\leq p\leq 4\,,
\textrm{ in Case~\case{b}}\,.
\end{cases}
\end{align*}
\epr

The cumbersome form involving a $\bfzeta$ is chosen to emphasize that the leading-order part of $\cN[\bV,\bfphi]$ has the tensorized form $\varphi\,\bF$ with $\bF$ periodic required to apply the linear estimates of the former subsection. Note also that as announced the involved periodic factors either belong to the range of $L_{\b0}$ or are $\d_j\ubU$ for some $j$. One may also prove estimates of the remainder in $L^{p}$ for any $p \in (1,\infty]$ but we omit those as useless for the rest of our analysis.

\begin{proof} 
Combining \eqref{e:dk} with the above bounds on $\nabla \bfphi$, $\bZ$ and $\nabla^{2} \bfphi$, one derives (with notation from Lemma \ref{lem:cancellation-separation})
\begin{align*}
(\bL_{\tbW}\bP&[\ubU,\Id-\bfphi]-\bL_{\tbW}\bP[\ubU,\Id])(\bV)-(\bfphi_t\cdot [\I_2-\nabla\bfphi]^{-1}\nabla) \bV\\
&=\transp{(\ubK\nabla\bfphi\nabla)}\left(\dD\bG(\ubU)\left(\dD_{\bK}\bU(\ubK)(\ubK\nabla\bfphi)\right)\right)-\left(\dD_\bK\bOm(\ubK)(\ubK\nabla\bfphi)\cdot\nabla\right)\left(\dD_{\bK}\bU(\ubK)(\ubK\nabla\bfphi)\right)\\
&\quad+\transp{\nabla}\left(\,
\left(\transp{\ubK}(\ubK\,\nabla\bfphi)
+\transp{\ubK}(\ubK\,\nabla\bfphi)\right)\nabla\left(\dD_{\bK}\bU(\ubK)(\ubK\nabla\bfphi)\right)\right)
+\bfr_0\\
\bP[\ubU&+\bV,\Id-\bfphi]-\bP[\ubU,\Id-\bfphi]
-\bL_{\tbW}\bP[\ubU,\Id-\bfphi](\bV)\\
&=\dD^{2} \! \bff (\ubU)(\dD_{\bK}\bU(\ubK)(\ubK\nabla\bfphi),\dD_{\bK}\bU(\ubK)(\ubK\nabla\bfphi))\\
&\quad+\transp{(\ubK \nabla)} \dD^{2} \! \bG(\ubU)  (\dD_{\bK}\bU(\ubK)(\ubK\nabla\bfphi],\dD_{\bK} \bU(\ubK)(\ubK\nabla\bfphi)) +\bfr_1\\
-(\bfphi_t \cdot &\nabla\bfphi\,[\I_2-\nabla\bfphi]^{-1}\nabla)\ubU\\
&\quad\,=(\transp{(\ubK\,\nabla\bfphi)}\dD_{\bK}\bfc(\ubK)(\ubK\nabla\bfphi)\cdot \,\nabla)\ubU +\bfr_2\\
\bP[\ubU,&\Id-\bfphi]
-\bP[\ubU,\Id]
-\bL_{\bPhi} \bP [\ubU,\Id](-\bfphi)\\
&=-L_{\b0}[\dD_{\bK}\bU (\ubK)(\ubK(\nabla\bfphi)^2)]
\,-(\transp{\ubK}\dD_{\bK}\bfc(\ubK)(\ubK(\nabla\bfphi)^2) \cdot \nabla)\ubU\\
&\quad\,+\transp{\nabla}\left(\,\transp{(\ubK\,\nabla\bfphi)}(\ubK\,\nabla\bfphi)\nabla\ubU\right)
+\bfr_3\\
\end{align*}
with $\bfr_j$, $j\in\{0,1,2,3\}$, decaying as claimed for $\bfr$ in Proposition~\ref{pr:nonlin-W}. The result then follows by summing the foregoing identities and using \eqref{e:d2k}.
\end{proof}

Thanks to the previous proposition one can also prove that, if $E_{0}$ is small enough,
\begin{align*}
\|\nabla^3\bfphi \|_{W^{1,4}} \leq C \,E_{0,3}
\begin{cases}
\frac{1}{(1+t)^{\frac{11}{8}}}\,,
&\textrm{ in Case~\case{a}} \,,\\
\frac{1}{(1+t)^{\frac{5}{4}}}\,,
&\textrm{ in Case~\case{b}}\,.
\end{cases}
\end{align*}
This last bound allows us to prove the bound on $\|\bcW(t,\cdot)-\bU^{\bcK(t,\cdot)}(\bPsi(t,\cdot))\|_{L^{\infty}}$ in Case~\case{a} in Remark~\ref{rk:log}. Note also that
\begin{align*}
\|\nabla^3\bfphi \|_{L^{r}} \leq
\frac{C \,E_{0,3}}{(1+t)^{\frac{3}{2}-\frac{1}{r}}}\,,\qquad
2 \leq r \leq \infty
&\textrm{ in Case~\case{b}}\,.
\end{align*}

\subsection{Averaged systems: proofs of Theorems \ref{th:whitham}, \ref{th:dispersive-linear-whitham} and \ref{th:subcritical-linear-whitham}}

We have now all elements in hands to prove our last round of main results, on averaged systems. We provide details for Theorem~\ref{th:whitham}, proofs of Theorems~\ref{th:dispersive-linear-whitham} and~\ref{th:subcritical-linear-whitham} following from similar computations.

We first observe that the existence of a global solution $\bPsiW$ to \eqref{e:W-phase-intro} stems from Proposition~\ref{p:Wphi}. Motivated by the estimates of Proposition~\ref{p:Wphi} and following the lines of Subsection~\ref{proof.th:mod-behavior}, we get
\bas
\|\bcK(t,\cdot)  - \bcKW(t,\cdot) \|_{L^{p}} \lesssim  &\|\nabla\bfphi(t,\cdot) -\nabla\Wbphi(t,\cdot) \|_{L^{p}}\\
&+ \|\nabla\Wbphi(t,\cdot) - \nabla\bphiW(t,\cdot) \|_{L^{p}}
+ \|\nabla\bfphi(t,\cdot)\|_{L^{2p}}^2 + \|\nabla\bphiW(t,\cdot)\|_{L^{2p}}^2,
\eas
where bounds on the second and fourth terms of the right-hand side are provided by Proposition~\ref{p:Wphi}, bounds on the third stem from either Theorem~\ref{th:dispersive} or Theorem~\ref{th:scalar}, and $\Wbphi$ satisfies $\Wbphi(0,\cdot) = \bfphi_{0}$ and Equation~\eqref{eq:Wphi}, that we write as
\[
\d_t\Wbphi =-\transp{\ubK}\dD_{\bK}\bfc(\ubK)(\ubK\nabla\Wbphi) +\ubLambda_0[\ubK\nabla](\ubK\nabla\Wbphi) + \cQ(\nabla \Wbphi) 
\]
with
\[
\cQ(\cK) := \frac12\dD_{\bK}^2\bOm(\ubK)(\ubK\cK,\ubK\cK) +\transp{(\ubK\,\cK)}\dD_{\bK}\bfc(\ubK)(\ubK\cK) -\transp{\ubK}\dD_{\bK}\bfc(\ubK)(\ubK\,(\cK)^2).
\]
There only remains to bound $\nabla\bfphi-\nabla\bfphi_{W}$.

On one hand, using notation near~\eqref{def:SigW0} and Lemma~\ref{l:spec-wn}, the Cauchy problem for $\Wbphi$ is equivalently written as
\[
\Wbphi(t) = \SigO(t)[\bfphi_{0}] + \int_{0}^{t} \SigO(t-\tau)[ \cQ(\nabla \Wbphi(\tau))] d\tau.
\]
On the other hand, the preliminary estimates of the present section and estimates of Appendix~\ref{s:WKB} yield
\[
\bfphi(t) = \SigO(t)[\bfphi_{0}] + \int_{0}^{t} \SigO(t-\tau)[\cQ(\nabla \bfphi (\tau))] d\tau + \bR(t)
\]
where the residual $\bR$ satisfies
\begin{align*}
\| \nabla \bR(t) \|_{L^{p}} \leq\begin{cases}\ds
\frac{C \big(E_{0,3}\,(1+ E_{0,3}) +\| \bV_{0} \|_{L^{\frac32}}\big)}{(1+t)^{\frac34-\frac32\,\frac1p + \frac{1}{2}}}\,, \qquad
&2 \leq p \leq \infty 
\textrm{, Case~\case{a}},\\[1em]\ds
\frac{C E_{0,3}}{(1+t)^{\frac12-\frac1p + \frac{1}{2}}} \,,
&2 \leq p \leq \infty 
\textrm{, Case~\case{b0}},\\[1em]\ds
\frac{C_{\eta} E_{0,3}\,(1+ E_{0,3})}{(1+t)^{\frac{1}{2}-\frac1p + \frac12 (1-\frac1p) - \eta}}\,,
&2 \leq p \leq \infty,\, \eta>0
\textrm{, Case~\case{b} but \case{b0} fails}.
\end{cases}
\end{align*}
Actually in the estimate of the remainder $\bR$, when Case~\case{b} holds but Subcase~\case{b0} fails, we have completed Proposition~\ref{pr:ubL0-scalar} with 
\[
\|\nabla (\SigLF-\SigOLF)(t)[\bfg]\|_{L^p}\leq
\frac{C_{r_{0}}}{(1+t)^{\frac{1}{q}-\frac{1}{p}+\frac12\left(\frac1q-\frac1p\right) - \frac12}}\,\| \nabla^2 \bfg\|_{L^{q}}
\]
for any $1 \leq q \leq 2 \leq p \leq \infty$ such that $\tfrac1q - \tfrac1p \geq \tfrac{1}{r_{0}} > \tfrac12$, whose proof is omitted as similar to other estimates of Proposition~\ref{pr:ubL0-scalar}. The proof is then concluded with a continuity argument on $\nabla\bfphi-\nabla\Wbphi$. To provide some details on the latter, we point out for instance that in Case~\case{b} (when \case{b0} fails) when $2\leq p\leq 3$, the $\int_0^{t/2}$ part of the integral is bounded using $L^{\frac32} \to L^p$ estimates whereas the $\int_{t/2}^t$ part is bounded with $L^{\min(p,3)} \to L^p$ estimates.

\appendix

\section{Spectral background}\label{s:spec}

\subsection{The Bloch transform}\label{ss:Bloch}

In the present subsection, we provide main properties of the suitable Bloch transform. We shall be rather bold concerning summation issues and meaning of equalities since, the actual resolution of these questions follows from a combination of the classical arguments, or even results, for the Fourier transform/series. In particular, everything is readily justified when applied to Schwartz-class functions and extensions to more general spaces follow from classical density arguments.

The Bloch transform is designed to ensure the following Bloch-wave decomposition of any function $g$ over $\R^2$
\[
g(\bfx) = \int_{[-\pi,\pi]^{2}} \eD^{\iD \bfxi \cdot \bfx}\,\check{g}(\bfxi,\bfx) \dD\bfxi
\]
where, for each Floquet parameter $\bfxi\in [-\pi,\pi]^2$, $\check{g}(\bfxi,\cdot)$ is $(\beD_1,\beD_2)$-periodic. It is explicitly given as
\[
\cB(g)(\bfxi,\bfx)
=\check{g}(\bfxi,\bfx)\,:=\,
\sum_{\bfp \in \Z^2} \eD^{2\iD\pi \bfp \cdot \bfx} 
\hat{g}(\bfxi + 2\pi \bfp)
\,=\,
\frac{1}{(2\pi)^{2}} \sum_{\bfq \in \Z^2} \eD^{-\iD \bfxi \cdot(\bfx + \bfq)} g(\bfx + \bfq)\,,
\]
where $\hat{\cdot}$ denotes the Fourier transform normalized by
\[
\cF(g)(\bfxi)=\hat{g}(\bfxi)\,:=\,\frac{1}{(2 \pi)^2} \int_{\R^2} \eD^{-\iD \bfxi \cdot \bfx} g(\bfx)\,\dD\bfx
\]
and the equivalence of both formula stems from the Poisson summation formula. 
%Note also that for each real number $\bfx \in [0,1]^2$, $\bfxi \mapsto \check{g}(\bfxi,\bfx) e^{\iD \bfx \cdot \bfxi}$ is $2\pi (\beD_1,\beD_2)$-periodic.

Elementary computational rules are
\begin{align*}
\widecheck{(\nabla_\bfx g)}(\bfxi,\bfx)&\,=\,(\nabla_\bfx+\iD\bfxi)(\check{g})(\bfxi,\bfx)\,,&\\
\widecheck{(g\,h)}(\bfxi,\bfx)&\,=\,h(\bfx)\,\check{g}(\bfxi,\bfx)\,,&
\textrm{when }h\textrm{ is $(\beD_1,\beD_2)$-periodic,}\\
\check{g}(\bfxi,\bfx)&\,=\,\hat{g}(\bfxi)\,,&
\textrm{when }\supp\hat{g}\subset[-\pi,\pi]^2\,,
\end{align*}
where $g$ is scalar, and the Bloch transform is applied coordinate-wise. We say that a function is slow when it satisfies the foregoing support condition on its Fourier transform. Note that, as a consequence of the above relations, when $g$ is slow and $h$ is $(\beD_1,\beD_2)$-periodic, $\widecheck{(g\,h)}(\bfxi,\bfx)\,=\,\hat{g}(\bfxi)\,h(\bfx)$, a relation particularly useful when extracting averaged dynamics from slow modulation behavior.

For any $s\in\N$, $2\pi$ times the Bloch transform provides a total isometry from\footnote{We omit to mark the space in which scalar maps are valued so as to omit the standard discussion between complex-valued maps and real-valued maps, whose Bloch transforms are characterized by an extra symmetry in the $\bfxi$ variable.} $H^s(\R^2)$ to $L^2([-\pi,\pi]^2;H^s_{\rm per}([0,1]^2))$ endowed with respective norms, equivalent to standard norms,
\begin{align*}
g&\mapsto\sqrt{\sum_{|\alpha|\leq s}\|\nabla_\bfx^\alpha g\|_{L^2(\R^2)}^2}\,,&
g&\mapsto\sqrt{\sum_{|\alpha|\leq s}\|\bfxi\mapsto\|(\nabla_\bfx+\iD\bfxi)^\alpha g(\bfxi,\cdot)\|_{L^2([0,1]^2)}\|_{L^2([-\pi,\pi]^2)}^2}\,,&
\end{align*}
with $H^s_{\rm per}([0,1]^2)$ denoting the closure for the $H^s([0,1]^2)$-topology of restrictions to $[0,1]^2$ of smooth $(\beD_1,\beD_2)$-periodic functions. Throughout the text we call these isometry properties Parseval identities. Interpolating between those and simple triangle inequalities one derives inequalities that we call Hausdorff-Young inequalities throughout the text,
\begin{align*}
\|g\|_{L^{p}(\R^2)}&\leq (2\pi)^{\frac{2}{p}} \|\check{g}\|_{L^{p'}([-\pi,\pi]^{2};L^{p}([0,1]^2))}\,,& 2\leq p\leq\infty\,,&\qquad\frac{1}{p}+\frac{1}{p'}=1\,,\\
 \|\check{g}\|_{L^p([-\pi,\pi]^2;L^{p'}([0,1]^2))}&\leq \frac{1}{(2\pi)^{\frac{2}{p'}}} 
 \|g\|_{L^{p'}(\R^2)}\,,& 2\leq p\leq\infty\,,&\qquad\frac{1}{p}+\frac{1}{p'}=1\,.
\end{align*}
By using explicit representations of derivatives when $s\in\N$ and interpolation, those also yield
\begin{align*}
\|g\|_{W^{s,p}(\R^2)}&\lesssim \|\check{g}\|_{L^{p'}([-\pi,\pi]^{2};W^{s,p}([0,1]^2))}\,,& s\in\R\,,\ 2\leq p\leq\infty\,,&\qquad\frac{1}{p}+\frac{1}{p'}=1\,,
\end{align*}
also referred to as Hausdorff-Young inequalities. 

Incidentally, let us point out that throughout the text we also use classical Parseval identities and Hausdorff-Young inequalities, adapted to the Fourier transform.

\subsection{Spectral perturbation}\label{ss:pert}

We gather here some standard facts, specialized to our present analysis, about spectral perturbation analysis as contained in \cite{Kato}. In particular we sketch a proof of \eqref{e:spec-decomp}. We warn the reader that in the present case one can use neither a spectral theorem for self-adjoint operators nor Evans' function arguments based on a spatial dynamics interpretation. Our present account remotely echoes the arguments sketched in \cite[p.30-31]{R} for plane waves.

As a relatively compact perturbation of $\transp{(\ubK\nabla)}(\ubK\nabla)$ acting on $L^2([0,1]^2;\R^n)$ with domain $H^2_{\rm per}(\R^2;\R^n)$, each $L_\bfxi$ has compact resolvents, hence discrete spectrum with finite multiplicity. For each $\bfxi_0$ and $\lambda_0\notin \sigma(L_{\bfxi_0})$,
\[
(\lambda_0\,\I-L_\bfxi)^{-1}\,=\,(\lambda_0\,\I-L_{\bfxi_0})^{-1}\,(\I-(L_{\bfxi_0}-L_{\bfxi})(\lambda_0\,\I-L_\bfxi)^{-1})
\]
provides a smooth representation of the resolvent $(\lambda_0\,\I-L_\bfxi)^{-1}$ when $\bfxi$ is sufficiently close to $\bfxi_0$. This is transferred to spectral projectors through Riesz' formula
\[
\Pi^{\Gamma}_\bfxi\,=\,\frac{1}{2\iD\pi} \int_{\Gamma} \left(\lambda\I- L_{\bfxi} \right)^{-1}\,\dD\lambda
\]
where $\Gamma$ is a simple positively-oriented curve. The spectral projector $\Pi^{\Gamma}_\bfxi$ projects on the sum of generalized eigenspaces associated with eigenvalues of $L_\bfxi$ inside $\Gamma$, its rank providing the sum of algebraic multiplicities of these eigenvalues. Incidentally note that it follows from the formula that $\Pi^{\Gamma}_\bfxi$ is the sum of the residues of the resolvent map $\lambda\mapsto \left(\lambda\I- L_{\bfxi} \right)^{-1}$ at eigenvalues contained inside $\Gamma$. Note moreover that the ranges of $\Pi^{\Gamma}_\bfxi$ and $(\Pi^{\Gamma}_\bfxi)^*$ are both valued in $H^{\infty}_{\rm per}(\R^2;\C^n)$. We refer the reader to \cite[Section~III.6]{Kato} for details concerning the foregoing arguments.

At this stage, we may already sketch a proof of \eqref{e:spec-decomp}. If for any $\bfxi\in[-\pi,\pi]^2$, $\lambda_0\notin \sigma(L_\bfxi)$, it follows by continuity over the compact $[-\pi,\pi]^2$ that $\sup_{\bfxi}\|(\lambda_0\,\I-L_\bfxi)^{-1}\|_{L^2\to L^2}<+\infty$ and, by the Parseval identity, that $\lambda_0\notin \sigma(L)$ with 
\[
(\lambda_0\,\I-L)^{-1}(\bfg)(\bfx)\,=\,\int_{[-\pi,\pi]^{2}} \eD^{\iD \bfxi \cdot \bfx}\,
(\lambda_0\,\I-L_\bfxi)^{-1}(\check{\bfg}(\bfxi,\cdot))(\bfx) \dD\bfxi\,.
\]
In the reverse direction assume that $\lambda_0\in\sigma(L_{\bfxi_0})$ for some $\bfxi_0$ and denote $\bfq^{0}$ a corresponding eigenvector. Then, using again Parseval identities, since $\bfq^{0}\in H^2_{\rm per}([0,1]^2)$, with $\delta>0$ sufficiently small, 
\[
\bfq_\delta^{0}(\bfx):= \int_{[-\pi,\pi]^{2}\cap B(\bfxi_0,\delta)} \eD^{\iD \bfxi \cdot \bfx}\,\bfq^{0}(\bfx) \dD\bfxi
\]
defines a nonzero $\bfq_\delta^{0}\in H^2(\R^2)$ such that
\[
\frac{\|(\lambda_0\,\I-L)\,\bfq_\delta^{0}\|_{L^2(\R^2)}}{\|\bfq_\delta^{0}\|_{L^2(\R^2)}} = \frac{\| {\bf1}_{[-\pi,\pi]^{2}\cap B(\bfxi_0,\delta)}(\bfxi) (L_{\bfxi_0} - L_\bfxi) (\bfq^{0})(\cdot) \|_{L_{\bfxi}^2([-\pi,\pi]^2;L^{2}([0,1]^2))}}{\| {\bf1}_{[-\pi,\pi]^{2}\cap B(\bfxi_0,\delta)}(\bfxi)  \bfq^{0} \|_{L_{\bfxi}^2([-\pi,\pi]^2;L^{2}([0,1]^2))}}
\stackrel{\delta\to0}{\longrightarrow} 0\,.
\]
Hence $\lambda_0\in\sigma(L)$. This concludes the proof of \eqref{e:spec-decomp}. Note that the same arguments apply if, for some $s\in\N$, one considers $L$ as an operator on $H^s(\R^2;\R^n)$ with domain $H^{s+2}(\R^2;\R^n)$ and each $L_\bfxi$ as an operator on $H^s_{\rm per}([0,1]^2;\C^n)$ with domain $H^{s+2}_{\rm per}([0,1]^2;\C^n)$. Incidentally, we observe moreover that it stems from elliptic regularity that the spectrum of each $L_\bfxi$ does not depend on which $H^s_{\rm per}([0,1]^2;\C^n)$ it is considered.

To go further and analyze the implicitly finite-dimensional spectral problems arising from perturbations in $\bfxi$, it is convenient to introduce coordinates. Let $\lambda_0$ be an eigenvalue of $L_{\bfxi_0}$ of multiplicity $m_0$ and $\Gamma_0$ a simple positively-oriented curve such that the intersection of $\sigma(L_{\bfxi_0})$ with its interior is $\{\lambda_0\}$. Pick $(\bfq_j^{(0)})_{1\leq j\leq m_0}$ a basis of the range of $\Pi^{\Gamma_0}_{\bfxi_0}$, and $(\tbq_j^{(0)})_{1\leq j\leq m_0}$ a dual basis of the range of $(\Pi^{\Gamma_0}_{\bfxi_0})^*$. One may extend those to $\bfxi$ near $\bfxi_0$ by 
\begin{align*}
\bfq_j^{\Gamma_0}(\bfxi,\cdot)&\,=\,U_\bfxi^{\Gamma_0}\,\bfq_j^{(0)}\,,&
\tbq_j^{\Gamma_0}(\bfxi,\cdot)&\,=\,((U_\bfxi^{\Gamma_0})^*)^{-1}\,\tbq_j^{(0)}\,,&
j=1,\cdots,m_0\,,
\end{align*}
provided that $(U_\bfxi^{\Gamma_0})_\bfxi$ is a smooth family of bounded invertible operators such that 
\begin{align*}
U_{\bfxi_0}^{\Gamma_0}&\,=\,\I\,,&
U_\bfxi^{\Gamma_0}\,\Pi^{\Gamma_0}_{\bfxi_0}&\,=\,\Pi^{\Gamma_0}_{\bfxi}\,U_\bfxi^{\Gamma_0}\,,\qquad \textrm{for }\bfxi\textrm{ near }\bfxi_0\,.
\end{align*}
Such a family is obtained by setting
\begin{align*}
U_\bfxi^{\Gamma_0}\,:=\,\left(\Pi^{\Gamma_0}_\bfxi\,\Pi^{\Gamma_0}_{\bfxi_0}
+(\I-\Pi^{\Gamma_0}_\bfxi)\,(\I-\Pi^{\Gamma_0}_{\bfxi_0})\right) 
\left(\I-(\Pi^{\Gamma_0}_\bfxi-\Pi^{\Gamma_0}_{\bfxi_0})^2\right)^{-\frac{1}{2}}\,,&
\qquad \textrm{for }\bfxi\textrm{ near }\bfxi_0\,.
\end{align*}
There are various ways to build such a family of operators, we follow here the construction in \cite[Subsection~I-4.6]{Kato} to which we refer for details. Let us simply recall that ${(\cdot)}^{-\frac12}$ is analytic on the open unit ball centered on $\I$ so that the above definition makes sense when $\bfxi$ is sufficiently close to $\bfxi_0$. Note that in this way one obtains for $\bfxi$ near $\bfxi_0$
\begin{align*}
\Pi^{\Gamma_0}_\bfxi&\,=\,\bp \bfq_1^{\Gamma_0}(\bfxi,\cdot)&\cdots&\bfq_{m_0}^{\Gamma_0}(\bfxi,\cdot)\ep\,\bp \langle\tbq_1^{\Gamma_0}(\bfxi,\cdot);\cdot\rangle_{L^2([0,1]^2;\C^n)}\\
\vdots\\
\langle\tbq_{m_0}^{\Gamma_0}(\bfxi,\cdot);\cdot\rangle_{L^2([0,1]^2;\C^n)}\ep\,,\\
L_\bfxi\,\Pi^{\Gamma_0}_\bfxi&\,=\,\bp \bfq_1^{\Gamma_0}(\bfxi,\cdot)&\cdots&\bfq_{m_0}^{\Gamma_0}(\bfxi,\cdot)\ep\,
\bD_\bfxi^{\Gamma_0}
\,\bp \langle\tbq_1^{\Gamma_0}(\bfxi,\cdot);\cdot\rangle_{L^2([0,1]^2;\C^n)}\\
\vdots\\
\langle\tbq_{m_0}^{\Gamma_0}(\bfxi,\cdot);\cdot\rangle_{L^2([0,1]^2;\C^n)}\ep\,,\\
\eD^{t\,L_\bfxi}\,\Pi^{\Gamma_0}_\bfxi&\,=\,\bp \bfq_1^{\Gamma_0}(\bfxi,\cdot)&\cdots&\bfq_{m_0}^{\Gamma_0}(\bfxi,\cdot)\ep\,
\eD^{t\,\bD_\bfxi^{\Gamma_0}}
\,\bp \langle\tbq_1^{\Gamma_0}(\bfxi,\cdot);\cdot\rangle_{L^2([0,1]^2;\C^n)}\\
\vdots\\
\langle\tbq_{m_0}^{\Gamma_0}(\bfxi,\cdot);\cdot\rangle_{L^2([0,1]^2;\C^n)}\ep\,,\qquad t\geq0\,,
\end{align*}
with 
\begin{align*}
\bD_\bfxi^{\Gamma_0}&\,:=\,
\bp \langle\tbq_\ell^{\Gamma_0}(\bfxi,\cdot);\,L_\bfxi\,\bfq_j^{\Gamma_0}(\bfxi,\cdot)\rangle_{L^2([0,1]^2;\C^n)}\ep_{1\leq j,\ell\leq m_0}\,.
\end{align*}

In particular, under Assumption~\cond{2}, we may apply the latter construction with $\bfxi_0=\b0$, $\lambda_0=0$, $m_0=2$, some convenient $\Gamma_0$ symmetric with respect to $0$, $(\bfq_1^{(0)},\bfq_2^{(0)})=(\d_1\ubU,\d_2\ubU)$, and throughout the text we denote
\begin{align*}
\Pi_\bfxi\,,&&\bfq_j^\bfxi(\cdot)\,,&\quad j=1,2\,,&\tbq_j^\bfxi(\cdot)\,,&\quad j=1,2\,,&\bD_\bfxi\,,
\end{align*}
the corresponding objects. Moreover we denote $\Sigma_\bfxi$ the range of $\Pi_\bfxi$ and  $\Sigma_\bfxi^*$ the range of $\Pi_\bfxi^*$. Note that the real symmetry is propagated through the construction, for any $\bfxi$ sufficiently small, $\overline{\Pi_\bfxi}=\Pi_{-\bfxi}$, $\overline{\bD_\bfxi}=\bD_{-\bfxi}$.

\br
Note that there is some freedom in the construction of $\bfq_j^\bfxi(\cdot)$, $\tbq_j^\bfxi(\cdot)$, $j=1,2$, hence of $\bD_\bfxi$. Yet, except in a few places where this is explicitly specified, which particular choice is made is essentially immaterial. A simple fact in this direction is that the first-order expansion of $\bD_\bfxi$ with respect to $\bfxi$ does not depend on this choice. 
\er

\br\label{rk:per-xi}
We have expounded spectral perturbation arguments by varying $\bfxi$ over the compact $[-\pi,\pi]^2$. Yet it is more intrinsic and, for some purposes, also more convenient when $\bfxi_0$ lies on the boundary of $[-\pi,\pi]^2$ in $\R^2$ to consider $\bfxi$ as varying over $\R^2/(2\pi\,\Z)^2$. With this point of view, translation by an element of the lattice $\bfeta\in(2\pi\,\Z)^2$ leaves the spectra invariant but affects generalized eigenspaces according to 
\begin{align*}
\bfq_\ell^{\Gamma_0}(\bfxi+\bfeta,\bfx)&=\eD^{-\iD\,\bfeta\cdot\bfx}\,\bfq_\ell^{\Gamma_0}(\bfxi,\bfx)\,,&
\tbq_\ell^{\Gamma_0}(\bfxi+\bfeta,\bfx)&=\eD^{-\iD\,\bfeta\cdot\bfx}\,\tbq_\ell^{\Gamma_0}(\bfxi,\bfx)\,.
\end{align*}
\er

\subsection{Diffusive stability}\label{ss:diff}

In the present subsection, we investigate equivalent formulations of Assumption~\cond{1}. Our main result is the following proposition.
 
 \bpr
 Assume \cond2. Then \cond1 is equivalent to the union of the following conditions:
\begin{enumerate}
\item[\cond{1'}] For any nonzero $\bfxi$, 
\[
\sigma(L_{\bfxi})\subset\left\{\ \lambda\ ;\ \Re(\lambda)<0\ \right\}\,.
\]
\item[\cond{1W}] The operator $\d_t-\bA(\nabla)$ is hyperbolic and there exists $\theta>0$ and $\xi_0>0$ such that for any $\bfxi\in[-\pi,\pi]^2$ satisfying $\|\bfxi\|\leq\xi_0$,
\[
\sigma(\bDW(\iD\bfxi))\subset\left\{\ \lambda\ ;\ \Re(\lambda)\leq-\theta\|\bfxi\|^2\ \right\}\,
\]
\end{enumerate}
where $\bA$ and $\bDW$ are defined below in \eqref{e:low-expansion}.
\epr

Let us first observe that it follows from standard elliptic estimates that there exist $\omega\in\R$, $C>0$ and $M\in\R^+$ such that for any $\bfxi\in[-\pi,\pi]^2$, 
\[
\sigma(L_\bfxi)\,\subset\,\left\{\,\lambda\in\C
\,;\,|\Im(\lambda)|\,<\,-C\,(\Re(\lambda)-\omega)\,\right\}
\]
and that, when $|\Im(\lambda)|\,\geq\,-C\,(\Re(\lambda)-\omega)$,
\[
\|(\lambda\,\I-L_\bfxi)^{-1}\|\,\leq\,\frac{M}{|\lambda-\omega|}\,.
\]
Combining this with standard analytic semigroup theory and a continuity argument in $\bfxi$ shows that $\cond{1}$ is equivalent to the following assertions (1) and $\cond{1'}$:
\begin{enumerate}
\item[(1)] for any $\xi_0>0$, there exist $\theta>0$ and $C>0$ such that for any $\bfxi\in[-\pi,\pi]^2$ satisfying $\|\bfxi\|\leq \xi_0$ and any $t\geq0$ 
\[
\lnor \eD^{t\,L_{\bfxi}} \rnor\,\leq\,C e^{-\theta t\|\bfxi\|^2}\,;
\]
\item[\cond{1'}]  for any nonzero $\bfxi$, 
\[
\sigma(L_{\bfxi})\subset\left\{\ \lambda\ ;\ \Re(\lambda)<0\ \right\}\,.
\] 
\end{enumerate}
Similarly, the same arguments can be applied to the restriction of $L_\bfxi$ to the range of $(\I-\Pi_\bfxi)$ when $\bfxi$ is sufficiently small. One can show that \cond{2} and \cond{1'} imply that there exist $\xi_0>0$, $\theta>0$ and $C>0$ such that for any $\bfxi\in[-\pi,\pi]^2$  satisfying $\|\bfxi\|\leq \xi_0$ and any $t\geq0$
\[
\lnor \eD^{t\,L_{\bfxi}}\,(\I-\Pi_\bfxi)\rnor\,\leq\,C e^{-\theta t\,\xi_0^2}
\,\leq\,C e^{-\theta t\,\|\bfxi\|^2}\,.
\]
For background on standard analytic semigroup theory used in the foregoing discussion the reader is referred to \cite{Pazy}.

Therefore, assuming~\cond{2}, condition~\cond1 is equivalent to \cond{1'} and
\begin{enumerate}
  \item[\cond{1''}] There exist $\xi_0>0$, $\theta>0$ and $C>0$ such that for any $\bfxi\in[-\pi,\pi]^2$  satisfying $\|\bfxi\|\leq \xi_0$ and any $t\geq0$
\[
\lnor \eD^{t\,\bD_{\bfxi}} \rnor\,\leq\,C e^{-\theta t\|\bfxi\|^2}\,.
\]
\end{enumerate}
We now focus on elucidating \cond{1''}. To do so we introduce the second-order expansion, 
\begin{align}\label{e:low-expansion}
\bD_\bfxi\,\stackrel{\bfxi\to0}{=}\,\underbrace{\bA(\iD\bfxi)\,+\,\bB(\iD\bfxi)}_{\bDW(\iD\bfxi)}\,+\,\cO(\|\bfxi\|^3)\,,
\end{align}
with 
\begin{align*}
\bA(\bfzeta)&\,=\,\bA_1\,\zeta_1+\bA_2\,\zeta_2\,,\\
\bB(\bfzeta)&\,=\,\bB_{1,1}\,\zeta_1^2+2\,\bB_{1,2}\,\zeta_1\,\zeta_2
+\bB_{2,2}\,\zeta_2^2\,,
\end{align*}
where $\bA_1$, $\bA_2$, $\bB_{1,1}$, $\bB_{1,2}$, $\bB_{2,2}$ belong to $\cM_2(\R)$. The next lemma contains a first reduction.

\bl\label{l:robust}
Assume \cond2. Then \cond{1''} is equivalent to 
\begin{enumerate}
\item[\cond{1W''}] There exist $\xi_0>0$, $\theta>0$ and $C>0$ such that for any $\bfxi\in[-\pi,\pi]^2$  satisfying $\|\bfxi\|\leq \xi_0$ and any $t\geq0$
\[
\lnor \eD^{t\,\bDW(\iD\bfxi)} \rnor\,\leq\,C e^{-\theta t\|\bfxi\|^2}\,.
\]
\end{enumerate}
Moreover \cond{1''} implies that $\d_t-\bA(\nabla)$ is an hyperbolic operator.
\el

\begin{proof}
Assuming \cond{1''} with some $(\xi_0,\theta,C)$, one proves \cond{1W''} for $(\tilde \xi_{0},\theta',C')$ where $\theta'$ may be chosen arbitrarily in $(0,\theta)$, $\tilde \xi_{0}$ small enough and $C'$ is tuned accordingly. This follows readily from a Gr\"onwall argument on 
\[
\sup_{0\leq s\leq t}\eD^{\theta'\,s\,\|\bfxi\|^2}\lnor \eD^{s\,\bDW(\iD\bfxi)} \rnor\,
\]
based on 
\[
\eD^{t\,\bDW(\iD\bfxi)} = \eD^{t\,\bD_{\bfxi}}-\int_{0}^{t}\eD^{(t-s)\,\bD_\bfxi} (\bD_{\bfxi}-\bDW(\iD\bfxi))\,\eD^{s\,\bDW(\iD\bfxi)}\,\dD s\,.
\]
The reverse implication is obtained by reversing the roles of $\bDW(\iD\bfxi)$ and $\bD_\bfxi$.

At last, it follows from a similar comparison argument that from \cond{1''} with some $(\xi_0,\theta,C)$, there exists $C,\omega>0$ such that for any $t \geq 0$ and any $\bfxi \in \R^{2}$ with $\|\bfxi\|\leq\xi_0$, 
\[
\lnor \eD^{t\,\bA(\iD\bfxi)} \rnor\ \leq C e^{\omega |\bfxi|^2 t}.
\]
By applying the previous inequality to $(\frac{t}{\eps}, \eps \bfxi)$ for $t>0$ and $\bfxi \in \R^{2}$ and letting $\eps$ goes to $0$, we get the hyperbolicity.
\end{proof}

This first reduction is extremely robust. Now we turn to arguments that use the dimension at hand. In this direction, for comparison, note that in dimension $1$, a first-order constant-coefficient hyperbolic system of two equations is either scalar or strictly hyperbolic. Yet, in general, hyperbolicity is equivalent neither to direction-wise hyperbolicity nor to Friedrichs symmetrizability. Nevertheless this is known to be true for constant-coefficient systems of two equations in arbitrary dimension, see the appendix in \cite{Strang_hyp}. For further related comments and basic background on multidimensional hyperbolic equations we refer the reader to \cite{Benzoni-Serre}.

Our analysis goes further by benefiting from the fact that we have essentially two equations in two dimensions.

\bl\label{l:hyp}
Let $\bA^0_1$, $\bA^0_2\in\cM_2(\R)$, and $\bA^0$ be defined by $\bA^0(\bfzeta):=\bA^0_1\,\zeta_1+\bA^0_2\,\zeta_2$. Then the hyperbolicity of $\d_t-\bA^0(\nabla)$ is equivalent to any of the following conditions
\begin{enumerate}
\item For any unitary $\bfxi_0$, $\d_t-\bA^0(\bfxi_0\,\d_x)$ is hyperbolic.
\item There exists $\bS\in\cM_2(\R)$ such that $\bS$ is symmetric positive definite and $\bS\,\bA^0_1$ and $\bS\,\bA^0_2$ are symmetric.
\item One of the two following conditions holds
\begin{enumerate}
\item For any unitary $\bfxi_0$, $\bA^0(\iD \bfxi_0)$ has real distinct eigenvalues.
\item There exists $\bP\in\cM_2(\R)$ invertible such that $\bP\,\bA^0_1\,\bP^{-1}$ and $\bP\,\bA^0_2\,\bP^{-1}$ are diagonal.
\end{enumerate} 
\end{enumerate}
\el 
 
 \br 
Condition (3)(a) can occur. For instance, one can take for any $\delta \neq 0$,
\begin{align*}
\bA^0_1&\,=\,\bp 0 & \delta \\\delta & 1 \ep \text{ , }&
\bA^0_2&\,=\,\bp 1 & \delta \\\delta & 0 \ep.
\end{align*}
 \er
 \begin{proof}
 The facts that on one hand hyperbolicity implies direction-wise hyperbolicity, and that on the other hand Friedrichs symmetrizable systems, strictly hyperbolic systems or systems of uncoupled scalar equations are indeed hyperbolic are standard elementary parts of the hyperbolic theory. We only need to prove that direction-wise hyperbolicity implies both Friedrichs symmetrizability and the third condition. Thus we assume the first condition.
 
To be more concrete, we introduce coordinates
\begin{align*}
\bA^0_j&\,=\,\bp a_j&b_j\\c_j&d_j\ep\,,\qquad j=1,2\,,&
\bA^0(\bfxi)&\,= \,
\bp \bfa\cdot\bfxi&\bfb\cdot\bfxi\\
\bfc\cdot\bfxi&\bfd\cdot\bfxi\ep\,.
\end{align*}
By elementary considerations the first condition is seen to be equivalent to the fact that for any $\bfxi \in\R^2$ one of the two following conditions holds
\begin{enumerate}[label=(\roman*)]
\item $\bA^0(\bfxi)$ has real distinct eigenvalues;
\item $(\bfa-\bfd)\cdot\bfxi=0$, $\bfb\cdot\bfxi=0$ and $\bfc\cdot\bfxi=0$.
\end{enumerate}
If Condition (ii) of the alternative holds for some nonzero $\bfxi_0$, then $(\bfa-\bfd)$, $\bfb$ and $\bfc$ are colinear, thus for some $\bfell\in\R^2$ and $\bA_0^0\in\cM_2(\R)$, for any $\bfxi \in \R^2$,
\[
\bA^0(\bfxi)\,=\, \tfrac12\,((\bfa+\bfd)\cdot\bfxi)\,\I_2
+ (\bfell\cdot\bfxi)\,\bA_0^0 \,,
\]
and, if $\bfell\neq\b0$, hyperbolicity in the direction $\bfell$ implies that $\bA_0^0$ is diagonalizable with real eigenvalues. Therefore, in this case, independently of whether $\bfell=\b0$ or not, we meet the second part of the third condition of the lemma. It is elementary to check that this implies Friedrichs symmetrizability. Indeed, if $\bP$ is the corresponding diagonalizing matrix, $\bS:=\bP^*\bP$ is a Friedrichs symmetrizer.

The only thing left is to check that in the present case strict hyperbolicity, Condition (i) for any non zero $\bfxi$, implies Friedrichs symmetrizability. In this direction, note first that strict hyperbolicity is equivalently written as, for any nonzero $\bfxi \in \R^{2}$,
\[
(\bfb\cdot\bfxi)\,(\bfc\cdot\bfxi)\,+\,\frac{((\bfa-\bfd)\cdot\bfxi)^2}{4}\,>\,0\,.
\]
This implies that $\bfb$ and $\bfa-\bfd$ are not colinear. Thus there exist real $\alpha$ and $\beta$ such that $\bfc=\alpha\,\bfb+\beta\,(\bfa-\bfd)$. Then, for any $\bfxi$,
\[
(\bfb\cdot\bfxi)\,(\bfc\cdot\bfxi)\,+\,\frac{((\bfa-\bfd)\cdot\bfxi)^2}{4}
\,=\,\alpha\,(\bfb\cdot\bfxi)^2
+\beta\,(\bfb\cdot\bfxi)\,((\bfa-\bfd)\cdot\bfxi)
+\tfrac14\,((\bfa-\bfd)\cdot\bfxi)^2
\]
so that the above sign condition is equivalent to $\alpha>\beta^2$. This implies that
\[
\bS\,:=\,\bp \alpha&-\beta\\-\beta&1\ep
\]
is a Friedrichs symmetrizer.
\end{proof}

From the latter lemmas, one deduces readily the following corollary. In the following, we denote by $\Re(\cdot)$ the self-adjoint part, $\Re(\bM)=(\bM+\bM^*)/2$. 

\bc\label{c:dissipative}
Assume \cond2. Then \cond{1''} (thus also \cond1) implies that $\d_t-\bA(\nabla)$ is Friedrichs symmetrizable. In the reverse direction, if $\bS$ is a Friedrichs symmetrizer of $\d_t-\bA(\nabla)$ such that for any unitary $\bfxi_0 \in \R^{2}$, $\Re(\bS\,\bB(\bfxi_0))$ is positive definite, then \cond{1''} holds. 
\ec

\br\label{rk:symbolic}
In order to apply Corollary~\ref{c:dissipative}, note that in the strictly hyperbolic case, $\bS$ is uniquely determined (up to multiplication by a positive constant) whereas in the case where the system consists in two uncoupled scalar equations, the set of allowed $\bS$ forms a $1$-dimensional family (up to multiplication by a positive constant) if $\bA$ is not identically scalar and a $3$-dimensional family otherwise. One may use the latter freedom to optimize positivity of $\Re(\bS\,\bB(\bfxi_0))$. The same level of $1$-dimensional freedom may also be obtained in the strictly hyperbolic case, provided that one uses symbolic symmetrizers instead of Friedrichs symmetrizer and notices that the corollary also holds for symbolic symmetrizers.
\er

By benefiting from the foregoing considerations, we derive the following lemma about low-frequency diffusivity of second-order systems of two equations.

\bl\label{l:special}
Assume \cond2. Then \cond{1''} is equivalent to any of the following propositions.
\begin{enumerate}
\item The operator $\d_t-\bA(\nabla)$ is hyperbolic and there exists $\theta>0$ and $\xi_0>0$ such that for any $\bfxi\in[-\pi,\pi]^2$ satisfying $\|\bfxi\|\leq\xi_0$,
\[
\sigma(\bDW(\iD\bfxi))\subset\left\{\ \lambda\ ;\ \Re(\lambda)\leq-\theta\|\bfxi\|^2\ \right\}\,.
\]
\item There exists $\theta>0$ such that, for any unitary $\bfxi_0 \in \R^{2}$, 
\begin{enumerate}
\item either $\d_t-\bA(\bfxi_0\,\d_x)$ is strictly hyperbolic, and, for any pair $(\bfell^0,\bV^0)$ such that $\langle\bfell^0;\bV^0\rangle_{\R^2}=1$, with $\bV^0$ an eigenvector of $\bA(\bfxi_0)$ and $\bfell^0$ an eigenvector of $\bA(\bfxi_0)^{\top}$ for the same eigenvalue, we have 
\[
\langle \bfell^0; \bB(\bfxi_0)\,\bV^0\rangle_{\R^2} \geq \theta\,;
\]
\item or $\d_t-\bA(\bfxi_0\,\d_x)$ is scalar and the eigenvalues of $\bB(\bfxi_0)$ have real part larger than $\theta$.
\end{enumerate}
\item There exists $\theta>0$ such that, for any unitary $\bfxi_0$, there exist $\xi_0>0$ and $C>0$ such that for any $\xi\in[0,\xi_0]$ and any $t\geq0$
\[
\lnor \eD^{t\,\bDW(\iD\xi\,\bfxi_0)} \rnor\,\leq\,C e^{-\theta t\,|\xi|^2}\,.
\]
\item There exists $\theta>0$ such that, for any unitary $\bfxi_0$, there exist $\xi_0>0$ and $C>0$ such that for any $\xi\in[0,\xi_0]$ and any $t\geq0$
\[
\lnor \eD^{t\,\bD_{\xi\,\bfxi_0}} \rnor\,\leq\,C e^{-\theta t\,|\xi|^2}\,.
\]
\end{enumerate}
\el

\br 
With this lemma in hand, one can characterize all the operators $\bDW(\nabla)$ satisfying a low-frequency diffusive stability. Note that a bad interaction between the hyperbolic part and the second order part can result in the absence of a low-frequency diffusivity. We illustrate this point with the following matrix
\[
\bD(\iD \bfxi) = \iD \bp \xi_{1} & 0 \\ 0 & 0 \ep + \| \bfxi \|^2 \bp 0 & 1 \\ -1 & -1 \ep
\]
where the second order part is diffusive and yet, $\bDW(\iD \bfxi)$ does not satisfy the second condition of the previous lemma. A spectral perturbation argument at low frequencies reveals that, for $\bfxi \in \R^{2}$ such that $\xi_{1} = \| \bfxi \|$, $\sigma(\bDW(\iD\bfxi)) = \{ \lambda_{1}^{\bfxi} , \lambda_{2}^{\bfxi} \}$ with $\lambda_{1}^{\bfxi} = \iD \| \bfxi \| + \iD \| \bfxi \|^{3} - \| \bfxi \|^{4} + \mathcal{O}(\| \bfxi \|^{5})$ and $\lambda_{2}^{\bfxi} =  - \| \bfxi \|^{2} + \mathcal{O}(\| \bfxi \|^{3})$ so that it is clear that \cond{1W''} is not satisfied.
 \er 
\begin{proof}
That \cond{1''} implies the first condition stems from Lemma~\ref{l:robust}. Now we show that the first condition implies the second. It is classical that hyperbolicity implies direction-wise hyperbolicity. 

Assume first that $\d_t-\bA(\bfxi_0\,\d_x)$ is strictly hyperbolic. Then, when $\xi$ is small, $\bDW(\iD\,\xi\,\bfxi_0)$ is smoothly diagonalizable with simple eigenvalues, with eigenvectors perturbing from those of $\bA(\iD \bfxi_0)$ and eigenvalues expanding as 
\[
\lambda\,\stackrel{\xi\to0}{=}\, \iD \lambda_0\, \xi\,-\,\xi^2\,\langle \bfell^0; \bB(\bfxi_0)\,\bV^0\rangle_{\R^2}\,+\,\cO(|\xi|^3)\,,
\]
where $\lambda_0$ is an eigenvalue of $\bA(\bfxi_0)$ and $(\bfell^0,\bV^0)$ is an associated dual pair of left-right eigenvectors. Therefore the first condition implies the second in the direction $\bfxi_0$ with a uniform $\theta$.

Assume now that $\d_t-\bA(\bfxi_0\,\d_x)$ is scalar, with characteristic speed $\bfc^0$. Then $\bDW(\iD\,\xi\,\bfxi_0)=\iD\,\xi \, \bfc^0\,-\,\xi^2\,\bB(\bfxi_0)$. It is clear that that in this case also the first condition implies the second in the direction $\bfxi_0$ with a uniform $\theta$.

Now we aim at proving that the second condition imply \cond{1''}. Thanks to a compactness argument and Lemma \ref{l:robust}, it is sufficient to prove that \cond{1W''} holds in the neighborhood of any direction $\bfxi_0$. If the strictly-hyperbolic part of the second condition holds at $\bfxi_0$, one can introduce $\bP\in\cM_2(\R)$ invertible that diagonalizes $\bA(\bfxi_0)$ and easily get eigenvalue expansions of $\bP \eD^{t\bDW(\iD\bfxi)} \bP^{-1}$ when $\tfrac{\bfxi}{\| \bfxi \|}$ is sufficiently close to $\bfxi_0$ and $\| \bfxi \|$ is small enough. It is also quite immediate if $\d_t-\bA(\nabla)$ is scalar, since then one may use dissipative symbolic symmetrizers adapted to the second-order part. Thanks to Lemma~\ref{l:hyp}, this means that we are only left with the analysis of the case when $\d_t-\bA(\nabla)$ is scalar in the direction $\bfxi_0$ but strictly hyperbolic in nearby directions. After a global diagonalization and a change of coordinates, one may assume that $\bfxi_0=\beD_1$, $\bA_{1}$ and $\bA_2$ are diagonal and that 
\[
\bB_{1,1}\,=\,\bp \alpha&\beta\\\gamma&\delta\ep \text{ , }
\]
with $\alpha > 0$, $\delta > 0$ and $\alpha\,\delta>\beta\,\gamma$. Note that we have used here that there is a uniform spectral gap for nearby directions to determine the signs of $\alpha$ and $\delta$, since otherwise the assumption in the direction $\bfxi_0$ would only yield $\alpha+\delta>0$ and $\alpha\,\delta>\beta\,\gamma$. As in Corollary~\ref{c:dissipative}, it is sufficient to find $\theta'>0$ and a real symmetric positive definite $\bS$ such that $\bS\,\bA(\bfxi)$ is symmetric and $\Re(\bS\,\bB(\bfxi))\geq \theta'$ when $\bfxi$ is unitary and sufficiently close to $\bfxi_0$. If $\beta\,\gamma\neq 0$, 
\[
\bS:=\bp \sqrt{\frac{|\gamma|}{|\beta|}}&0\\0&\sqrt{\frac{|\beta|}{|\gamma|}}\ep
\] 
fits the requirement. If $\beta=0$ (resp. $\gamma=0$), then $\alpha >0$ and $\delta >0$ and 
\begin{align*}
\bS&:=\bp M&0\\0&1\ep\,,&
\Big(\textrm{resp. }\bS&:=\bp 1&0\\0&M\ep\,,\Big)
\end{align*}
with $M>0$ sufficiently large, does the job. Note that when checking the requirements we are using that when $\alpha_1>0$, $\alpha_4>0$, $\alpha_2\in\R$, $\alpha_3\in\R$, the matrix
\[
\bp \alpha_1&\alpha_2\\\alpha_3&\alpha_4 \ep
\] 
has positive real part provided that $\alpha_1\,\alpha_4>\left(\frac{\alpha_2+\alpha_3}{2}\right)^2$. This achieves the proof that the second condition implies \cond{1''}.

The second condition, being direction-wise, is clearly equivalent to the third, whereas the equivalence of the third and fourth conditions follow by a perturbation argument as in Lemma~\ref{l:robust}.

This concludes the proof. 
\end{proof}

\section{Profile equations}\label{s:profile}

\subsection{Local structure}\label{ss:profile}

The present subsection is devoted to the proof of the following proposition ensuring that Assumption~\cond{2} encodes sufficient information to elucidate the structure of nearby periodic waves.

\bpr\label{p:structure}
Assume~\cond{2}. Then there exist $\eps_0>0$ and a smooth map 
\[
B(\ubK,\eps_0)\to H^2_{\rm per}([0,1]^2;\R^n)\times\R^2\,,\quad
\bK\mapsto (\bU^\bK(\cdot),\bfc(\bK))
\]
such that, for any wavematrix $\bK\in B(\ubK,\eps_0)$, $(\bK,\bU^\bK(\cdot),\bfc(\bK))$ solves \eqref{stand_eq} and for any $(\bU,\bfc)\in H^2_{\rm per}([0,1]^2;\R^n)\times\R^2$ such that $(\bK,\bU,\bfc)$ solves \eqref{stand_eq} and 
\begin{align*}
\|\bfc-\ubc\|&\leq \eps_0\,,&
\inf_{\bfvarphi_0\in\R^2}\|\bU-\ubU(\cdot+\bfvarphi_0)\|_{H^2_{\rm per}([0,1]^2;\R^n)}&\leq \eps_0\,,
\end{align*}
one has $\bfc=\bfc(\bK)$ and there exists $\bfvarphi\in\R^2$ such that $\bU=\bU^\bK(\cdot+\bfvarphi)$. Moreover the map $\bK\mapsto \bU^\bK$ is valued in $H^\infty_{\rm per}([0,1]^2;\R^n)$ and, for any $s\in\N$, there exists $0<\eps_0'\leq\eps_0$ such that it is smooth as a map from $B(\ubK,\eps_0')$ to $H^s_{\rm per}([0,1]^2;\R^n)$. 
\epr

The proof follows the Lyapunov-Schmidt reduction. We first show that we can factor out translational invariance. To do so, we may apply the Implicit Function Theorem to the map 
\[
H^2([0,1]^2;\R^n)\times\R^2\to\R^2\,,\qquad(\bU,\bfvarphi)\longrightarrow
(\langle\tbq_j^{(0)};\bU(\cdot-\bfvarphi)\rangle_{L^2([0,1]^2;\R^n)})_{j=1,2}
\]
in a neighborhood of $(\ubU,\b0)$, where  $(\tbq_1^{(0)},\tbq_2^{(0)})$ is the basis of $\Sigma_0^*$ in duality with $(\bfq_1^{(0)},\bfq_2^{(0)})=(\d_1\ubU,\d_2\ubU)$. Indeed the map is $\cC^1$ and, at $(\ubU,\b0)$, its differential map with respect to $\bfvarphi$ is $-\I$. By using translational invariance, this implies that there exist $\eps>0$ and $C>0$ such that if $(\bU,\bfvarphi_0)$ is such that
\[
\|\bU-\ubU(\cdot+\bfvarphi_0)\|_{H^2_{\rm per}([0,1]^2;\R^n)}\leq \eps
\]
then there exists $\bfvarphi$ such that $\tbU=\bU(\cdot-\bfvarphi)$ satifies
\begin{align*}
\|\tbU-\ubU\|_{H^2_{\rm per}([0,1]^2;\R^n)}&\leq C\,\|\bU-\ubU(\cdot+\bfvarphi_0)\|_{H^2_{\rm per}([0,1]^2;\R^n)}\,,&
(\I-\Pi_0)\,(\tbU-\ubU)&=0\,.
\end{align*}
It is thus sufficient to prove a genuine uniqueness under the assumption that $\bU-\ubU$ is small and $(\I-\Pi_0)\,(\bU-\ubU)=0$.

Let us denote by $L_{\b0}^\dagger$ the inverse of $L_{\b0}$ restricted to the range of $(\I-\Pi_0)$. With the extra constraint $(\I-\Pi_0)\,(\bU-\ubU)=0$, Equation~\eqref{stand_eq} is equivalent to
\begin{align*}
\bU&\,=\,\ubU
-L_{\b0}^\dagger[(\I-\Pi_0)\,\cR]\\
\bfc&\,=\,\transp{(\bK^{-1}\ubK)}\ubc
-\transp{(\bK^{-1})}
\,\bp \langle\tbq_1^{(0)};\cdot\rangle_{L^2([0,1]^2;\R^n)}\\
\langle\tbq_2^{(0)};\cdot\rangle_{L^2([0,1]^2;\R^n)}\ep\,\cR
\end{align*}
with
\begin{align*}
\cR&=\left(\transp{(\bK\nabla)}(\bK\nabla)-\transp{(\ubK\nabla)}(\ubK\nabla)\right)\bU
+((\transp{\bK}\bfc-\transp{\ubK}\ubc)\cdot \nabla)(\bU-\ubU)\\
&\quad+\transp{(\bK\nabla)}\bG(\bU)
-\transp{(\ubK\nabla)}\left(\bG(\ubU)+\dD \! \bG(\ubU)(\bU-\ubU)\right)\\
&\quad+\bff(\bU)-\bff(\ubU)-\dD \! \bff(\ubU)(\bU-\ubU))\,.
\end{align*}
The proof is then achieved by another application of the Implicit Function Theorem.

\subsection{Profile variations}\label{ss:algebraic}

We collect here some algebraic relations between the expansions of $L_\bfxi$ when $\bfxi$ is small and the derivatives of wave profiles, obtained by differentiating profile equation \eqref{stand_eq}. 

To prepare comparisons, we expand Bloch symbols $L_\bfxi$ as
\be\label{def:Bloch-expand}
L_\bfxi \bV\ =\ L_{\b0}\bV\ +\ \transp{(L^{(1)}\bV)}\iD(\ubK\bfxi)
\ -\ \|\ubK\bfxi\|^2\, \bV
\ee
where
\[
L^{(1)}\bV \ :=\ 2 \ubK\nabla \bV + \dD \! \bG(\ubU)(\bV) + \ubc \transp{\bV}.
\]

By design and invariance by translation, for any $\bfvarphi$, $(\bK,\bU^\bK(\cdot+\bfvarphi),\bfc(\bK))$ solves \eqref{stand_eq}. Differentiating this with respect to $\bfvarphi$ gives
\be\label{e:dx}
 L_{\b0}\,\d_j\ubU=0\,,\qquad j=1,2,
\ee
while differentiating it with respect to $\bK_j$ in the direction $\bfeta$ leads to 
\be\label{e:dk}
L_{\b0}[\dD_{\bK_j}\bU (\ubK)(\bfeta)]\,+\,\transp{(L^{(1)} \d_j\ubU)}\bfeta + (\transp{\ubK}\dD_{\bK_j}\bfc(\ubK)(\bfeta) \cdot \nabla)\ubU =0\,,\quad\textrm{for }j=1,2.
\ee
Finally, with $\bOm(\bK):=-\transp{\bK}\bfc(\bK)$, differentiating the same relation first with respect to $\bK_j$ in the direction $\bfeta$ then with respect to $\bK_\ell$ in the direction $\bfzeta$ leads to 
\ba\label{e:d2k}
L_{\b0}[\dD^{2}_{\bK_j,\bK_\ell}&\bU(\ubK)(\bfeta,\bfzeta)] 
+ \dD^{2} \! \bff (\bU)(\ubK)(\dD_{\bK_j}\bU(\ubK)(\bfeta),\dD_{\bK_\ell} \bU(\ubK)(\bfzeta))\\
&+ \transp{(\ubK \nabla)} \dD^{2} \! \bG(\ubU)  (\dD_{\bK_j}\bU(\ubK)(\bfeta) , \dD_{\bK_\ell} \bU(\ubK)(\bfzeta))\\
&+\d_{\ell} \transp{\left( \dD \! \bG(\ubU)  (\dD_{\bK_j} \bU(\ubK)(\bfeta) ) \right)} \bfzeta + \d_j \transp{\left( \dD \! \bG(\ubU)  (\dD_{\bK_\ell} \bU(\ubK)(\bfzeta) ) \right)} \bfeta\\
&+2 (\transp{\ubK}\bfzeta\cdot\nabla) \d_{\ell}(\dD_{\bK_j} \bU(\ubK)(\bfeta)) 
+ 2 (\transp{\ubK} \bfeta\cdot\nabla) \d_{j}(\dD_{\bK_\ell} \bU(\ubK)(\bfeta))+ 2 (\transp{\bfeta} \bfzeta)\,\d_{j\ell}^{2} \ubU\\
&=\left( \dD_{\bK_j } \bOm(\ubK)(\bfeta)\cdot\nabla \right) \dD_{\bK_\ell} \bU(\ubK)(\bfzeta)
+\left( \dD_{\bK_\ell} \bOm(\ubK)(\bfzeta) \cdot \nabla \right) \dD_{\bK_j} \bU(\ubK)(\bfeta)\\
&\quad+\left(\dD^{2}_{\bK_j,\bK_\ell} \bOm(\ubK)(\bfeta,\bfzeta) \cdot\nabla \right) \ubU\,,\hspace{8em} j=1,2\,,\ \ell=1,2\,.
\ea

\section{Phase related estimates}\label{s:phase-like}

We gather here some estimates associated with the presence of a spatially-dependent phase modulation.

\subsection{Sobolev-like estimates}\label{ss:Sobolev}

In the present subsection, we bound $\bfphi$ and $\nabla\bfphi$ in terms of $\Delta\bfphi$. 

Up to immaterial perpendicular rotations, the reconstruction of $\nabla\bfphi$ from $\bfphi$ coincides with the Biot-Savart law that recovers divergence-free vector fields from their curl. For this reason, the estimates gathered here are essentially special cases of harmonic analysis estimates commonly used in the analysis of incompressible fluid mechanics. However, for the sake of consistency with the rest of our analysis, we have decided to provide simplified versions of the latter so as to prove them solely from Young and Hausdorff-Young inequalities and elliptic regularity in Calder\'on-Zygmund form. For sharper estimates, we refer the reader to either \cite[Annexe~C]{Rodrigues-these} or \cite[Section~1.1]{Rodrigues-density-dependent}. Yet we do make some comments and remarks involving more advanced functional spaces and we refer the reader to \cite[Part~1]{LM_NS} for the necessary background.

Obviously, bounds are affordable only in regimes where we already know that $\Delta\bfphi$ determines $\nabla\bfphi$ or even $\bfphi$. When $\bfphi$ is a tempered distribution, $\Delta\bfphi$ determines $\bfphi$ up to a polynomial and a further condition is needed to ensure uniqueness. In our case, the reconstruction implicitly hinges on the uniqueness result that the only harmonic tempered distribution that belongs to $\Span\left(\bigcup_{\substack{1\leq p<\infty\\1\leq q\leq\infty}}L^{p,q}(\R^2)\right)$ is the zero function. In the foregoing, $L^{p,q}$ denotes Lorentz spaces, whose Lebesgue spaces $L^p=L^{p,p}$ are special cases, and the uniqueness follows from the fact that no nonzero polynomial belongs to the latter span. As a consequence, when $\bfphi$ is a tempered distribution,
\begin{itemize}
\item the knowledge of $\Delta\bfphi$ and the condition $\nabla\bfphi\in\Span\left(\bigcup_{\substack{1\leq p<\infty\\1\leq q\leq\infty}}L^{p,q}(\R^2;\cM_2(\R))\right)$ determines $\nabla\bfphi$ completely thus it also determines $\bfphi$ up to a constant function; 
\item for any fixed $\bfphi_0$, the knowledge of $\Delta\bfphi$ and the condition $\bfphi-\bfphi_0\in\Span\left(\bigcup_{\substack{1\leq p<\infty\\1\leq q\leq\infty}}L^{p,q}(\R^2;\R^2)\right)$ determines $\bfphi$ thus also $\nabla\bfphi$.
\end{itemize}
The first extra condition is consistent with the way we recover $\nabla\bfphi$ since it ensures that $\Delta\bfphi\in L^1$ implies $\nabla\bfphi \in L^{2,\infty}(\R^2;\cM_2(\R))$. Moreover, it follows from Propositions~\ref{pr:lin-local} and~\ref{pr:lin-mod} and Lemma~\ref{l:time-bl} that both extra conditions are propagated by the time evolution.

Since the reconstruction is done component-wise, we may reduce to the consideration of a scalar $\phi$. We begin by recovering $\nabla\phi$ from $\Delta\phi$. Given some $d\in L^1(\R^2)$, we define $\bfv:=\nabla\Delta^{-1}d$ such that $\Curl \bfv\equiv 0$, $\Div\bfv=d$, by
\begin{align}\label{eq:BS-Fourier}
\widehat{\bfv}(\bfxi)\,=\,\frac{\iD\,\bfxi}{\|\bfxi\|^{2}} \widehat{d}(\bfxi)\,,
\end{align}
or equivalently through
\begin{align}\label{eq:BS}
\bfv&:=\cG\star d\,,&
\cG(\bfy)&=\frac{1}{2\pi}\frac{\bfy}{\|\bfy\|^2}\,,\quad\bfy\in\R^2\,.
\end{align}

\bpr\label{p:BS}
\begin{enumerate}
\item There exist $(C_p)_{2<p<\infty}$ such that for any $d\in (L^1\cap L^2)(\R^2)$, $\bfv$ defined by \eqref{eq:BS-Fourier} lies in $\cap_{2<p<\infty}L^p(\R^2;\R^2)$ and
\begin{align*}
\|\bfv\|_{L^p(\R^2;\R^2)}&\leq\,C_p\,\|d\|_{(L^1\cap L^2)(\R^2)}\,,
&2<p< \infty\,.
\end{align*}
\item For any $1\leq p_0<2<p_1\leq\infty$, there exists $C_{p_0,p_1}$ such that $\bfv$ defined by \eqref{eq:BS-Fourier} from $d$ satisfies 
\begin{align*}
\|\bfv\|_{L^\infty(\R^2;\R^2)}&\leq\,C_{p_0,p_1}\,\|d\|_{(L^{p_0}\cap L^{p_1})(\R^2)}\,.
\end{align*}
\item For any $1<p<\infty$, there exists $C_p$ such that $\bfv$ defined by \eqref{eq:BS-Fourier} from $d$ satisfies 
\begin{align*}
\|\nabla\bfv\|_{L^p(\R^2;\cM_2(\R))}&\leq\,C_p\,\|d\|_{L^p(\R^2)}\,.
\end{align*}
\item Assume that $\|\,\cdot\,\|\,d\in L^1(\R^2)$ and $d\in L^p(\R^2)$ for some $p>1$. Then $\bfv$ defined by \eqref{eq:BS-Fourier} from $d$ belongs to $L^2(\R^2;\R^2)$ if and only if $\int_{\R^2}d=0$.
\end{enumerate}
\epr

We prove the last point only to justify a remark of the introduction.

\begin{proof}
Since $\nabla \bfv=\nabla\nabla \Delta^{-1} d$, the third estimate stems directly from Calder\'on-Zygmund theory that ensures that $\nabla (-\Delta)^{-\frac12}$ acts boundedly on $L^p$, $1<p<\infty$.

The first estimate follows from Hausdorff-Young and H\"older inequalities since $\bfxi\mapsto \|\bfxi\|^{-1}$ belongs to $\left(\cap_{2<p\leq\infty}L^p\right)+\left(\cap_{1\leq p<2}L^p\right)$ and $\|\,\widehat{d}\ \|_{L^2\cap L^\infty}\lesssim \|d\|_{L^1\cap L^2}$. The second estimate follows from Young inequality since $\cG\in \left(\cap_{2<p\leq\infty}L^p\right)+\left(\cap_{1\leq p<2}L^p\right)$.

To prove the last point, we first observe that for any $1\leq q\leq p$, 
\[
\|d\|_{L^q(\R^2)}
\lesssim
\|d\|_{L^p(\R^2)}+\||\cdot|\,d\|_{L^1(\R^2)}\,,
\]
so that we may use any of these norms in the argument. Now, it follows from Hausdorff-Young and H\"older inequalities that the $L^2$ norm of the high-frequency part of $\bfv$ is controlled by $\|d\|_{L^{\min(\{p,2\})}(\R^2)}$. Moreover, from the pointwise bound
\[
\left\|\widehat{\bfv}(\bfxi)-\frac{\iD\,\bfxi}{\|\bfxi\|^{2}} \widehat{d}(0)\right\|
\lesssim \|\nabla_\bfxi\,\widehat{d}\ \|_{L^\infty(\R^2;\R^2)}
\lesssim \||\cdot|\,d\|_{L^1(\R^2)}
\]
one deduces that $\widehat{\bfv}$ is locally square-integrable if and only if $\widehat{d}(0)=0$. Hence the result.
\end{proof}

We now turn to the reconstruction of $\phi$ from $\Delta\phi$. Given some $d\in L^1(\R^2)$, we define $\phi:=\Delta^{-1}d$ such that $\Delta\phi=d$ by
\begin{align}\label{eq:BS0-Fourier}
\widehat{\phi}(\bfxi)\,=\,{\textrm p.v.}\left(\frac{1}{\|\bfxi\|^{2}}\right)\widehat{d}(\bfxi)\,,
\end{align}
or equivalently through
\begin{align}\label{eq:BS0}
\phi&:=\cG_0\star d\,,&
\cG_0(\bfy)&=\frac{1}{4\pi}\ln\left(\|\bfy\|^2\right)\,,\quad\bfy\in\R^2\,.
\end{align}
Note that conventions are consistent in the sense that $\nabla(\Delta^{-1}d)=(\nabla\Delta^{-1})d$.

\bpr\label{p:BS0}
\begin{enumerate}
\item Assume that $\chi$ is a smooth compactly supported functions equal to $1$ in a neighborhood of $\b0$. Then there exist $(C_{p,q})_{\substack{1\leq q\leq p\leq\infty\\(q,p)\neq (1,\infty)}}$ such that if $\phi$ is defined by \eqref{eq:BS0-Fourier} from $d$ and $\phi_{HF}$ is defined by $\widehat{(\phi_{HF})}=(1-\chi)\,\phi$, there holds
 \begin{align*}
&\|\phi_{HF}\|_{L^p(\R^2)}\leq\,C_{p,q}\,\|d\|_{L^q(\R^2)}\,,
&1\leq q\leq p\leq\infty\,,&\quad(q,p)\neq (1,\infty)\,.
\end{align*}
\item For any $1<p_0\leq\infty$ and $1<p_1\leq\infty$, there exists $C_{p_0,p_1}$ such that $\phi$ defined by \eqref{eq:BS0-Fourier} from $d$ satisfies
\begin{align*}
\big\|\phi\,-\,\left(\int_{\R^2}d\right)\,\frac{1}{4\pi}\ln_+\left(\|\cdot\|^2\right)\Big\|_{L^\infty(\R^2)}
\leq C_{p_0,p_1}\,\left(
\|d\|_{L^{p_0}(\R^2)}+\left\|\|\cdot\|^{2\,\left(1-\tfrac{1}{p_1}\right)}\,\left(\ln(2+\|\cdot\|)\right)^3\,d\right\|_{L^{p_1}(\R^2)}\right)
\end{align*}
and, in particular, $\phi$ belongs to $L^\infty(\R^2)$ if and only if $\int_{\R^2}d=0$.
\end{enumerate}
\epr

Here also, we prove the last point only to justify a remark of the introduction.

\begin{proof}
Since $\bfxi\mapsto \|\bfxi\|^{-2}\,(1-\chi(\bfxi))$ belongs to $\bigcap_{1<q\leq 2}W^{2,q}$, its inverse Fourier transform belongs to $\bigcap_{2\leq r<\infty}L^r((1+\|\cdot\|^2))$ thus to $\bigcap_{1\leq s<\infty}L^s$. Therefore the first estimate follows from Young inequalities. 

We now turn to the second bound. We first observe that 
\[
\|d\|_{L^1}\lesssim \|\ln(2+\|\cdot\|)\,d\|_{L^1}
\lesssim \|d\|_{L^{p_0}}+\left\|\|\cdot\|^{2\,\left(1-\tfrac{1}{p_1}\right)}\,\left(\ln_+\|\cdot\|\right)^3\,d\right\|_{L^{p_1}}
\]
When $\|\bfx\|\leq 1$, splitting
\[
\phi(\bfx)\,=\,\int_{\R^2}\,\cG_0(\bfy)\,d(\bfx-\bfy)\,\dD\bfy
\]
according to whether $\|\bfy\|\leq 1$ or $\|\bfy\|\geq 1$, one deduces that
\[
|\phi(\bfx)|\,\lesssim\,\|d\|_{L^{p_0}}+\|\ln(2+\|\cdot\|)\,d\|_{L^1}
\]
since $\cG_0$ belongs locally to any $L^q$ with $1\leq q<\infty$. We now assume that $\|\bfx\|\geq 1$ and split the integral in 
\[
\phi(\bfx)-\,\left(\int_{\R^2}d\right)\,\frac{1}{4\pi}\ln\left(\|\bfx\|^2\right)\,=\,\frac{1}{2\pi}\int_{\R^2}\,\ln\left(\frac{\|\bfy\|}{\|\bfx\|}\right)\,d(\bfx-\bfy)\,\dD\bfy
\]
according to $\|\bfy\|\leq \tfrac12\|\bfx\|$, $\tfrac12\|\bfx\|\leq\|\bfy\|\leq 2\|\bfx\|$ or $2\|\bfy\|\leq \|\bfx\|$. The contribution from $\tfrac12\|\bfx\|\leq\|\bfy\|\leq 2\|\bfx\|$ is bounded by a multiple of $\|d\|_{L^1}$. In the regime $\|\bfy\|\leq \tfrac12\|\bfx\|$, we also have $\|\bfx-\bfy\|\geq \tfrac12\|\bfx\|$ so that its contribution is bounded by a multiple of $\left\|\|\cdot\|^{2\,\left(1-\tfrac{1}{p_1}\right)}\,d\right\|_{L^{p_1}}$ since
\[
\Big\|\ln\left(\frac{\|\cdot\|}{\|\bfx\|}\right)\Big\|_{L^{p_1'}(\|\cdot\|\leq \tfrac12\|\bfx\|)}
\,\leq\,\|\bfx\|^{2\,\left(1-\tfrac{1}{p_1}\right)}\,.
\]
At last, in the regime $\|\bfy\|\geq 2\|\bfx\|$, we also have $\|\bfx-\bfy\|\geq \tfrac12\|\bfy\|$ so that this contribution is bounded by a multiple of $\left\|\|\cdot\|^{2\,\left(1-\tfrac{1}{p_1}\right)}\,\left(\ln_+\|\cdot\|\right)^3\,d\right\|_{L^{p_1}}$ since
\[
\Big\|\|\cdot\|^{-2\,\left(1-\tfrac{1}{p_1}\right)}\,
\left(\ln\|\cdot\|\right)^{-2}\Big\|_{L^{p_1'}(\|\cdot\|\geq 1)}
\,<+\infty\,.
\]
\end{proof}

Note that an argument similar to the one used to prove the first estimate yields for any $\phi$, if we denote $\phi_{LF}:=\phi-\phi_{HF}$,
\begin{align}\label{r:estim_lowfreq}
\| \phi_{LF} \|_{L^q}
&\lesssim\|\phi\|_{L^p}
&1\leq p\leq q\leq\infty\,,
\end{align}
that we use without mention throughout the text.

The condition $d\in L^1$ only ensures $\phi=\Delta^{-1}d\in BMO$. Two classical ways to restore $\phi\in L^\infty$ are to assume that $d$ belongs either to the real Hardy space $\cH^1$ or to the homogeneous Besov space $\dot{B}^{0}_{1,1}$. We stress that this is consistent with Proposition~\ref{p:BS0} since both $\cH^1$ and $\dot{B}^{0}_{1,1}$ are included in the space of integrable functions with zero integral.

\br\label{r:sobolev_gradient}
We also have:
\begin{enumerate}
\item For any $r \in [1,2[$, there exists a constant $C_{r}$ such that for any locally integrable function $\phi$, there exists a constant $\phi_{\infty}$ so that
\[
\| \phi - \phi_{\infty} \|_{L^{\frac{2r}{2-r}}(\R^{2})} \leq C_{r} \| \nabla \phi \|_{L^{r}(\R^{2})}.
\]
\item there exist $(C_p)_{2 \leq p<\infty}$ such that for any function $\phi$ that vanishes at infinity
\begin{align*}
\|\phi \|_{L^p(\R^2)}&\leq\,C_p\,\| \nabla \phi \|_{(L^1\cap L^2)(\R^2;\R^{2})}\,,
&2 \leq p< \infty\,,
\end{align*}
\item for any $1\leq p_0<2<p_1\leq\infty$, there exists $C_{p_0,p_1}$ such that for any function $\phi$ that vanishes at infinity
\begin{align*}
\| \phi \|_{L^\infty(\R^2)}&\leq\,C_{p_0,p_1}\,\| \nabla \phi \|_{(L^{p_0}\cap L^{p_1})(\R^2;\R^{2})}\,,
\end{align*}
\item there exists $C>0$, such that for any function $\bfphi$ that vanishes at infinity 
 \begin{align*}
&\|\phi_{HF}\|_{L^q(\R^2)}\leq\,C\, \| \nabla \phi \|_{L^q(\R^2)}\,,
&1\leq q\leq \infty\,.
\end{align*}
\end{enumerate}
The proofs of the last three point is similar to the previous ones. The main new ingredient is the fact that 
\[
\widehat{\phi} (\bfxi) = - \iD \sum_{j} \frac{\xi_{j}}{\| \bfxi \|^{2}} \widehat{\partial_{j} \phi} (\bfxi) := - \iD \sum_{j} h_{j}(\bfxi) \widehat{\partial_{j} \phi} (\bfxi)
\]
where $h_{j}$ belongs to $\left(\cap_{2<p\leq\infty}L^p\right)+\left(\cap_{1\leq p<2}L^p\right)$  and $\nabla_{\bfxi} \left( (1-\chi) h_{j} \right)$ belongs to $\bigcap_{1<q \leq 2}W^{1,q}$ so that the inverse Fourier transform of $(1-\chi) h_{j}$ belongs to $\bigcap_{2\leq r<\infty}L^r(( \|\cdot\|+\|\cdot\|^2))$ thus to $\bigcap_{1\leq s<2}L^s$.

The first point is a consequence of an homogeneous Poincar\'e-Wirtinger type inequality.
\er

\subsection{Change of variables}\label{ss:var-change}

We store here a basic estimate, useful to invert $\Id-\bfphi(t,\cdot)$ and quantify its impact on bounds. It is almost identical to the first half \cite[Lemma~2.7]{JNRZ-conservation}. 

\bl\label{l:var-change}
Assume that $\bfphi\,:\ \R^2\to\R^2$ is a Lipschitz function such that $\|\nabla\bfphi\|_{L^\infty(\R^2;\R^2)}<1$. Then $\Id-\bfphi$ is invertible and for any $1 \leq p \leq \infty$,
\begin{align*}
\|A-B\circ(\Id-\bfphi)^{-1}\|_{L^p}
&\leq (1+\|\nabla\bfphi\|_{L^\infty(\R^2;\R^2)})^{\frac{2}{p}}
\,\|A\circ(\Id-\bfphi)-B\|_{L^p}\,,\\
\|A\circ(\Id-\bfphi)-B\|_{L^p}
&\leq \frac{1}{(1-\|\nabla\bfphi\|_{L^\infty(\R^2;\R^2)})^{\frac{2}{p}}}
\,\|A-B\circ(\Id-\bfphi)^{-1}\|_{L^p}\,.
\end{align*}
Finally, if $\bfphi_{1}, \bfphi_{2}\,:\ \R^2\to\R^2$ are  Lipschitz functions such that we have $\|\nabla\bfphi_{1} \|_{L^\infty(\R^2;\R^2)}<1$ and $\|\nabla\bfphi_{2} \|_{L^\infty(\R^2;\R^2)}<1$, then for any $1 \leq p \leq \infty$,
\[
\|(\Id-\bfphi_{1})^{-1}-(\Id-\bfphi_{2})^{-1} \|_{L^p}
\leq \frac{(1+\|\nabla\bfphi\|_{L^\infty(\R^2;\R^2)})^{\frac{2}{p}}}{1-\|\nabla\bfphi\|_{L^\infty(\R^2;\R^2)}} \|\bfphi_{1}-\bfphi_{2} \|_{L^p}.
\]
\el

\begin{proof}
The invertibility is a direct consequence of the Banach fixed point argument. The first two estimates follow from a change of variable. The last one is a consequence of the equality
\[
(\Id-\bfphi_{1})^{-1}-(\Id-\bfphi_{2})^{-1} = (\bfphi_{1}-\bfphi_{2}) \circ (\Id-\bfphi_{1})^{-1} + \bfphi_{2} \circ (\Id-\bfphi_{1})^{-1} - \bfphi_{2} \circ (\Id-\bfphi_{2})^{-1}.
\]
\end{proof}

The second --- and less trivial --- half of \cite[Lemma~2.7]{JNRZ-conservation}, estimating $\|A-B\circ(\Id+\bfphi)\|_{L^p}$ in terms of $\|A\circ(\Id-\bfphi)-B\|_{L^p}$ and $\|\nabla\bfphi\|_{L^p}$, is of no use here because it requires $\bfphi$ to be bounded.

Note that, though we do not bother to state those, it is clear from the proof that variants involving regularity in $\bfx$ or in passive variables --- such as $t$ --- also hold.

\section{Geometrical optics}\label{s:WKB}

In the present appendix, we show how to guess from formal geometrical optics considerations the conclusions about modulation behavior and averaged dynamics that our analysis proves rigorously following different paths. 

\subsection{Formal derivation of averaged equations}\label{ss:formal}

To begin our formal process, let us consider the slow/fastly-oscillatory \emph{ansatz}  
\be\label{e:ansatz}
\bcW^{(\eps)}(t,\bfx)
\,=\,\bcU^{(\eps)}\left(\eps\,t,\eps\,\bfx;
\frac{\bPsi^{(\eps)}(\eps\,t,\eps\,\bfx)}{\eps}\right)
\ee 
with, for any $(T,\bY)$, $\bfzeta\mapsto\cU^{(\eps)}(T,\bY;\bfzeta)$ $(\beD_1,\beD_2)$-periodic and, as $\eps\to0$,
\begin{align*}
\bcU^{(\eps)}(T,\bY;\bfzeta)&=
\bcU_{(0)}(T,\bY;\bfzeta)+\eps\,\bcU_{(1)}(T,\bY;\bfzeta)+\cO(\eps^2)\,,\\
\bPsi^{(\eps)}(T,\bY)&=
\bPsi_{(0)}(T,\bY)+\eps\,\bPsi_{(1)}(T,\bY)+\cO(\eps^2)\,.
\end{align*}
Requiring \eqref{e:ansatz} to solve \eqref{e:rd-intro} up to a remainder of size $\cO(\eps)$ is equivalent to $\bfzeta\mapsto\bcU_{(0)}(T,\bY;\bfzeta)$ being a scaled periodic traveling wave of profile. Explicitly,
\[
0\,=\,\transp{(\bcK_{(0)}\nabla_{\bfzeta})}\left((\bcK_{(0)}\nabla_{\bfzeta})\bcU_{(0)}\right)+\transp{(\bcK_{(0)}\nabla_{\bfzeta})}\left(\bG(\bcU_{(0)})\right)-(\bOm_{(0)}\cdot \nabla_{\bfzeta}) \bcU_{(0)} + \bff(\bcU_{(0)})\,,
\]
with local parameters (depending on slow variables $(T,\bY)$) related to phases by
\begin{align*}
\d_T\bPsi_{(0)}&=\bOm_{(0)}\,,&
\nabla_\bY \bPsi_{(0)}&=\bcK_{(0)}\,.
\end{align*}
Choosing a wave parametrization as in Proposition~\ref{p:structure}, this is solved by imposing the slow-modulation form
\begin{align*}
\bcU_{(0)}(T,\bY;\bfzeta)
&=\bU^{\bcK_{(0)}(T,\bY)}(\bfzeta)\,,&
\bcK_{(0)}(T,\bY)
&=\nabla_\bY \bPsi_{(0)}(T,\bY)\,.
\end{align*}
jointly with the slow evolution equation
\be\label{e:W-phase}
\d_T\bPsi_{(0)}\,=\,\bOm(\nabla_\bY \bPsi_{(0)})\,.
\ee

System~\eqref{e:W-phase} fails to capture dissipative effects, because they are high-order with respect to slow expansions. As far as large-time analysis is concerned, one could just correct System~\eqref{e:W-phase} with an artificial semilinear second-order term enforcing a good description of slow/low-Floquet expansions up to second-order. For an example of the latter we refer the reader to \cite{R_linKdV} (with dispersion instead of diffusion). Yet, instead, we follow \cite{Noble-Rodrigues} and show how going on with the formal identification provides a relevant higher-order correction.

Requiring \eqref{e:ansatz} to solve \eqref{e:rd-intro} up to a remainder of size $\cO(\eps^2)$ provides, besides the foregoing equalities, the extra constraint 
\begin{align*}
0&\,=\,\transp{(\bcK_{(0)}\nabla_{\bfzeta})}\left((\bcK_{(0)}\nabla_{\bfzeta})\bcU_{(1)}\right)+\transp{(\bcK_{(0)}\nabla_{\bfzeta})}\left(\dD \! \bG(\bcU_{(0)})(\bcU_{(1)}) \right)-(\bOm_{(0)}\cdot \nabla_{\bfzeta}) \bcU_{(1)} + \dD \! \bff(\bcU_{(0)})(\bcU_{(1)})\\
&+\transp{(\bcK_{(1)}\nabla_{\bfzeta})}\left((\bcK_{(0)}\nabla_{\bfzeta})\bcU_{(0)}\right)
+\transp{(\bcK_{(0)}\nabla_{\bfzeta})}\left((\bcK_{(1)}\nabla_{\bfzeta})\bcU_{(0)}\right)
+\transp{(\bcK_{(1)}\nabla_{\bfzeta})}\left(\bG(\bcU_{(0)})\right)-(\bOm_{(1)}\cdot \nabla_{\bfzeta}) \bcU_{(0)}\\
&+\transp{\nabla_{\bY}}\left((\bcK_{(0)}\nabla_{\bfzeta})\bcU_{(0)}\right)
+\transp{(\bcK_{(0)}\nabla_{\bfzeta})}\left(\nabla_{\bY}\bcU_{(0)}\right)
+\transp{\nabla_\bY}\left(\bG(\bcU_{(0)})\right)-\d_T\bcU_{(0)}\,,
\end{align*}
with $\bcK_{(1)}:=\nabla_\bY \bPsi_{(1)}$, $\bOm_{(1)}:=\d_t\bPsi_{(1)}$. Denoting $L^{\bcK_0}_{\b0}$ the linearized operator $L$ in the variable $\bfzeta$ and corresponding to the profile $\bcU_{(0)}=\bU^{\bcK_0}$ and using relations from Subsection~\ref{ss:algebraic}, the constraint is equivalently written as 
\begin{align*}
0&\,=\,L^{\bcK_0}_{\b0}[\bcU_{(1)}-\dD_{\bK}\bU(\bcK_{(0)})(\bcK_{(1)})]
-\left(\left(\bOm_{(1)}-\dD_{\bK}\Omega(\bcK_{(0)})(\bcK_{(1)})\right)\cdot \nabla_{\bfzeta}\right) \bcU_{(0)}\\
&+\transp{\nabla_{\bY}}\left((\bcK_{(0)}\nabla_{\bfzeta})\bcU_{(0)}\right)
+\transp{(\bcK_{(0)}\nabla_{\bfzeta})}\left(\nabla_{\bY}\bcU_{(0)}\right)
+\transp{\nabla_\bY}\left(\bG(\bcU_{(0)})\right)-\d_T\bcU_{(0)}\,,
\end{align*}
Introducing\footnote{We warn the reader of the notational inconsistency $\tbq_j^{\ubK}=\tbq_j^{\b0}$.} $\tbq_1^{\bK_0}$, $\tbq_2^{\bK_0}$ a basis of the kernel of the adjoint of $L^{\bK_0}_{\b0}$, in duality with $\d_{\zeta_1}\bcU_{(0)}=\d_1\bU^{\bcK_0}$, $\d_{\zeta_2}\bcU_{(0)}=\d_2\bU^{\bcK_0}$, we deduce as a necessary constraint
\begin{align*}
\bOm_{(1)}-\dD_{\bK}\Omega(\bcK_{(0)})&(\bcK_{(1)})\\
=\bp \langle\tbq_1^{\bcK_0};\cdot\rangle_{L^2_{per}}\\
\langle\tbq_2^{\bcK_0};\cdot\rangle_{L^2_{per}}\ep
&\Big[\transp{\nabla_{\bY}}\left((\bcK_{(0)}\nabla_{\bfzeta})\bcU_{(0)}\right)
+\transp{(\bcK_{(0)}\nabla_{\bfzeta})}\left(\nabla_{\bY}\bcU_{(0)}\right)\\
&\qquad+\transp{\nabla_\bY}\left(\bG(\bcU_{(0)})\right)-\dD_{\bK}\bU(\bcK_{(0)})(\nabla_\bY\bOm(\bcK_{(0)}))\Big]\,,
\end{align*}
that we denote in abstract form
\be\label{e:W-phase2}
\d_T\bPsi_{(1)}\,=\,\dD_{\bK}\bOm(\nabla_\bY \bPsi_{(0)})(\nabla_\bY \bPsi_{(1)})
+\bLambda^{\nabla_\bY \bPsi_{(0)}}[\nabla_\bY](\nabla_\bY \bPsi_{(0)})\,,
\ee
with
\be\label{def:blambda-1}
\bLambda^{\bK}[\nabla_\bY](\tbK)
\,=\,\sum_{j=1}^2\sum_{\ell=1}^2\sum_{m=1}^2\sum_{p=1}^2 \bLambda_{p,m}^{j,\ell}(\bK)\,\d_j(\tbK_{p,m})\,\beD_\ell\,
\ee
where 
\begin{align}\label{def:blambda-2}
\bLambda_{p,m}^{j,\ell}(\bK)
&=\,\delta_{\ell,m}\,\delta_{j,p}
\,+\,\langle\tbq_\ell^{\bK};
\transp{(L^{(1)}_{\bK}[\d_{\bK_{p,m}}\bU(\bK)])}\beD_j
\,+\,
\sum_{r=1}^2\dD_{\bK_r}\bU(\bK)(\transp{\bK}\d_{\bK_{p,m}}\bfc_r(\bK))\rangle_{L^2_{per}}\,,
\end{align}
$L^{(1)}_{\bK}$ being associated with $L^{\bK}$ through \eqref{def:Bloch-expand}, as $L^{(1)}$ with $L$.  

Introducing $\bPsi=\bPsi_{(0)}+\eps\,\bPsi_{(1)}$ and grouping together \eqref{e:W-phase}-\eqref{e:W-phase2} yield, up to $\cO(\eps^2)$ terms that we discard,
\[
\d_T\bPsi\,=\,\bOm(\nabla_\bY \bPsi)
+\bLambda^{\nabla_\bY \bPsi}[\eps\,\nabla_\bY](\nabla_\bY \bPsi)\,.
\]
Going back to original $(t,\bfx)$-variables, the upshot of the formal analysis is that we may expect
\begin{align*}
\bcW(t,\bfx)&\approx \bU^{\bcK(t,\bfx)}\left(\bPsi(t,\bfx)\right)\,,&
\bcK(t,\bfx)&=\nabla\bPsi(t,\bfx)\,,
\end{align*}
with $\bPsi$ satisfying 
\be\label{e:W-phase-formal}
\d_t\bPsi\,=\,\bOm(\nabla_\bfx\bPsi)
+\bLambda^{\nabla_\bfx\bPsi}[\nabla_\bfx](\nabla_\bfx\bPsi)\,.
\ee
Alternatively, one may observe that the slow evolution obeys
\be\label{e:W-wn-formal}
\d_t\bcK\,=\,\nabla_\bfx\left(\bOm(\bcK)\right)
+\nabla_\bfx\left(\bLambda^{\bcK}[\nabla_\bfx](\bcK)\right)\,,
\ee
with $\bcK$ curl-free.

Whereas the formal arguments expounded so far do contain some form of large-time considerations since implicitly here the time variable $t$ lives in an interval of length $\cO(1/\eps)$ with $\eps\to0$, it is not specialized to the situation at stake in the rest of the paper where $t\to\infty$ and $\bcK$ is sufficiently close to $\ubK$. In the present paper, we consider cases where nonlinear terms are at worst critical from the point of view of time decay so that it is only necessary to retain nonlinear terms with the worst decay rates if one aims at a leading-order description. Moreover, at the level of wavevectors, the decay is inherently the one of conservative hyperbolic-parabolic systems near constant states so that every extra spatial derivative is expected to bring an extra $t^{-\tfrac12}$ decay. With this in mind, for our purposes we expect that it should be sufficient to retain from \eqref{e:W-phase-formal} either
\[
\d_t\bPsi\,=\,\ubOm+\dD_{\bK}\bOm(\ubK)(\nabla_\bfx\bPsi-\ubK)
+\frac12\dD_{\bK}^2\bOm(\ubK)(\nabla_\bfx\bPsi-\ubK,\nabla_\bfx\bPsi-\ubK)
+\bLambda^{\ubK}[\nabla_\bfx](\nabla_\bfx\bPsi)\,,
\]
or, alternatively, if one prefers to keep a compact form with the same level of approximation
\be\label{e:W-phase-semilinear}
\d_t\bPsi\,=\,\bOm(\nabla_\bfx\bPsi)
+\bLambda^{\ubK}[\nabla_\bfx](\nabla_\bfx\bPsi)\,.
\ee
We stress that we regard the rigorous justification of System~\eqref{e:W-phase-semilinear} from \eqref{e:W-phase-formal} as a routine task in the sense that the proof of the above formal claims that quadratic second-order terms (or even cubic first-order terms) could be discarded would follow from a direct inspection of a Duhamel formula if we already knew that System~\eqref{e:W-phase-formal} was well-posed in some dissipative sense, including some form of high-frequency damping estimates similar to (but possibly weaker than) those of Lemmas~\ref{l:energy-estimate} and~\ref{l:energy-estimate-Lp}. 

In cases subcritical from the point of view of time decay, almost by definition, we expect to be allowed to go even further and simply retain the linear version
\be\label{e:W-phase-linear}
\d_t\bPsi\,=\,\ubOm+\dD_{\bK}\bOm(\ubK)(\nabla_\bfx\bPsi-\ubK)
+\bLambda^{\ubK}[\nabla_\bfx](\nabla_\bfx\bPsi)\,.
\ee

\subsection{Asymptotic equivalence of hyperbolic-parabolic systems}\label{ss:artificial}

The issue we want to address now is that the derivation of \eqref{e:W-phase-formal}, in a slow expansion regime, brings relevant information only of low-frequency type. In particular, it could well be that despite the fact that \eqref{e:W-phase-formal} contains a correct large-time low-frequency description of \eqref{e:rd-intro}, the system is ill-posed because of high-frequency instabilities having nothing to do with the original system. The issue is ubiquitous in the theory and we refer the reader to \cite{Noble-Rodrigues,R_linKdV} for closely related discussions.

To be more precise, we observe that in the low-frequency regime the leading-order part is the first-order hyperbolic part and the second-order part only brings corrections. As a consequence, in a direction where strict hyperbolicity is met for the first-order part, Assumption \cond{1} is reduced to a sign condition on two coefficients of the four-dimensional second-order operator whereas the high-frequency properties involve the missing coefficients. To give a concrete example, note that 
\begin{align*}
\d_t u_1\,+\,\d_{x_1}u_1&\,=\,(\d_{x_1}^2+\d_{x_2}^2)\,u_1+2\,\d_{x_1}^2\,u_2\\
\d_t u_2\,-\,\d_{x_1}u_2&\,=\,(\d_{x_1}^2+\d_{x_2}^2)\,u_2+2\,\d_{x_1}^2\,u_1
\end{align*}
exhibits both a diffusive low-frequency behavior, compatible with \cond{1}, and a violent high-frequency ill-posedness. Consistently, the only situation where we are able to deduce good high-frequency properties is when the first-order part is strictly hyperbolic in no direction, that is, in Subcase~\case{b0} when the first-order part is scalar.

For this reason, we show here how to replace \eqref{e:W-phase-semilinear} with a well-posed system expected to share, at leading-order, the same large-time dynamics. The discussion is parallel to the one in \cite[Appendix~B.2]{JNRZ-conservation} and extends in various ways, including the class of systems considered and the sharpness of estimates proved, the analysis about \emph{artificial viscosity systems} in \cite{Hoff_Zumbrun-NS_compressible_pres_de_zero,Rodrigues-compressible} (discussed further in \cite{Rodrigues-these} and \cite[Appendix~A]{R}). Even if System~\eqref{e:W-phase-formal} were known to be well-posed in some dissipative sense, there would be a gain in simplicity --- but a loss in explicitness --- in replacing \eqref{e:W-phase-semilinear} in the way expounded here since the systems introduced below are to be semilinear, genuinely parabolic, with first-order and second-order parts commuting at the linear level. This is precisely the commutation property that enables one to extend the good low-frequency properties to the whole dynamics. 

The issue is linear in essence so that our task is to identify a $\ubLambda_0[\nabla]$ such that one could replace \eqref{e:W-phase-linear} with
\[
\d_t\bPsi\,=\,\ubOm+\dD_{\bK}\bOm(\ubK)(\nabla_\bfx\bPsi-\ubK)
+\ubLambda_0[\nabla](\nabla\bPsi)\,.
\]
To begin with, we take a step back from the foregoing discussion and continue the study of Appendix~\ref{ss:diff} --- with notational conventions introduced there  --- so as to prove that the linearized evolution contained in \eqref{e:W-phase-linear} does reproduce correctly the averaged low-Floquet evolution. To ease comparisons we write $\bPsi$ in a co-moving frame and in a perturbative form
\be\label{bfpsi}
\bPsi(t,\bfx)\,=\,\transp{\ubK}\,(\bfx-t\,\ubc)
+\bfpsi\left(t,\transp{\ubK}\,(\bfx-t\,\ubc)\right)\,.
\ee
This turns \eqref{e:W-phase-linear} into
\be\label{e:W-phase-lin}
\d_t\bfpsi\,=\,
-\transp{\ubK}\dD_{\bK}\bfc(\ubK)(\ubK\nabla_\bfx\bfpsi)
+\bLambda^{\ubK}[\ubK\nabla_\bfx](\ubK\nabla_\bfx\,\bfpsi)
\ee
Note that Lemmas~\ref{l:spec-wn} and~\ref{l:spec-2nd} contain that the latter system is equivalently written as $\d_t\bfpsi\,=\,\bDW(\nabla_\bfx)\bfpsi$. Accordingly we introduce evolution operators $(\SigW(t))_{t\geq0}$ and $(\SigLF(t))_{t\geq0}$ defined as 
\begin{align}\label{def:SigW}
\widehat{(\SigW(t)\bfphi)}(\bfxi)&\,=\,\chi(\bfxi)\,\beD^{t\,\bDW(\iD\bfxi)}\,\widehat{\bfphi}(\bfxi)\,,&
\widehat{(\SigLF(t)\bfphi)}(\bfxi)&\,=\,\chi(\bfxi)\,\beD^{t\,\bD_\bfxi}\,\widehat{\bfphi}(\bfxi)\,.
\end{align}
Note that the low-frequency cut-off is needed in the definition of $\SigW$ because of the ill-posedness issues already mentioned and in the definition of $\SigLF$ since $\bD_\bfxi$ is not even defined when $\bfxi$ is not small.

We now come back to the question of identifying $\ubLambda_0$ such that one may replace \eqref{e:W-phase-lin} with
\be\label{e:W-phase-lin0}
\d_t\bfpsi\,=\,
-\transp{\ubK}\dD_{\bK}\bfc(\ubK)(\ubK\nabla_\bfx\bfpsi)
+\ubLambda_0[\ubK\nabla_\bfx](\ubK\nabla_\bfx\,\bfpsi)
\ee
and derive parabolic behavior without altering large-time low-frequency dynamics. Let us anticipate the choices of $\ubLambda_0$ detailed below and define for later comparisons
\be\label{def:SigW0}
\bDO(\bfeta)\bfphi:=\bA(\bfeta)\bfphi+\ubLambda_0[\ubK\bfeta](\ubK\bfeta\,\transp{\bfphi})
\ee
and
\bas
\widehat{(\SigO(t)[\bfphi])}(\bfxi)&\,=\,\beD^{t\,\bDO(\iD\bfxi)}\,\widehat{\bfphi}(\bfxi)\,,\\
\widehat{(\SigOLF(t)[\bfphi])}(\bfxi)&\,=\,\chi(\bfxi)\,\beD^{t\,\bDO(\iD\bfxi)}\,\widehat{\bfphi}(\bfxi)\,,\\
\widehat{(\SigOHF(t)[\bfphi])}(\bfxi)&\,=\,(1-\chi(\bfxi))\,\beD^{t\,\bDO(\iD\bfxi)}\,\widehat{\bfphi}(\bfxi)\,.
\eas

In Subcase~\case{b0}, we may simply set $\ubLambda_0=\bLambda^{\ubK}$ so that $\SigOLF=\SigW$.

\bpr\label{pr:ubL0-scalar0}
Assume \cond1-\cond2 and Subcase~\case{b0}, and define $\ubLambda_0$ by $\ubLambda_0=\bLambda^{\ubK}$. 
\begin{enumerate}
\item For any $\alpha\in\N^2$ and any $\ell\in\N$, there exists $C_{\alpha,\ell}$ such that for any  $2 \leq p \leq +\infty$, $1\leq q \leq p$, and any $t\geq0$
\begin{align*}
\|\,\d_\bfx^\alpha\,\d_t^\ell\,(\SigLF-\SigOLF)(t)[\bfg]\|_{L^p}
&\leq \frac{C_{\alpha,\ell}}{(1+t)^{\frac{|\alpha|+\ell}{2}+\frac{1}{q}-\frac{1}{p}+\frac12}}\,\|\bfg\|_{L^{q}}\,.
\end{align*}
\item For any $\alpha\in\N^2$ and any $\ell\in\N$, there exists $C_{\alpha,\ell}$ such that for any $2 \leq p \leq +\infty$, and any $t\geq0$
\begin{align*}
\|\,\d_\bfx^\alpha\,\d_t^\ell\,(\SigLF-\SigOLF)(t)[\bfphi]\|_{L^p}
&\leq \frac{C_{\alpha,\ell}}{(1+t)^{\frac{|\alpha|+\ell}{2}+\frac12-\frac{1}{p}}}\,\|\Delta\bfphi\|_{L^1}\,.
\end{align*}
\item For any $\alpha\in\N^2$ and any $\ell\in\N$, there exists $C_{\alpha,\ell}$ such that for any  $2 \leq p \leq +\infty$, $1\leq q \leq p$, and any $t\geq0$
\begin{align*}
\|\,\d_\bfx^\alpha\,\d_t^\ell\,\SigOLF(t)[\bfg]\|_{L^p}
&\leq \frac{C_{\alpha,\ell}}{(1+t)^{\frac{|\alpha|+\ell}{2}+\frac{1}{q}-\frac{1}{p}}}\,\|\bfg\|_{L^{q}}\,.
\end{align*}
\item For any $\alpha\in\N^2$, any $\ell\in\N$ and any $2 \leq p \leq +\infty$ such that $|\alpha|+\ell-\frac2p>0$, there exists $C_{p,\alpha,\ell}$ such that for any $t\geq0$
\begin{align*}
\|\,\d_\bfx^\alpha\,\d_t^\ell\,\SigOLF(t)[\bfphi]\|_{L^p}
&\leq \frac{C_{p,\alpha,\ell}}{(1+t)^{\frac{|\alpha|+\ell}{2}-\frac{1}{p}}}\,\|\Delta\bfphi\|_{L^1}\,.
\end{align*}
\item There exists $\theta>0$ such that for any $\alpha\in\N^2$ and any $\ell\in\N$, there exists $C_{\alpha,\ell}$ such that for any  $2 \leq p \leq +\infty$, $1\leq q \leq p$, and any $t\geq0$
\begin{align*}
\|\,\d_\bfx^\alpha\,\d_t^\ell\,\SigOHF(t)[\bfg]\|_{L^p}
&\leq \frac{C_{\alpha,\ell}\,\eD^{-\theta\,t}}{(\min(\{1,t\}))^{\frac{|\alpha|+2\,\ell}{2}+\frac{1}{q}-\frac{1}{p}}}\,\|\bfg\|_{L^{q}}\,.
\end{align*}
\item There exists $\theta>0$ such that for any $\alpha\in\N^2$ and any $\ell\in\N$, there exists $C_{\alpha,\ell}$ such that for any  $2 \leq p \leq +\infty$, and any $t\geq0$
\begin{align*}
\|\,\d_\bfx^\alpha\,\d_t^\ell\,\SigOHF(t)[\bfphi]\|_{L^p}
&\leq \frac{C_{\alpha,\ell}\,\eD^{-\theta\,t}}{(\min(\{1,t\}))^{\left(\frac{|\alpha|+2\,\ell}{2}-\frac{1}{p}\right)_+}}\,\|\Delta\bfphi\|_{L^1}\,.
\end{align*}
\end{enumerate}
\epr

\begin{proof}
Let us recall that
\be\label{e:LF-W}
\eD^{t\,\bD_{\bfxi}}-\eD^{t\,\bDW(\iD\bfxi)}
=\int_{0}^{t}\eD^{(t-s)\,\bD_\bfxi} (\bD_{\bfxi}-\bDW(\iD\bfxi))\,\eD^{s\,\bDW(\iD\bfxi)}\,\dD s\,
\ee
where $\lnor \bD_{\bfxi}-\bDW(\iD\bfxi) \rnor \lesssim \| \bfxi \|^{3}$. This is sufficient to deduce from Hausdorff-Young inequalities the $L^q\to L^p$ bounds when $1\leq q \leq 2 \leq p \leq +\infty$. However formula~\eqref{e:LF-W} is also well-adapted to the arguments of Subsection~\ref{ss:scalar0}. When considering the case $p=q=+\infty$, we get an integral with a form similar to \eqref{int:scalarcase} and one can use Lemma \ref{l:Lp}-(2). Altogether this yields the first two sets of inequalities.

The last estimates follow from Hausdorff-Young inequalities and arguments of Subsection \ref{ss:linear_estim_mod} through Green functions representations.
\end{proof}

When Case~\case{b} holds but Subcase~\case{b0} fails, we define $\ubLambda_0$ through
\be\label{def:ubL0-scalar}
\ubLambda_0[\ubK\bfeta](\ubK\bfeta\,\transp{\bfphi})
\,:=\,
\bP^{-1}\diag((\bP\bB(\bfeta)\bP^{-1})_{1,1},(\bP\bB(\bfeta)\bP^{-1})_{2,2})\bP\,\bfphi
\ee
where $\bP$ diagonalizes $\bA_1$ and $\bA_2$. We stress that this definition does not depend on $\bP$ (since two convenient $\bP$s only differ by a multiplication from the left by a diagonal matrix) and that it reduces System~\eqref{e:W-phase-lin0} to two uncoupled scalar transport-diffusion equations (in a suitable basis). 

\bpr\label{pr:ubL0-scalar}
Assume \cond1-\cond2 and Case~\case{b} but with Subcase~\case{b0} failing. Define $\ubLambda_0$ by \eqref{def:ubL0-scalar}. 
\begin{enumerate}
\item For any $\alpha\in\N^2$, any $\ell\in\N$, and any  $2 \leq p \leq +\infty$, $1\leq q \leq p$, there exists $C_{\alpha,\ell,p,q}$ such that for any $t\geq0$
\begin{align*}
\|\,\d_\bfx^\alpha\,\d_t^\ell\,(\SigLF-\SigOLF)(t)[\bfg]\|_{L^p}
&\leq 
\begin{cases}
\frac{C_{\alpha,\ell,p,q}}{(1+t)^{\frac{|\alpha|+\ell}{2}+\frac{1}{q}-\frac{1}{p}+\frac12\left(\frac1q-\frac1p\right)}}\,\|\bfg\|_{L^{q}}\,,&(p,q)\neq (\infty,1)\,,\\
\frac{C_{\alpha,\ell,p,q}\,\ln\left(2+t\right)\,}{(1+t)^{\frac{|\alpha|+\ell}{2}+\frac{3}{2}}}\,\| \bfg\|_{L^{q}}\,,&(p,q)=(\infty,1)\,.
\end{cases}
\end{align*}
\item For any $\alpha\in\N^2$, any $\ell\in\N$, and any  $2 \leq p \leq +\infty$, there exists $C_{\alpha,\ell,p}$ such that for any $2 \leq p \leq +\infty$, and any $t\geq0$
\begin{align*}
\|\,\d_\bfx^\alpha\,\d_t^\ell\,(\SigLF-\SigOLF)(t)[\bfphi]\|_{L^p}
&\leq 
\begin{cases}
\frac{C_{\alpha,\ell,p}}{(1+t)^{\frac{|\alpha|+\ell}{2}-\frac{1}{p}+\frac12\left(1-\frac1p\right)}}\,\|\Delta\bfphi\|_{L^1}\,,&p\neq\infty\,,\\
\frac{C_{\alpha,\ell,p}\,\ln\left(2+t\right)}{(1+t)^{\frac{|\alpha|+\ell}{2}+\frac{1}{2}}} \,\|\Delta\bfphi\|_{L^1}\,,&p=\infty\,.
\end{cases}
\end{align*}
\item For any $\alpha\in\N^2$ and any $\ell\in\N$, there exists $C_{\alpha,\ell}$ such that for any  $2 \leq p \leq +\infty$, $1\leq q \leq p$, and any $t\geq0$
\begin{align*}
\|\,\d_\bfx^\alpha\,\d_t^\ell\,\SigOLF(t)[\bfg]\|_{L^p}
&\leq \frac{C_{\alpha,\ell}}{(1+t)^{\frac{|\alpha|+\ell}{2}+\frac{1}{q}-\frac{1}{p}}}\,\|\bfg\|_{L^{q}}\,.
\end{align*}
\item For any $\alpha\in\N^2$, any $\ell\in\N$ and any $2 \leq p \leq +\infty$ such that $|\alpha|+\ell-\frac2p>0$, there exists $C_{p,\alpha,\ell}$ such that for any $t\geq0$
\begin{align*}
\|\,\d_\bfx^\alpha\,\d_t^\ell\,\SigOLF(t)[\bfphi]\|_{L^p}
&\leq \frac{C_{p,\alpha,\ell}}{(1+t)^{\frac{|\alpha|+\ell}{2}-\frac{1}{p}}}\,\|\Delta\bfphi\|_{L^1}\,.
\end{align*}
\item There exists $\theta>0$ such that for any $\alpha\in\N^2$ and any $\ell\in\N$, there exists $C_{\alpha,\ell}$ such that for any  $2 \leq p \leq +\infty$, $1\leq q \leq p$, and any $t\geq0$
\begin{align*}
\|\,\d_\bfx^\alpha\,\d_t^\ell\,\SigOHF(t)[\bfg]\|_{L^p}
&\leq \frac{C_{\alpha,\ell}\,\eD^{-\theta\,t}}{(\min(\{1,t\}))^{\frac{|\alpha|+2\,\ell}{2}+\frac{1}{q}-\frac{1}{p}}}\,\|\bfg\|_{L^{q}}\,.
\end{align*}
\item There exists $\theta>0$ such that for any $\alpha\in\N^2$ and any $\ell\in\N$, there exists $C_{\alpha,\ell}$ such that for any  $2 \leq p \leq +\infty$, and any $t\geq0$
\begin{align*}
\|\,\d_\bfx^\alpha\,\d_t^\ell\,\SigOHF(t)[\bfphi]\|_{L^p}
&\leq \frac{C_{\alpha,\ell}\,\eD^{-\theta\,t}}{(\min(\{1,t\}))^{\left(\frac{|\alpha|+2\,\ell}{2}-\frac{1}{p}\right)_+}}\,\|\Delta\bfphi\|_{L^1}\,.
\end{align*}
\end{enumerate}
\epr

\begin{proof}
A small variation on the proof of Proposition~\ref{pr:ubL0-scalar0} provides a version of the proposition where in the first two estimates $\SigOLF$ is replaced with $\SigW$ (and the time decay is actually stronger). Thus we only need to explain how to bound $\SigW-\SigOLF$.
We denote
\[
\bGamma(t,z) = \int_{[-\pi,\pi]^{2}} e^{\iD \bfz \cdot \bfxi} \left( \eD^{t\,\bDW_{\bfxi}}-\eD^{t\,\bDO(\iD\bfxi)} \right) \chi(\bfxi) \dD \bfxi := \int_{[-\pi,\pi]^{2}} e^{\iD \bfz \cdot \bfxi} \bfm(t,\bfxi) \dD \bfxi
\]
and we have to bound $\| \bGamma(t,\cdot) \|_{L^{r}(\R^{2})}$ for $1 \leq r \leq +\infty$. As encoded in \eqref{e:Sig-aux}-\eqref{e:Sig-aux2}, the analysis of Subsection~\ref{ss:scalar} (applied with $\eD^{t\,\bDW(\iD\bfxi)}$ replacing $\eD^{t\,\bD_{\bfxi}}$) is actually already written in terms of comparisons with the evolution of \eqref{e:W-phase-lin0} so that $\bfm$ satisfies estimates similar to \eqref{e:Sig-aux}. The arguments expounded there complete the proof when $1 \leq r \leq 2$ and also $2<r<+\infty$, since, using Hausdorff-Young inequalities, we only have to bound $\|\bfm(t,\cdot)\|_{L^{\frac{r}{r-1}}}$. There is one detail worth mentioning, in \eqref{e:Sig-aux} we have absorbed a factor $\|\bfxi\|^2\,t$ in the exponential but it is useful to keep it apparent when bounding with $\|\Delta\bfphi\|_{L^1}$. Furthermore when $r=+\infty$, we can bound $\| \bfm(t,\cdot) \|_{L^{1}}$ by splitting the integration domain into three areas corresponding to $|\xi_\mypar|\leq |\xi_\perp|^2$,  $|\xi_\perp|^2\leq|\xi_\mypar| \leq (1+t)|\xi_\perp|^2$ and $(1+t)|\xi_\perp|^2\leq |\xi_\mypar|$, the factor $\eD^{-\theta\,t\,|\xi_\mypar|^2}$ being bounded by $1$ in the first two areas and, when $t\geq 1$, by a multiple of $(t\,|\xi_\mypar|^2)^{-\eta}$ for some $\eta>0$ in the last area. 
\end{proof} 

We now turn to Case~\case{a}. With polar coordinates conventions of Subsection~\ref{ss:dispersive}, including identification of $\bfxi$ with $(r,\omega)$, let us recall that there exist complex-valued maps $\lambda_1$, $\lambda_2$ and complementary projector-valued maps $\pi_1$, $\pi_2$, all smooth in polar coordinates, such that
\[
\eD^{t\,\bD_\bfxi}\,=\,\sum_{j=1}^2\eD^{t\,\lambda_j(\bfxi)}\,\pi_j(\bfxi)
\]
with, for $j=1,2$, $\lambda_j$ continuous at $\bfxi=\b0$ with value $0$, and for some $\theta>0$,
\[
\Re(\lambda_j(\bfxi))\,\leq -\theta\,\|\bfxi\|^2\,,
\]
and $\d_r\lambda_j+\d^2_\omega\d_r\lambda_j$ nowhere vanishing. We define 
\begin{align*}
\alpha_j(\bfxi)&:=\frac{1}{\iD}\,r\,(\d_r\lambda_j)_{|r=0}(\bfxi)\,,&
\beta_j(\bfxi)&:=\frac12\,r^2\,(\d_r^2\lambda_j)_{|r=0}(\bfxi)\,,&\\
\lambO_j(\bfxi)&:=\iD\,\alpha_j(\bfxi)+\beta_j(\bfxi)\,,&
\piO_j(\bfxi)&:=(\pi_j)_{|r=0}(\bfxi)\,,&
\end{align*}
where ${}_{|r=0}$ means that instead of evaluating at $(r,\omega)$ we evaluate at $(0,\omega)$. Note that $\alpha_j$, $\beta_j$ and $\piO_j$ are defined over $\R^2$, with respective homogeneity $1$, $2$ and $0$, that $\alpha_j$ is real valued with $\d_r\alpha_j+\d^2_\omega\d_r\alpha_j$ bounded away from zero, that $\Re(\d_r^2\beta_j)$ is bounded from above away from zero and that
\begin{align*}
\lambda_j(\bfxi)&\stackrel{\bfxi\to\b0}{=}\lambO_j(\bfxi)+\cO(\|\bfxi\|^3)\,,&
\pi_j(\bfxi)&\stackrel{\bfxi\to\b0}{=}\piO_j(\bfxi)+\cO(\|\bfxi\|)\,,&
\bA(\iD\bfxi)&=\sum_{j=1}^2\iD\,\alpha_j(\bfxi)\,\piO_j(\bfxi)\,.%\,,\quad\|\bfxi\|\leq \xi_0\,.
\end{align*}
Now we define $\ubLambda_0$ through
\be\label{def:ubL0-dispersive}
\ubLambda_0[\ubK\bfeta](\ubK\bfeta\,\transp{\bfphi})
\,:=\,
\sum_{j=1}^2 \,\beta_j(\bfeta)\,\piO_j(\bfeta)\,\bfphi.
\ee
Replacing $\bD_{\bfxi}$ with $\bDW(\iD\bfxi)$ in the foregoing definition of $\ubLambda_0$ would not change the value of $\ubLambda_0$ since it only involves Taylor expansions captured by $\bDW(\iD\bfxi)$. In particular, $\ubLambda_0$ may equally be obtained from \eqref{e:W-phase-linear}. 

\br\label{rk:ubL0}
We stress that definition \eqref{def:ubL0-dispersive} is conceptually similar to \eqref{def:ubL0-scalar} in that both define the new second-order operator as being, in a frame diagonalizing the first-order expansion, the diagonal part of the second-order expansion. A strong difference is that since here the diagonalization is given by a Fourier multiplier instead of a constant matrix, the resulting operator is an homogeneous second-order multiplier instead of a differential operator. Incidentally we point out that similar analyses in the literature are for the moment restricted to either one-dimensional situations \cite{Liu-Zeng} or isotropic situations \cite{Hoff_Zumbrun-NS_compressible_pres_de_zero,Rodrigues-compressible}, thus to cases where extra cancellations do give back a differential operator.
\er

\bpr\label{pr:ubL0-dispersive}
Assume \cond1-\cond2 and Case~\case{a}. Define $\ubLambda_0$ by \eqref{def:ubL0-dispersive}. 
\begin{enumerate}
\item For any $\alpha\in\N^2$ and any $\ell\in\N$, there exists $C_{\alpha,\ell}$ such that for any  $2 \leq p \leq +\infty$, $1\leq q \leq 2$, and any $t\geq0$
\begin{align*}
\|\,\d_\bfx^\alpha\,\d_t^\ell\,(\SigLF-\SigOLF)(t)[\bfg]\|_{L^p}
&\leq \frac{C_{\alpha,\ell}}{(1+t)^{\frac{|\alpha|+\ell}{2}+\frac{1}{q}-\frac{1}{p}+ \frac12 + \frac12 \min \left( \left\{ \frac12 -\frac{1}{p} , \frac{1}{q} - \frac12 \right\}\right)}}\,\|\bfg\|_{L^{q}}\,.
\end{align*}
\item For any $\alpha\in\N^2$ and any $\ell\in\N$, there exists $C_{\alpha,\ell}$ such that for  any $2 \leq p \leq +\infty$ and any $t\geq0$
\begin{align*}
\|\,\d_\bfx^\alpha\,\d_t^\ell\,(\SigLF-\SigOLF)(t)[\bfphi]\|_{L^p}
&\leq \frac{C_{\alpha,\ell}}{(1+t)^{\frac{|\alpha|+\ell}{2}+\frac34-\frac32\frac{1}{p}}}\,\|\Delta\bfphi\|_{L^1}\,.
\end{align*}
\item For any $\alpha\in\N^2$ and any $\ell\in\N$, there exists $C_{\alpha,\ell}$ such that for any  $2 \leq p \leq +\infty$, $1\leq q \leq 2$, and any $t\geq0$
\begin{align*}
\|\,\d_\bfx^\alpha\,\d_t^\ell\,\SigOLF(t)[\bfg]\|_{L^p}
&\leq 
\frac{C_{\alpha,\ell}}{(1+t)^{\frac{|\alpha|+\ell}{2}+\frac{1}{q}-\frac{1}{p} +\,\frac12 \min \left( \left\{ \frac12 -\frac{1}{p} , \frac{1}{q} - \frac12 \right\}\right)}}
\,\|\bfg\|_{L^{q}}\,.
\end{align*}
\item For any $\alpha\in\N^2$, any $\ell\in\N$ and any $2 \leq p \leq +\infty$ such that $|\alpha|+\ell-\tfrac{2}{p}>0$, there exists $C_{p,\alpha,\ell}$ such that for any $t\geq0$
\begin{align*}
\|\,\d_\bfx^\alpha\,\d_t^\ell\,\SigOLF(t)[\bfphi]\|_{L^p}
&\leq \frac{C_{p,\alpha,\ell}}{(1+t)^{\frac{|\alpha|+\ell}{2}+\frac14-\frac32\frac{1}{p}}}\,\|\Delta\bfphi\|_{L^1}\,.
\end{align*}
\item There exists $\theta>0$ such that for any $\alpha\in\N^2$ and any $\ell\in\N$, there exists $C_{\alpha,\ell}$ such that for any  $2 \leq p \leq +\infty$, $1\leq q \leq p$, and any $t\geq0$
\begin{align*}
\|\,\d_\bfx^\alpha\,\d_t^\ell\,\SigOHF(t)[\bfg]\|_{L^p}
&\leq \frac{C_{\alpha,\ell}\,\eD^{-\theta\,t}}{(\min(\{1,t\}))^{\frac{|\alpha|+2\,\ell}{2}+\frac{1}{q}-\frac{1}{p}}}\,\|\bfg\|_{L^{q}}\,.
\end{align*}
\item There exists $\theta>0$ such that for any $\alpha\in\N^2$ and any $\ell\in\N$, there exists $C_{\alpha,\ell}$ such that for any  $2 \leq p \leq +\infty$, and any $t\geq0$
\begin{align*}
\|\,\d_\bfx^\alpha\,\d_t^\ell\,\SigOHF(t)[\bfphi]\|_{L^p}
&\leq \frac{C_{\alpha,\ell}\,\eD^{-\theta\,t}}{(\min(\{1,t\}))^{\left(\frac{|\alpha|+2\,\ell}{2}-\frac{1}{p}\right)_+}}\,\|\Delta\bfphi\|_{L^1}\,.
\end{align*}
\end{enumerate}
\epr

\begin{proof}
The two first sets of estimates are derived by applying the arguments of Subsection~\ref{ss:dispersive} to operators arising from the decomposition
\begin{align*}
\eD^{t\,\bD_{\bfxi}}
-\eD^{t\,\bDO(\iD\bfxi)}
&\,=\,\sum_{j=1}^2\eD^{t\,\lambda_j(\bfxi)}\,\left(\pi_j(\bfxi)-\piO_j(\bfxi)\right)\\
&\ +\,
\sum_{j=1}^2
\int_{0}^{t}\eD^{(t-s)\,\lambda_j(\bfxi)+s\,\lambO_j(\bfxi)} (\lambda_j(\bfxi)-\lambO_j(\bfxi))\,\piO_j(\bfxi)\,\dD s\,.
\end{align*}
The third and fourth sets of estimates are obtained by applying directly to $\SigOLF$ the arguments of Subsection~\ref{ss:dispersive}, whereas the last ones follow from Hausdorff-Young inequalities and arguments of Subsection \ref{ss:linear_estim_mod} through Green functions representations.
\end{proof}

\subsection{Implicit change of variables}\label{ss:changevariables}

The last transformation we would like to perform on formally derived equations is to convert equations on a $\bPsi$ or a $\bfpsi$ into equations for a $\bfphi$ related to $\bPsi$, $\bfpsi$ by 
\be\label{e:PsiPhi}
\bPsi(t,\bfx)
=(\Id+\bfpsi(t,\cdot))\left(\transp{\ubK}\,(\bfx-t\,\ubc)\right)
=(\Id-\bfphi(t,\cdot))^{-1}\left(\transp{\ubK}\,(\bfx-t\,\ubc)\right)\,,
\ee
so as to get closer to the phase introduction in the stability analysis. Our purpose is similar to the one in \cite[Appendix~B.3]{JNRZ-conservation} but we stress that for planar waves of reaction-diffusion-advection systems as considered in \cite{JNRZ-RD1,JNRZ-RD2} this discussion may easily be overlooked since systems for $\bfphi$ and $\bfpsi$ are essentially the same.

In Case~\case{b}, System~\eqref{e:W-phase-intro}
\[
\d_t\bPsiW=\bOm(\nabla\bPsiW)+\ubLambda_0[\nabla](\nabla\bPsiW)
\]
is differential and thus explicitly expressed in terms of $\d_t\bPsiW$, $\nabla \bPsiW$ and $\nabla^2\bPsiW$. Therefore the relevant algebraic manipulations stem directly from the fact when $\bPsi$, $\bfpsi$, and $\bfphi$ are related through \eqref{e:PsiPhi} one derives
\begin{align*}
\dD_\bfx\bPsi(t,\bfx)(\bfeta)
&=\transp{\ubK}\bfeta
+\dD_\bfx\bfpsi(t,\transp{\ubK}\,(\bfx-t\,\ubc))(\transp{\ubK}\bfeta)\,,\\
\nabla_\bfx\bPsi(t,\bfx)(\bfeta)
&=\ubK+\ubK\,\nabla_\bfx\bfpsi(t,\transp{\ubK}\,(\bfx-t\,\ubc))\,,\\
\dD_\bfx^2\bPsi(t,\bfx)(\bfeta,\bfzeta)
&=\dD_\bfx^2\bfpsi(t,\transp{\ubK}\,(\bfx-t\,\ubc))(\transp{\ubK}\bfeta,\transp{\ubK}\bfzeta)\,,\\
\d_t\bPsi(t,\bfx)
&=\left(\I+\dD_\bfx\bfpsi(t,\transp{\ubK}\,(\bfx-t\,\ubc))\right)(\ubOm)+\d_t\bfpsi(t,\transp{\ubK}\,(\bfx-t\,\ubc))
\end{align*}
and
\begin{align*}
\dD_\bfx\bfphi(t,\bfx+\bfpsi(t,\bfx))(\bfeta)
&=\dD_\bfx\bfpsi(t,\bfx)((\I+\dD_\bfx\bfpsi(t,\bfx))^{-1}\bfeta)\,,\\
\nabla_\bfx\bfphi(t,\bfx+\bfpsi(t,\bfx))
&=(\I+\nabla_\bfx\bfpsi(t,\bfx))^{-1}\nabla_\bfx\bfpsi(t,\bfx)\,,\\
\dD_\bfx^2\bfphi(t,\bfx+\bfpsi(t,\bfx))(\bfeta,\bfzeta)
&=(\I+\dD_\bfx\bfpsi(t,\bfx))^{-1}
\dD_\bfx^2\bfpsi(t,\bfx)((\I+\dD_\bfx\bfpsi(t,\bfx))^{-1}\bfeta,(\I+\dD_\bfx\bfpsi(t,\bfx))^{-1}\bfzeta)\,,\\
\d_t\bfphi(t,\bfx+\bfpsi(t,\bfx))
&=(\I+\dD_\bfx\bfpsi(t,\bfx))^{-1}\d_t\bfpsi(t,\bfx)\,.
\end{align*}
In Case~\case{b}, the upshot of the computations is that when $\bphiW$ is defined from $\bPsiW$, solving \eqref{e:W-phase-intro} with $\nabla_\bfx\bPsiW-\ubK$ sufficiently small, 
\begin{align*}
\d_t\bphiW
&=-\transp{\ubK}\dD_{\bK}\bfc(\ubK)(\ubK\nabla_\bfx\bphiW)
+\ubLambda_0[\ubK\nabla_\bfx](\ubK\nabla_\bfx\,\bphiW)
+\frac12\dD_{\bK}^2\bOm(\ubK)(\ubK\nabla_\bfx\bphiW,\ubK\nabla_\bfx\bphiW)\\
&\quad
+\transp{(\ubK\,\nabla_\bfx\,\bphiW)}\dD_{\bK}\bfc(\ubK)(\ubK\nabla_\bfx\bphiW)
-\transp{\ubK}\dD_{\bK}\bfc(\ubK)(\ubK\,(\nabla_\bfx\bphiW)^2)
+\brW
\end{align*}
with $\brW$ pointwise bounded by $\|\nabla_\bfx\bphiW\|\,\|\nabla_\bfx^2\bphiW\|+\|\nabla_\bfx\bphiW\|^3$, with consistent bounds for its derivatives. This suggests that it it sufficient to consider a solution to 
\begin{align}\label{eq:Wphi}
\d_t\Wbphi
&=-\transp{\ubK}\dD_{\bK}\bfc(\ubK)(\ubK\nabla_\bfx\Wbphi)
+\ubLambda_0[\ubK\nabla_\bfx](\ubK\nabla_\bfx\,\Wbphi)
+\frac12\dD_{\bK}^2\bOm(\ubK)(\ubK\nabla_\bfx\Wbphi,\ubK\nabla_\bfx\Wbphi)\\\nonumber
&\quad
+\transp{(\ubK\,\nabla_\bfx\,\Wbphi)}\dD_{\bK}\bfc(\ubK)(\ubK\nabla_\bfx\Wbphi)
-\transp{\ubK}\dD_{\bK}\bfc(\ubK)(\ubK\,(\nabla_\bfx\Wbphi)^2)
\end{align}
with the same initial data, when working on the $\bfphi$-side.

We now explain how to extend the above computations to Case~\case{a}, when System~\eqref{e:W-phase-intro} is not differential. The main issue is the lack of smoothness of involved Fourier multipliers. A convenient way to bypass this difficulty is to use as an intermediate step $\bPsiaux$ solving
\[
\d_t\bPsiaux=\bOm(\nabla\bPsiaux)
+\chi((\iD\ubK)^{-1}\nabla)\bLambda^{\ubK}[\nabla](\nabla\bPsiaux)
+(1-\chi((\iD\ubK)^{-1}\nabla))\Div((\transp{\ubK}\ubK)^{-1}\nabla\bPsiaux)
\]
with the same initial data as $\bPsiW$. The equation for $\bPsiaux$ is still non local, because of frequency cut-off operators, but non locality is encoded by smooth multipliers. In order to identify the leading-order part of the equation for $\bphiaux$ defined from $\bPsiaux$, the only missing argument is contained in the following lemma.

\bl
Let $\kappa\in[0,1)$, $s\in\N$ and $\tchi$ smooth and compactly supported. There exists $C_{\kappa,s,\tchi}$ such that for any $\bfpsi$ such that $\|\nabla\bfpsi\|_{L^\infty}\leq \kappa$, for any $\bF$, any $|\alpha|\leq s$, and any $1\leq p\leq \infty$,
\begin{align*}
\big\|\,\d^\alpha\big[(\tchi(\iD^{-1}\nabla)\,[\bF\circ(\Id-\bfpsi)^{-1})]\circ(\Id-\bfpsi)
&-\,\tchi(\iD^{-1}\nabla)\,\bF\big]\,\big\|_{L^p}\\
&\leq C_{\kappa,s,\tchi}
\sum_{\ell=0}^{|\alpha|}\,\left\|\nabla^\ell\bF\right\|_{L^p}
\,\left\|\nabla^{1+|\alpha|-\ell}\bfpsi\right\|_{L^\infty}\,.
\end{align*}
\el

\begin{proof}
Note that 
\begin{align*}
(\tchi(\iD^{-1}\nabla)\,[\bF\circ(\Id-\bfpsi)^{-1})](\bfx-\bfpsi(\bfx))
-\,\tchi(\iD^{-1}\nabla)\,[\bF](\bfx)
\,=\,\int_{\R^2}\,\Gamma(\bfx,\bfy)\,\bF(\bfy)\,\dD\,\bfy
\end{align*}
with
\begin{align*}
\Gamma(\bfx,\bfy)
&=\,\int_{\R^2}\,\eD^{\iD\,\bfxi\cdot (\bfx-\bfy)}
\,\left(\,\frac{\det(\I-\nabla\bfpsi(\bfy))}{\det(\I-\bM_{\bfpsi}(\bfx,\bfy))}
\tchi\left((\I-\bM_{\bfpsi}(\bfx,\bfy))\bfxi\right)
-\tchi(\bfxi)\,\right)\,\dD\bfxi
\end{align*}
where
\[
\bM_{\bfpsi}(\bfx,\bfy)
\,=\,\int_0^1\,\nabla\bfpsi(\bfy+\tau\,(\bfx-\bfy))\,\dD \tau\,.
\]
From here, simplest techniques used throughout the text to provide pointwise bounds achieve the proof of the lemma.
\end{proof}

By combining the lemma with algebraic computations expounded above one derives that as long as $\nabla_\bfx\bPsiaux-\ubK$ remains sufficiently small, 
\begin{align*}
\d_t\bphiaux
&=-\transp{\ubK}\dD_{\bK}\bfc(\ubK)(\ubK\nabla_\bfx\bphiaux)
+\chi(\iD^{-1}\nabla)\bLambda^{\ubK}[\ubK\nabla_\bfx](\ubK\nabla_\bfx\,\bphiaux)
+(1-\chi(\iD^{-1}\nabla))\Delta\bphiaux\\
&\quad
+\frac12\dD_{\bK}^2\bOm(\ubK)(\ubK\nabla_\bfx\bphiaux,\ubK\nabla_\bfx\bphiaux)\\
&\quad
+\transp{(\ubK\,\nabla_\bfx\,\bphiaux)}\dD_{\bK}\bfc(\ubK)(\ubK\nabla_\bfx\bphiaux)
-\transp{\ubK}\dD_{\bK}\bfc(\ubK)(\ubK\,(\nabla_\bfx\bphiaux)^2)
+\braux
\end{align*}
with $\braux$ bounded in $L^p$ by $\|\nabla_\bfx\bphiaux\|_{L^\infty}\,\|\nabla_\bfx^2\bphiaux\|_{L^p}+\|\|\nabla_\bfx\bphiaux\|^3\|_{L^p}$, with consistent bounds for its derivatives. As a final intermediate step, $\bphiaux$ may be compared to $\auxbphi$ solving the same equation without the remainder $\braux$, and starting from the same initial data. In turn, $\auxbphi$ is readily compared to $\Wbphi$ solution of \eqref{eq:Wphi}.

The final result is as follows.

\bpr\label{p:Wphi} There exist $\eps_0>0$, $(C_{(p_0)})_{p_0>2}$ and $(C_{\ell})_{\ell\in\N,\,\ell\geq1}$ such that if for some sublinear\footnote{As in Theorem \ref{th:nonlinear_stab} we mean that $\bfphi_0$ may differ from $\Delta^{-1}(\Delta \bfphi_0)$ by a constant function but not by a non-constant affine function.} $\bfphi_0$
\[
\|\Delta \bfphi_{0}\|_{(L^{1} \cap L^{4})(\R^2;\R^2)}\,\leq\,\eps_0
\]
then, there exist a unique global solution $\Wbphi$ to \eqref{eq:Wphi} with initial datum $\Wbphi(0,\cdot) = \bfphi_{0}$ and a unique global solution $\bPsiW$ to \eqref{e:W-phase-intro} with initial datum $\bPsiW(0,\cdot)= \bPsi_{0}$ given by 
\[
\bPsi_0(\bfx)
=(\Id-\bfphi_0)^{-1}\left(\transp{\ubK} \bfx \right)\,,
\]
such that, with
\[
r_{p,s}:=\begin{cases}
\frac34-\frac32\,\frac1p+\frac{s}{2}&\quad\textrm{ in Case~\case{a}}\\
\frac12-\frac1p+\frac{s}{2}&\quad\textrm{ in Case~\case{b}}
\end{cases}
\]
there holds for any $t\geq0$,
\begin{align*}
\|\d^\alpha\Wbphi \|_{L^4}
&\leq \frac{C_{|\alpha|}\,\|\Delta \bfphi_{0}\|_{L^{1} \cap W^{|\alpha|-2,4}}}{(1+t)^{r_{p,|\alpha|-1}}}\,,&\quad |\alpha|\geq 2\,,\\
\|\d^\alpha\Wbphi \|_{L^p}
&\leq \frac{C_{|\alpha|+1}\,\|\Delta \bfphi_{0}\|_{L^{1} \cap W^{|\alpha|-1,4}}}{(1+t)^{r_{p,|\alpha|-1}}}\,,&\quad |\alpha|\geq 2\,,\qquad 2\leq p\leq \infty\,,\\
\|\nabla\Wbphi \|_{L^p}
&\leq \frac{C_{(p_0)}\,\|\Delta \bfphi_{0}\|_{L^{1} \cap L^{4}}}{(1+t)^{r_{p,0}}}\,,&\quad 2<p_0\leq p\leq \infty\,,\\
\end{align*}
and, with $\bphiW$ defined from $\bPsiW$ and $\bPsi_W$ from $\Wbphi$ through \eqref{e:PsiPhi}, for any $t\geq0$,
\begin{align*}
\|\d^\alpha(\Wbphi-\bphiW) \|_{L^p}
&\leq \frac{C_{|\alpha|+2}\,\|\Delta \bfphi_{0}\|_{L^{1} \cap W^{|\alpha|,4}}}{(1+t)^{r_{p,|\alpha|}+\frac12}}\,,&\quad |\alpha|\geq 0\,,\qquad 2\leq p\leq \infty\,,\\
\|\d^\alpha(\bPsi_W-\bPsiW) \|_{L^p}
&\leq \frac{C_{|\alpha|+2}\,\|\Delta \bfphi_{0}\|_{L^{1} \cap W^{|\alpha|,4}}}{(1+t)^{r_{p,|\alpha|}+\frac12}}\,,&\quad |\alpha|\geq 0\,,\qquad 2\leq p\leq \infty\,.
\end{align*}
\epr

We skip the proof of Proposition~\ref{p:Wphi} because we have already discussed the main ingredients of the proof, and the remaining parts are very similar to arguments detailed elsewhere along the text. Let us simply point out that the existence part stems from a simple fixed-point argument, since equations are essentially of semilinear parabolic type. Note that the last estimate follows from Lemma \ref{l:var-change}.

\newpage

\bibliographystyle{alphaabbr}
\bibliography{Ref_RD2D}

\end{document}